\documentclass{gtart}
\usepackage{amssymb,rlepsf}
\let\Bbb\mathbb

\input labelfig

\input gtoutput
\volumenumber{2}\papernumber{10}\volumeyear{1998}
\pagenumbers{221}{332}\published{6 January 1999}
\proposed{Robion Kirby}
\seconded{Gang Tian, Tomasz Mrowka}
\received{2 February 1998}\revised{20 November 1998}
\accepted{3 January 1999}

\let\ss\scriptstyle
\let\relabela\adjustrelabel
\let\med\medskip

\renewcommand{\theequation}{\thesection.\arabic{equation}}

\newcommand{\x}{\times}

\renewcommand{\a}{\alpha}
\renewcommand{\b}{\beta}
\renewcommand{\d}{\delta}
\newcommand{\dt}{\cdot}
\newcommand{\D}{\Delta}
\newcommand{\e}{\varepsilon}
\newcommand{\g}{\gamma}
\newcommand{\G}{\Gamma}
\renewcommand{\k}{\kappa}
\renewcommand{\l}{\lambda}
\renewcommand{\L}{\Lambda}

\newcommand{\R}{{\mathbb R}}
\newcommand{\s}{\sigma}
\newcommand{\Sig}{\Sigma}

\renewcommand{\th}{\theta}
\renewcommand{\O}{\Omega}
\renewcommand{\o}{\omega}
\newcommand{\z}{\zeta}
\renewcommand{\i}{\infty}
\newcommand{\p}{\partial}
\def\itemnum#1{\hspace\fill{\rm (#1)} \addtocounter{equation}{1}}

\begin{document}

\title{The structure of pseudo-holomorphic subvarieties\\\vspace{-1.4mm}\\
for a degenerate almost complex structure\\
and symplectic form on $S^1 \times B^3$}
\covertitle{The structure of pseudo-holomorphic subvarieties\\
for a degenerate almost complex structure\\
and symplectic form on $S^1 \times B^3$}
\asciititle{The structure of pseudo-holomorphic subvarieties\\
for a degenerate almost complex structure\\
and symplectic form on S^1 X B^3}
\shorttitle{The structure of pseudo-holomorphic subvarieties}

\author{Clifford Henry Taubes}

\address{Department of Mathematics\\
Harvard University\\
Cambridge, MA 02138, USA}

\email{chtaubes@abel.math.harvard.edu}

\begin{abstract}
A self-dual harmonic 2--form on a 4--dimensional
Riemannian manifold is symplectic where it does not vanish. Furthermore,
away from the form's zero set, the metric and the 2--form give a
compatible almost complex structure and thus pseudo-holomorphic
subvarieties. Such a subvariety is said to have finite energy when the
integral over the variety of the given self-dual 2--form is finite.
This article proves a regularity theorem for such finite energy
subvarieties when the metric is particularly simple near the form's zero
set. To be more precise, this article's main result asserts the
following: Assume that the zero set of the form is non-degenerate and
that the metric near the zero set has a certain canonical form.
Then, except possibly for a finite set of points on the zero set,
each point on the zero set has a ball neighborhood which intersects the
subvariety as a finite set of components, and the closure of each
component is a real analytically embedded half disk whose boundary
coincides with the zero set of the form.
\end{abstract}

\asciiabstract{%
A self-dual harmonic 2-form on a 4-dimensional
Riemannian manifold is symplectic where it does not vanish. Furthermore,
away from the form's zero set, the metric with the 2-form give a
compatible almost complex structure and thus pseudo-holomorphic
subvarieties. Such a subvariety is said to have finite energy when the
integral over the variety of the given self-dual 2-form is finite.
This article proves a regularity theorem for such finite energy
subvarieties when the metric is particularly simple near the form's zero
set. To be more precise, this article's main result asserts the
following: Assume that the zero set of the form is non-degenerate and
that the metric near the zero set has a certain canonical form.
Then, except possibly for a finite set of points on the zero set,
each point on the zero set has a ball neighborhood which intersects the
subvariety as a finite set of components, and the closure of each
component is a real analytically embedded half disk whose boundary
coincides with the zero set of the form.}


\maketitlepage

\section{Introduction}
\setcounter{equation}{0}

The purpose of this paper is to prove a regularity theorem for
pseudo-holomor\-phic curves on the 4--manifold $X = S^1\x  B^3$ with a
 degenerate
symplectic structure. Here, $S^1$ is the circle and $B^3$ is the standard
3--ball
in ${\R}^3$. Let $t$ be the standard coordinate on the circle, and
let $(x,
y, z)$ be Euclidean coordinates on ${\R}^3$. Consider now the 2--form:
\begin{equation}
\o=dt\wedge    (x dx+y dy-2z dz)+ 
x dy\wedge     dz-y dx\wedge dz-2z dx\wedge dy
\end{equation}
One can check that this form is closed; moreover, $\o$ is self-dual with
respect to the obvious product flat metric on $S^1\x  B^3$. Thus, $\o$
 is
non-degenerate where $\o$ is non-zero. Indeed,
\begin{equation}
\o\wedge\o=2\dt(x^2+y^2+4\dt z^2)\dt dt\wedge dx\wedge dy\wedge dz
\end{equation}
and so $\o$ vanishes only on the circle $Z = S^1\x\{0\}$.

On the complement of the circle $Z$, the form $\o$ is non-degenerate and so
symplectic. On this complement, $\o$ and the flat, product metric
define an
almost complex structure, $J$, which is given as follows: First,
introduce
\begin{equation}
g=(x^2+y^2+4\dt z^2)^{\frac 12}
\end{equation}
which is the norm of $\omega$.

Then, set:
\begin{itemize} %1.4
 \item\quad $J\dt\p_t\ =-g^{-1}(x\dt\p_x+y\dt\p_y-2\dt
z\dt\p_z)$
 \item\quad $J\dt\p_x\ =g^{-1}(x\dt\p_t+y\dt\p_z+2\dt
z\dt\p_y)$
 \item\quad $J\dt\p_y\ =g^{-1}(y\dt\p_t-x\dt\p_z-2\dt
z\dt\p_x)$
 \item\quad $J\dt\p_z\ =g^{-1}(-2\dt
z\dt\p_t+x\dt\p_y-y\dt\p_x)$. \itemnum{1.4}
\end{itemize}

By the way, note that $J$ is $\o$--compatible as $\sqrt{2} \o( \dt  ,
J(\dt))/|\o |$ gives back
the product, flat metric. Note as well that $J$ is singular on $Z$.

This example of the self-dual, closed form $\o$ on the 4--manifold
$S^1\x  B^3$
serves as a model for singular symplectic forms defined by closed, self-dual
forms on general 4--manifolds. The relevance of the $(S^1\x  B^3,\o)$ model in
this more general context is remarked upon further in subsection 1.c, below.

Now consider:

\medskip
{\bf Definition 1.1}\qua Let $Y$ be a smooth 4--manifold (in this
case,
$S^1\x  (B^3-\{0\}$) with symplectic form $\o$ and compatible almost
complex structure
$J$. A subset $C\subset  Y$ will be called a {\sl pseudo-holomorphic
subvariety}  if the
following is true:
\begin{itemize}
\item  $C$ is closed.

\item There exists a $C^\i$  complex curve $C_0$ together
 with
a proper, $J$--pseudo-holomorphic map $f\co  C_0 \to  Y$ whose image is $C$
and which is an embedding on the complement of a countable set of points in
$C_0$. The complex curve $C_0$ will be called the smooth model
for $C$.
\end{itemize}
Call $C$ a {\sl pseudo-holomorphic submanifold} when $C$ is non-singular.
\med

The next definition introduces the notion of finite energy for a
pseudo-holomor\-phic subvariety.

\medskip
{\bf Definition 1.2}\qua Let $C$ be a pseudo-holomorphic subvariety as
in Definition 1.1. Say that $C$ has {\sl finite energy} when the
integral over $C$
of $\o$ is finite.  (This is equivalent to the assertion that the
integral over
the smooth model $C_0$ of $f^*\o$ is finite.)
\med

(Note that $\o|_{TC} > 0$ so that there is no question here with
conditional
convergence of the integral.)

The main theorem in this article (Theorem 1.8, below) describes the nature
of a finite energy, pseudo-holomophic subvariety near $Z$. However, part of
the story is simply this:

\medskip
{\sl Fix a finite energy, pseudo-holomorphic subvariety $C$. Then, all but
finitely many points in $Z$ have a neighborhood in $S^1\x
B^3$ in
which $C$ has finitely many components. Moreover, each of these
components is a smoothly embedded half-disk (that is, the image of
$\{w\in{\Bbb C}: |w| < 1$ and im$(w) \geq 0\}$).
Furthermore, the boundary line (where im$(w) = 0$) of each such
half-disk lies on $Z$.}

\subsection*{1.a\qua Examples}

Before turning to the precise statement of the main theorem, some examples
are in order.

The first examples below are the simplest cases of finite energy,
pseudo-holomor\-phic submanifolds whose closure in $S^1\x  B^3$ is a submanifold
with boundary where the latter coincides with $Z$.

\medskip
{\bf Example 1.3}\qua Fix a point $\nu$ on the boundary of unit disk in
the $(x, y)$
plane (where $z = 0$) in ${\R}^3$. Let $\l_v$ denote the complement of the
origin in the line segment between the origin in ${\R}^3$ and $\nu$. Then
\begin{equation} %1.5
C=S^1\x \l_v
\end{equation}
is a smooth, pseudo-holomorphic subvariety. See Figure 1.

\begin{figure}[htb!]
\centerline{%\ShowGrid 
\SetLabels\small
(.53*0.37){$C$}\\
(.37*0.42){$t$}\\
(.08*0.01){$Z$}\\
(.58*0.54){$y$}\\
(.8*0.8){$x=y=z=0$}\\
(.46*0.97){$z$--axis}\\
\endSetLabels 
\AffixLabels{\epsfxsize 2.5in \epsfbox{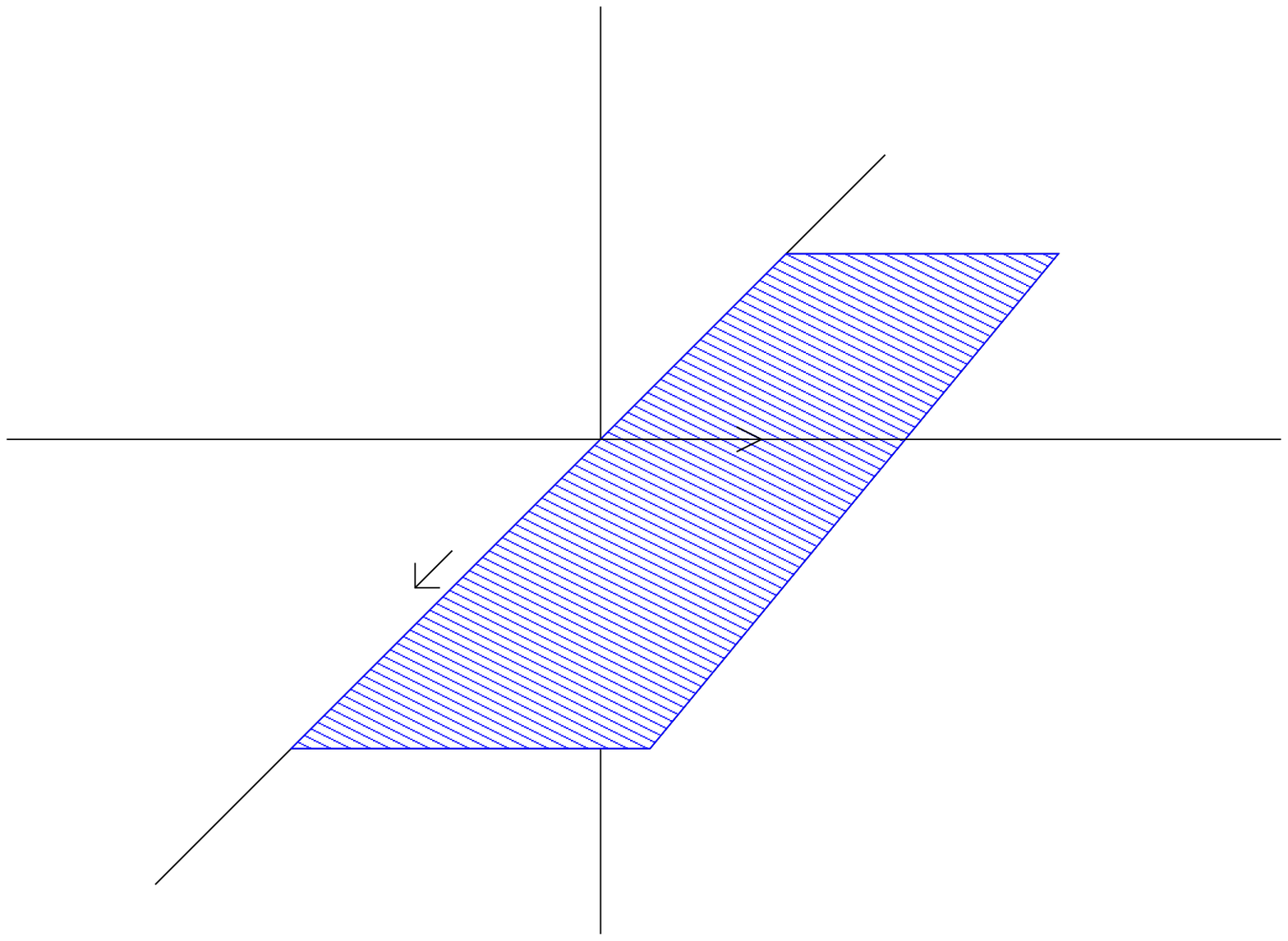}}}
\nocolon\caption{}
\end{figure}

\med{\bf Example 1.4}\qua Let $\nu = (0, 0 ,\pm 1) \in\R^3$ and define
$\l_\nu$ as before, and then define $C$ as in (1.5).  See Figure 2.

\begin{figure}[htb!]
\centerline{%\ShowGrid 
\SetLabels\small
(.34*0.47){$C$}\\
(.18*0.12){$t$}\\
(.13*-.02){$Z$}\\
(.78*0.32){$xy$--plane}\\
(.88*0.83){$x=y=z=0$}\\
(.46*0.93){$z$--axis}\\
\endSetLabels 
\AffixLabels{\epsfxsize 2.5in \epsfbox{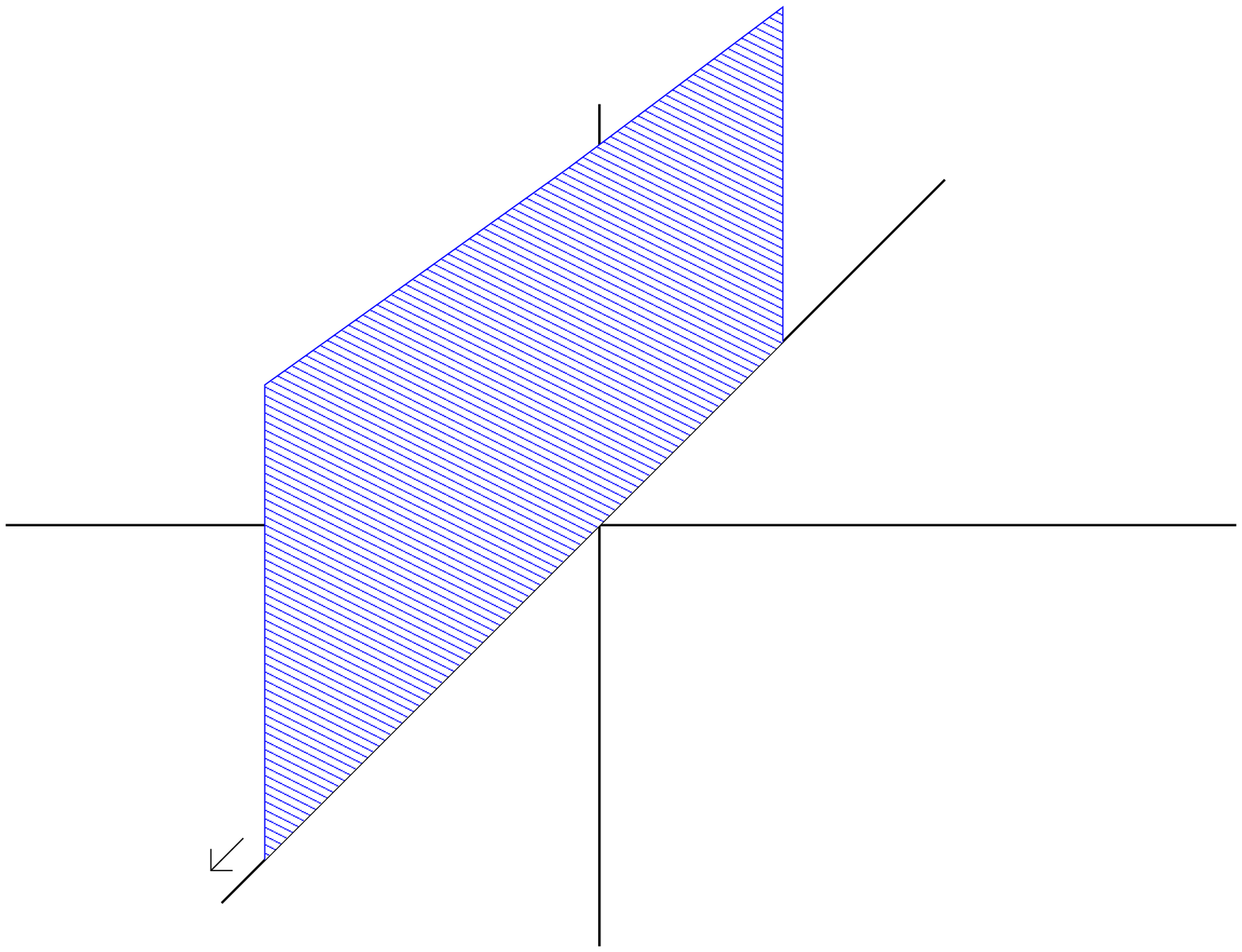}}}
\nocolon\caption{}
\end{figure}

Note that in both of these examples, $C$ is parametrized in the
obvious way
by points $(t, s) \in S^1\x  (0, 1]$, and the pull-back to $C$ of the form
 $\o$ is
$dt  \wedge s\dt ds$ when $\eta$ lies in the $(x, y)$ plane, and the
pull-back is
$2\dt dt\wedge  s\dt ds$ when $\eta$ is on the $z$--axis. Thus, in
either case, $C$ has
finite energy. Indeed, the integral of $\o$ over $C$'s intersection
 with the product of
$S^1$ and the ball of radius $r$ and center 0 in $\R^3$ is equal to
$r^2/2$ in the
first case, and $r^2$ in the second. (The similarities between the
cases in Example 1.3 and those from Example 1.4 are superficial, as
can be seen in Theorem 1.8 below.)

\med
The next examples are cones. Each has closure which intersects $Z$ in a single
point. For these examples, let $(\rho,\phi)$ denote the standard polar
coordinates
in the $(x, y)$ plane so that $x = \rho\dt\cos\phi$ and $y =
\rho\dt\sin\phi$.

\med
{\bf Example 1.5}\qua Fix an origin $t_0 \in S^1$. Let $c \in (0, 2/3^{3/2}]$
be a
constant and let $u(\dt)$ be a positive function on $\R$ which
satisfies the differential equation
\begin{equation}%1.6
u'' + \textstyle{\frac 49} \ u - \textstyle{\frac 29} \ cu^{-2}=0
\end{equation}
with the normalization $\frac{9}{4} (u' )^2 + u^2 + c/u = 1$.
(Here, $'$  denotes differentiation.) Note that all such
functions are periodic and there exist solutions with period in $2\pi\Bbb Q$.
Let $u$ be one of the latter and write its period as $2\pi a/b$
where $a$ and $b$ are positive integers which are relatively prime.
 Introduce
the pair of pseudo-holomorphic submanifolds $\Sigma_{\pm}\subset S^1\x B^3$ by
using $s\in (0, \i)$ and $\tau \in [0, 2\pi b]$ as coordinates and
using the parametrization rule:

\begin{itemize}
\item \quad$ t=t_0\pm s \textstyle\frac 32 u'(\tau)$
\item  \quad$z=\pm su(\tau)$
\item  \quad$\rho = sc^{\frac 12} |u(\tau)|^{-\frac 12}$
\item \quad$ \phi=\tau$ \hspace\fill (1.7)
\end{itemize}
In the extreme case where $c = 2/3^{3/2}$, the solution to (1.6) is $u\equiv
 1/3^{1/2}$. The corresponding submanifold from (1.7) lies in the time
slice $t
= t_0$ and coincides with either the $z > 0$ (for $\Sigma_+$) or the $z < 0$
(for $\Sigma_-$) half of the cone $x^2 + y^2 - 2\dt z^2 = 0$ in
$B^3$. See Figure 3.

\begin{figure}[htb!]
\centerline{\relabelbox\small
\epsfxsize 2.5in \epsfbox{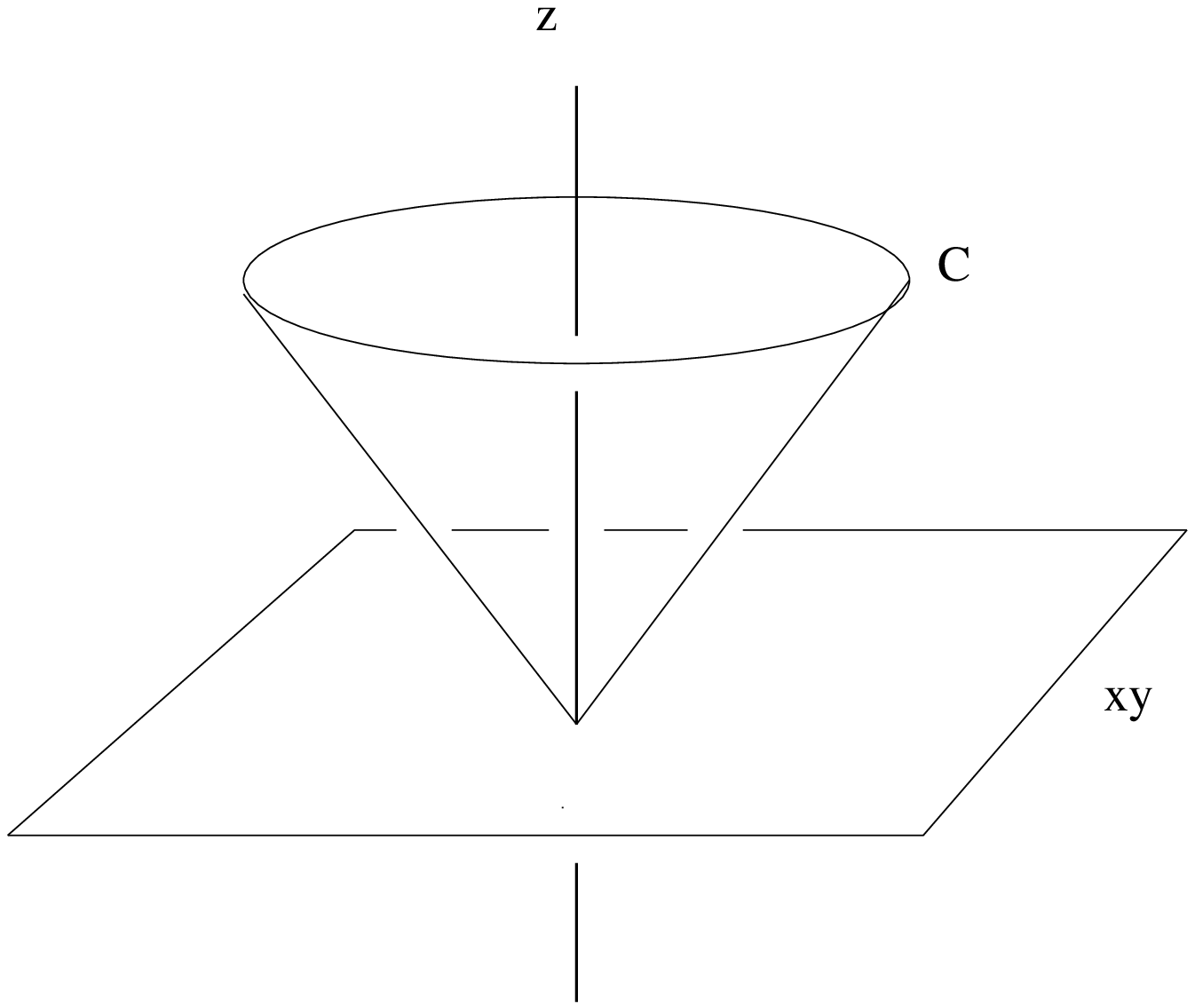}
\relabel {z}{$z$--axis}
\relabel {C}{$C$}
\relabel {xy}{$xy$--plane}
\endrelabelbox\hbox to 10mm{\hfil}}
\nocolon\caption{}
\end{figure}

\med
(Note that the cones in (1.7) are embedded because of the assumption that $a$
and $b$ are relatively prime.) The fact that $\Sig_{\pm}$  in (1.7) is
pseudo-holomorphic can be checked directly, or the reader can refer to an
argument which is given in the proof of Proposition 3.2 in Section 3. The
other assertions about $\Sig_{\pm}$  and the assertions about the
solutions to
(1.6) are proved as a part of Lemma 1.9 in subection 1.d, below.

These last examples also have finite energy. Indeed, the ball $B_r \subset
S^1\x
B^3$ of radius $r$ about the point $(t_0, 0)$ contains the intersection of
$\Sig_{\pm}$
with a neighborhood of $Z$; meanwhile, the fact that $\Sig_{\pm}$ is a cone
while $\o$ has homogeneous dependence on the coordinates implies that the
integral of $\o$ over $\Sig_{\pm}\cap B_r$ is equal to $r^3 \d_c$.
 Here, $\d_c$
is a constant which depends on $c$. Note however that there is a
$c$--independent
constant $\z\geq 1$
 with the property that $\d_c \geq \z^{-1}$. (This last assertion is
restated as Lemma 1.10 in subsection 1.d.)

For the next example, introduce the following pair of functions on $B^3$:
\begin{itemize} %1.8
\item\quad $f= 2^{-1}\dt (x^2+y^2-2\dt z^2)$
\item\quad  $h= z\dt (x^2+y^2)$ \hspace\fill (1.8)
\end{itemize} 
The next example includes smooth pseudo-holomorphic subvarieties which miss
an entire neighborhood of $Z$.

\med
{\bf Example 1.6}\qua Let $c$ and $c'$  be constants and let $C$ denote
the submanifold in $S^1 \x  (B - 0)$ which is described by one of the
following.  See Figure 4.
\begin{itemize} %1.9
\item\quad $\phi = c \ \ {\mbox{and}} \ \ h=c' $ 
\item\quad $ t =c \ \ {\mbox{and}} \ \ f=c'>0$ 
\item\quad  $t = c \ \ {\mbox{and}} \ \ f=c'\leq 0
\ \ {\mbox{and}} \ \ \pm z >0 $  \hspace\fill (1.9)
\end{itemize}

\begin{figure}[t!]
\centerline{\relabelbox\small
\epsfxsize 2in \epsfbox{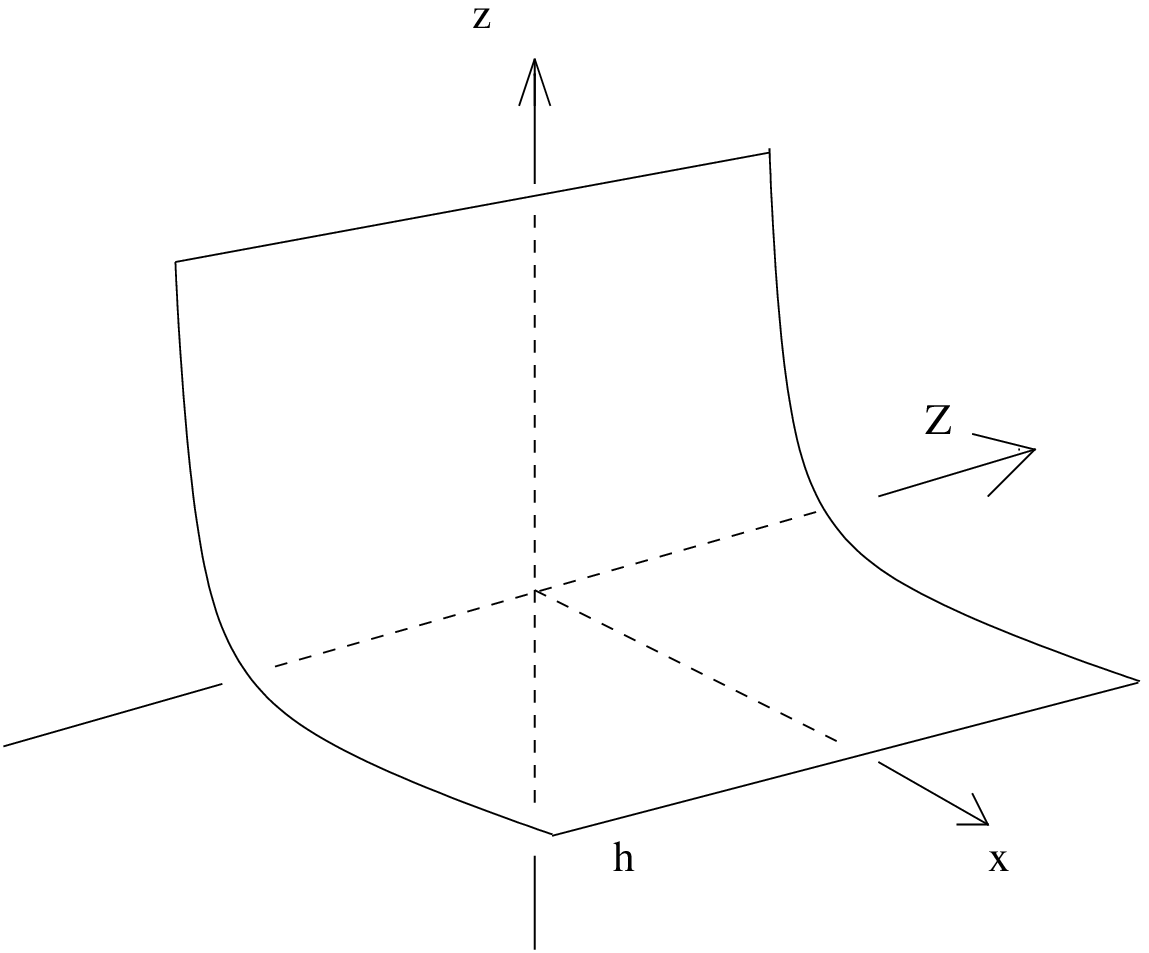}
\relabel {z}{$z$--axis}
\relabel {x}{$x$--axis}
\relabela <-3pt, 0pt> {Z}{$Z$}
\relabela <-7pt, -2pt> {h}{$\ss h=c'>0$}
\endrelabelbox\qquad 
\relabelbox\small
\epsfxsize 2in \epsfbox{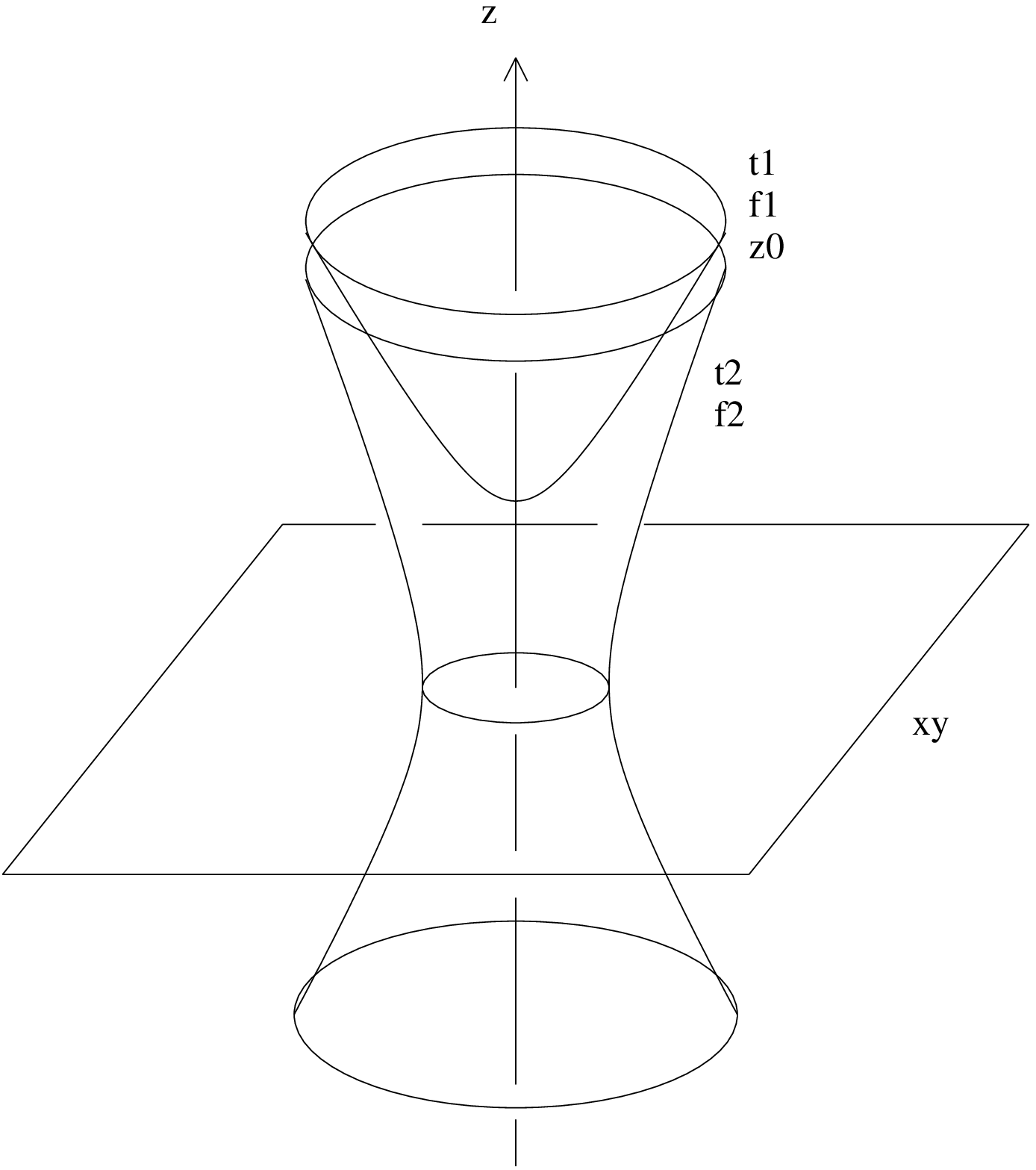}
\relabel {xy}{$xy$--plane}
\relabela <-3pt, -0pt> {z}{$z$--axis}
\relabela <0pt, 2pt> {f1}{$\ss f=c'<0$}
\relabela <0pt, -2pt> {f2}{$\ss f=c'>0$}
\relabel {t2}{$\ss t=c$}
\relabela <0pt, 4pt> {t1}{$\ss t=c$}
\relabel {z0}{$\ss z>0$}
\endrelabelbox}
\nocolon\caption{}
\end{figure}

\med
Note that the submanifolds in the first line of (1.9) whose closure in
$X$ intersect $Z$ are precisely the $C$'s from Examples 1.3 and
1.4. Likewise, those in the second line of (1.9) whose closure in $X$
intersect $Z$ are precisely the $C$'s from the extreme case 
$c = 2/3^{3/2}$ in Example 1.5.

By the way, as the constants $c$ and $c'$ vary in (1.9) and Figure 4,
the pseudoholomorphic subvarieties provide a pair of mutually
orthogonal, 2--dimensional foliations of $S^1 \times (B^3 - 0)$ by
finite energy, pseudoholomorphic subvarieties.  Indeed, the existence
of such a pair of foliations is related to a very useful fact: 
the form $\omega$ and the flat metric are particularly nice
in the
$(t,f, h, \phi)$ coordinates. Indeed,
\addtocounter{equation}{3}
\begin{equation} %1.10
\o=dt\wedge df+d\phi\wedge dh ,
\end{equation}
while the flat metric is
\begin{equation} % 1.11
dt^2+g^{-2}\dt df^2+g^{-2}\dt\rho^{-2}\dt dh^2+\rho^2\dt d\phi^2.
\end{equation}
(The fact that (1.9) are pseudo-holomorphic follows immediately from (1.10)
and (1.11).)

The final example, below, gives additional cases of a pseudo-holomorphic
subvariety whose closure in some neighborhood of each point in $Z$ is a union
of embedded half-disks. Note that these examples are connected, but even so,
some do not have submanifold with boundary closures in any open set which
contains all of $Z$.

\med
{\bf Example 1.7}\qua Let $q >0$ and $p$ be relatively prime integers and
let $\alpha \in S^1$. Let $C$ denote the image of the map from $S^1\x (0,
1]$ into
$S^1\x (B - 0)$ which sends $(\phi, s)$ to
\begin{equation}  %1.12
(q\dt\phi, s\dt\cos(p\dt\phi+\a),s\dt\sin(p\dt\phi+\a),z(s)).
\end{equation}

\begin{figure}[htb!]
\centerline{%\ShowGrid 
\SetLabels\small
(.1*1.01){$C$}\\
(1*0){$t$}\\
(.0*.73){$Z$}\\
(.8*0.72){$x$}\\
(.91*0.55){$y$}\\
\endSetLabels 
\AffixLabels{\epsfxsize 3in \epsfbox{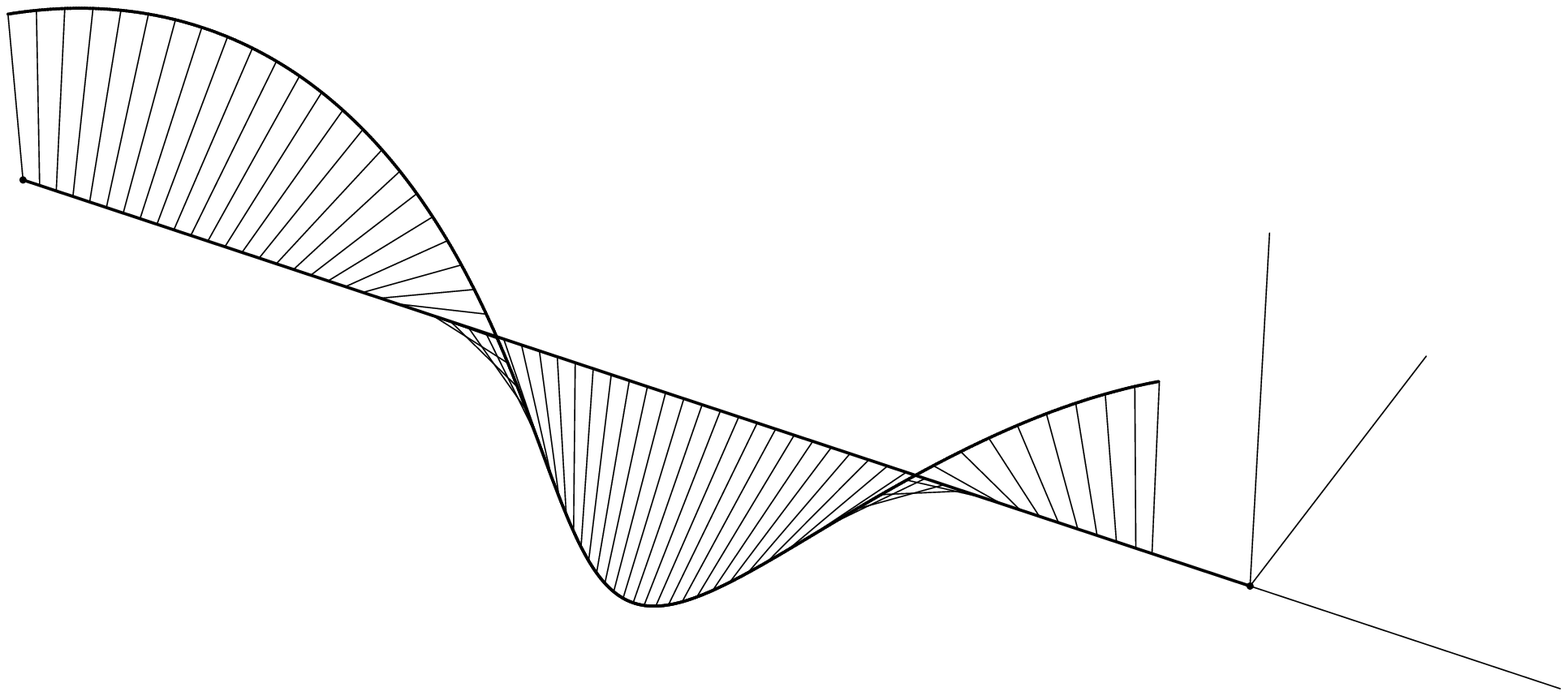}}}
\nocolon\caption{}
\end{figure}

See Figure 5. Here, $z(\dt )$ is the smooth function on $(0, 1]$ which obeys:
\begin{itemize} %1.13
\item\quad  $(q+2p\dt z)\dt s\dt z'+2q\dt z=p\dt s^2$
\item\quad $z(s)=\frac{p}{4q}\dt s^2+{\cal O}(s^4)$
 \hspace\fill (1.13)
\end{itemize}
Note that in the case $p = 0$, this $C$ is the same as in
Example 1.3.

\subsection*{1.b\qua The main theorem}

The following theorem summarizes the principle results of this article:

\med
{\bf Theorem 1.8}\qua{\sl{Let} $C \subset S^1 \x (B^3 =\{0\}$) {be a finite
energy, pseudo-holomorphic subvariety. Then $C$ has finite area and each
 point in
its closure in} $S^1 \x  B^3$ {has well defined tangent cones
up to the
rotation} $\phi\to \phi + {\rm constant}$. {Furthermore, there is a finite set}
$Z_s\subset  Z = S^1\x\{0\}$ {and a positive integer $N$ with
the following
significance:  Let} $\tau
\in Z - Z_s$. {Then $\tau$ centers an open ball in $S^1 \x
B^3$ whose
intersection with $C$ is a disjoint union of $N$ embedded components.
Furthermore, the closure of each such component is the image of the
standard
half disk via a real analytic embedding which sends the straight edge to
$Z$ and
the interior to the complement of $Z$. Finally, this embedding is described by
one of the following two possibilities:}

$\bullet$\qua{The image is the half disk} $\{(t, x, y,
z): (t - \tau)^2 +
z^2\leq \e^2$, $x = y = 0$ {and either $z\geq  0$ or $z\leq 0\}$.
 This
is a half disk in the $x = y = 0$ plane and a part of the submanifold in
Example
1.4.}

$\bullet$\qua{The image is that of the half disk} $\{(t,
\rho) \in\R^2:
 (t - \tau)^2 + \rho^2 \leq\e^2$ {and} $\rho\geq 0\}$. {The
embedding of
this half disk is defined by functions} 
$(\phi=\phi(t, \rho), \  z = z(t,
\rho)\}$ {via the parametrization} $(x = \rho\cos\phi, \ y = \rho
\sin\phi,
z)$. {Here,} $(\rho, z)$ {are real analytic functions of $(t, \rho)$
which obey} $\phi(t, \rho) = \phi_0(t) + {\cal O}(\rho^2)$ {and}
$z(t,\rho) =
4^{-1}\phi'_0  (t)\dt\rho^2 + {\cal O}(\rho^4)$ {near $\rho = 0$. In
particular, this portion of $C$ is, locally, a perturbation of an Example 1.3
model. (See also Example 1.7.)}}\med

The set $Z_s \subset Z$ will be called the singular set. The behavior of
$C$ near
$Z_s$ can also be described in some detail. The simplest sort of
singularity is
given in Example 1.5, and multi-sheeted versions of Example 1.5 can also
occur.
However, there is also a second, less concentrated type of singularity which
can be viewed as a perturbation of the $h =$ constant ($\neq 0$) and
 $\phi =$
constant example from Example 1.6. In the latter, $\phi$ and $h$
become (slowly
varying) functions of $t$ and $f$ where $h$ vanishes at a point where $f =
0$ and
$t = \tau\in Z_s$.
Some of the structure of $C$ near $Z_s$ is elucidated in Proposition
4.1,
below. In any event, the detailed structure of $C$ near $Z_s$ and the
question of
eliminating $Z_s$ by perturbations will be the subject of a planned
sequel to
this article.

As a parenthetical remark, note that Hofer, Wysocki and Zehnder in \cite{HWZ1},
\cite{HWZ2} have studied (using different techniques) pseudo-holomorphic
curves on
$M^3 \x\R$ in the case when the symplectic form arises from a contact
structure on $M^3$. See also \cite{HK} for a more recent application of this
approach.

The rest of this article is more or less centered around the proof of
Theorem 1.8, although various constructions are presented with regard to
pseudo-holomorphic subvarieties in $S^1 \x(B^3 - \{0\}$) which may have
some broader utility. Aside from subsections 1.c and 1.d, below, there are
Sections 2--9 plus a lengthy appendix. With regard to the proof of Theorem
1.8, the reader should note that the assertions of the theorem are either
reasserted as parts of propositions below or are their direct corollaries.
In particular, the assertion that $C$ has finite area is part of
Proposition
3.2. The fact that the tangent cones are unique up to rotation follows from
Proposition 4.1.  The assertion that there are at most finitely many
singular points on $Z$ follows from Propositions 5.1 and 9.1. Finally, the
description of $C$ near a regular point in $Z$ comes from Proposition
7.1.

In any event, here is a table of contents for this article:

Section 1\qua {\sl Introduction:}\qua  
Definitions, examples and the main theorem

Section 2\qua {\sl Finite energy and finite local energy:}\qua  
The derivation of a key monotonicity formula for the local energy in a
ball centered on $Z$

Section 3\qua {\sl Limits:}\qua 
The ultra-local structure near $Z$ as studied via
successive dilations and their limits

Section 4\qua {\sl The multiplicity of the limit:}\qua A uniqueness
theorem for dilation limits

Section 5\qua {\sl The integers $n_{\pm}$:}\qua  
A proof that there are only finitely
many singularities like Example 1.5

Section 6\qua {\sl The semi-continuity of $p$ and $q$ as functions on $Z$:}\qua
How local sheet number can change

Section 7\qua  {\sl Regular points:}\qua  The discussion of the first bullet in
Theorem 1.8

Section 8\qua {\sl A proof that $C_{\rho}$ is smooth:}\qua  
The discussion of the second bullet in Theorem 1.8

Section 9\qua {\sl The number of singular points:}\qua  
The proof that there are only finitely many

{\sl Appendix:}\qua Proofs of some assertions about multi-sheeted
examples

Before proceeding to the subsequent subsections of Section~1 or to Section~2,
please be advised of three conventions that are used in the subsequent
sections of this article. First, many equations involve constants whose
precise value is not terribly relevant to the discussion. Such constants are
typically denoted below by $\z$. Thus, the value of $\z$ will
 change from appearance to appearance.

The second convention involves metrics on a pseudo-holomorphic subvariety $C$:
The convention here is that the metric at the smooth points of $C$ is that
which is induced by $C$'s embedding into $S^1 \x  B^3$.

The third convention concerns the notation for the restriction to $C$ of a
differential form from $S^1 \x  B^3$: No special notation will be used.
That is,
the restricted form and the original form will not be notationally
distinguished. The context will usually make things clear.

With regard to the previous two conventions, the reader should take special
care to note that the norm on $C$ of a restricted form may differ
 from its
norm in $S^1 \x  B^3$ although both norms will be written alike. (Note that the
former norm is never greater than the latter.)

\subsection*{1.c\qua A comment about a Riemannian generalizaton}

This subsection constitutes a digression of sorts in order to provide the
promised brief remarks concerning the relevance of Theorem 1.8 in the wider
context of degenerate symplectic forms on 4--manifolds which come from
closed, self-dual 2--forms. (See also \cite{T1}.)

To start the digression, remark that if $X$ is a compact, oriented,
Riemannian
4--manifold, then there are non-trivial, harmonic self-dual 2--forms if the
intersection form on $H_2(X; \Bbb Z)$ is not negative definite. (Celebrated
theorems of Donaldson \cite{D1} and Freedman \cite{Fr} assert that a smooth, simply
connected 4--manifolds with negative definite intersection form is
homeomorphic to the multiple connect sum of $S^4$ with some number of
${\Bbb C{\Bbb P}}^2$'s having the non-complex orientation.)

With the preceding understood, let $X$ be a compact, oriented, 4--manifold
whose intersection form is not negative definite and let $\o_X$ be a
non-trivial, closed and self-dual 2--form on $X$. As $\o_X \wedge\o_X = |\o_X|^2
d$vol, this form is symplectic where it does not vanish. For some time, a
number
of geometers have been aware (see, eg \cite{T2} and \cite{D2}) of the fact
 that when
the metric on X is suitably generic, then there is a closed, self-dual
2--form $\omega_X$ whose zero set is cut out transversely as a section of the
self-dual 3--plane subbundle of $\L^2T^*X$. (See, eg \cite{Le} or \cite{Ho1} for a
 published
proof of this fact.) When such is the case, then $\o_X^{-1}(0)$ is a union of
embedded circles in $X$. Furthermore, a tubular neighborhood of such a circle
can be parametrized as $S^1 \x  B^3$ with coordinates $t \in S^1 = \R/\Bbb Z$
and $(x_1, x_2, x_3) \in B^3$ so that $\o_X$ has the form
\addtocounter{equation}{1}
\begin{equation} %1.14
\o_x=dt\wedge A_{ij} x^idx^j+A_{ij}x^i\e^{jkm}dx^k\wedge dx^m+
{\cal O}(|x|^2)
\end{equation}
where $\{A_{ij} = A_{ij}(t)\}_{i,j=1,\dots,3}$ are the components of a
$t$--dependent, symmetric, trace zero and non-degenerate $3\x  3$
 matrix. (In
(1.14), repeated indices are summed. Also, $\{\e^{jkm}\}$ are the components
of the completely anti-symmetric, $3\x 3$ matrix with $\e^{123} = 1$.)
Note that
by changing the orientation of $S^1$ if necessary, one can assume that
$A_{ij}(t)$ has, at each $t$, precisely one negative eigenvalue and
two positive
eigenvalues.

After passing to a double cover of $S^1\x B^3$ (if necessary), one can
find a
further change of coordinates (of the form $x_j \to  B_{jk} x_k$, where
$(B_{jk})$ is a $t$--dependent, invertible matrix) so that in these new
coordinates, $\o_X = \o+{\cal O}(\Sig_j |x_j|^2$), where $\o$ is as
 presented
in (1.2). Thus, up to a diffeomorphism (possibly on the double cover of $S^1\x
B^3$) the form $\o_X$ is given by (1.2) to leading order in the distance to
$\o_X^{-1}(0)$. (Note that this double cover is required when the negative
eigenspace of the matrix $(A_{ij})$ in (1.14) defines the unoriented line
bundle
over $S^1$. As pointed out by Luttinger \cite{Lu}, an example for the latter case
occurs when $\o_X$ is the self-dual form on $\R^4 = {\Bbb C}^2$ for the
standard Euclidian metric whose expression in complex coordinates $(w_1,
w_2)$ is
\begin{eqnarray}
&& i (1-|w_1|^2+|w_2|^2)(dw_1\wedge d\bar w_1+dw_2\wedge d\bar w_2) \nonumber\\
&& +i(R^{-1}w_2+\bar w_1\bar w_2)dw_1\wedge dw_2 \nonumber\\
&& -i (R^{-1}\bar w_2 +w_1 w_2)d\bar w_1\wedge d\bar w_2.
\end{eqnarray}
Here, take $R>>1$.  The zero set of this form consists of the unit
circle in the plane where $w_2 = 0$ and another circle where
$|w_1 | = R$.  The negative eigenspace of the matrix $(A_{ij})$
for the former component is unoriented.)

Moser's trick can be used (the details are written down in \cite{Ho2}) to make an
additional coordinate change so that $\o=\o_X$ near where $\o_X = 0$
(perhaps on
the double cover of $S^1\x B^3$). However, the latter change of
coordinates may
not be ${\Bbb C}^{\i}$.

Alternately, $\o_X$ can be homotoped through forms with fixed,
non-degenerate
vanishing set which are symplectic on the complement of this fixed vanishing
set, and so that the end result is a closed 2--form which agrees with (1.2)'s
form $\o$ (up to a diffeomorphism on $S^1\x B^3$ or its double
cover) near where
$\o_X = 0$. Note that this homotopy can be made so as to have non-trivial
support in any given neighborhood of fixed vanishing set. Furthermore, one
can prove directly that the form at the end of the homotopy is self-dual
with respect to a smooth metric, $g_X$, on $X$ which, on the
 given $S^1\x B^3$ or
its double cover, is the standard metric near $S^1\x\{0\}$. (See also
\cite{Ho2} for this last point.)

A form $\o_X$ as just described has a compatible almost complex
structure, $J_X$,
whose restriction to the given $S^1\x B^3$ (or to its double cover) is
equal to
$J$ from (1.4) in the coordinates where the form is given by (1.2).
Indeed,
take the metric $g_X$ from the preceding paragraph and set
$J_X = \sqrt{2}g_X^{-1}\o_X/|\o_X|$.

Given any $\o_X$--compatible almost complex structure $J_X$ for $Y
\!=\! X \!-\!
\o_X^{-1}(0)$, the notion of a finite energy, pseudo-holomorphic subvariety is
defined as in Definitions 1.1 and 1.2 using $\o = \o_X$ and $J = J_X$.

By restriction to a neighborhood of $\o_X^{-1}(0)$, Theorem 1.8
describes such finite energy, pseudo-holomorphic subvarieties in the
case when there are coordinates $(t, x_1, x_2, x_3)$ for each
component of a tubular neighborhood (or its double cover) of
$\o_X^{-1}(0)$ with the property that both $\o_X = \o$ from (1.2) and
$J_X= J$ from (1.4). Theorem 1.8 says nothing directly in the case
where $\o_X$ and $J_X$ is not as strongly constrained near
$\o_X^{-1}(0)$.

The study of finite energy, pseudo-holomorphic subvarieties in the
more general case will be deferred to a (planned) sequel to this
article.  However, note that \cite{T5} proves that pseudoholomorphic
subvarieties exist in $X-Z$ when $X$ has a nonzero Seiberg-Witten
invariant.  Also \cite{T1} and \cite{T6} contain various speculative
remarks about the significance of the pseudoholomorphic varieties
studied here in the context of 4--manifold differential topology.

\subsection*{1.d\qua More about Example 1.5}

The purpose of this subsection is to prove various assertions in the
discussion surrounding the solutions to (1.6) and the corresponding
pseudo-holomorphic submanifolds.  These assertions are restated as
Lemmas 1.9 and 1.10.  This subsection can be skipped at first reading.

\med {\bf Lemma 1.9}\qua{\sl{Let} $c \in (0, 2/3^{3/2}]$. {Then, up
to the
action of translation on $\R$, there is a unique solution}, $u$, {to} (1.6)
{which obeys}
$\frac{9}{4} (u' )^2 + u^2 + c/u = 1$ {and is somewhere
positive.
Furthermore, $u$ is a periodic function of $t$ and there exist $c \neq
2/3^{3/2}$ for which the period is a rational multiple of $\pi$. Indeed, the
assignment to $c$ of the period of the associated solution $u$
 defines a
non-constant, analytic function on} $(0, 2/3^{3/2})$.}

\med
{\bf Lemma 1.10}\qua{\sl {There is a constant $\z\geq 1$ with the following
significance: Suppose that} $c \in (0, 2/3^{3/2}]$ {is such that the
corresponding solution to} (1.6) {has rational periods. Let $\Sig_{\pm}$
denote the corresponding submanifold as described in} (1.7).  {Let $r > 0$
and let $B_r \subset S^1\x  B^3$ denote the ball with center $(t_0, 0)$ and
radius $r$. Then the integral of $\o$ over $\Sig_{\pm}\cap B_r$ has the form
$r^3 \d_c$ where} $\d_c\geq  \z^{-1}$.}

\med

The remainder of this subsection is occupied with the proofs of these two
lemmas.

\med

{\bf Proof of Lemma 1.9}\qua To find a solution to (1.6), consider first the
same equation where $c = 1$. Write the latter as the dynamical system:
\begin{itemize} %1.16
\item\quad   $u'=p$
\item\quad $p'= -\textstyle{\frac 49} u+
  \textstyle{\frac 29} u^{-2}$ \hspace\fill (1.16)
\end{itemize}

Standard existence and uniqueness theorems for integrating vector fields on
the plane imply that each point in the $(u, p)$ plane where $u > 0$
sits on a
unique trajectory of (1.16). Moreover, the function $\frac{4}{9} p^2 + u^2 +
1/u$ is constant on each trajectory. (Just differentiate and use
 (1.16).)
This constant is necessarily positive on a trajectory where $u > 0$. Then,
the constancy of $\frac{4}{9} p^2 + u^2 + 1/u$ along trajectories
implies
that each trajectory is compact and thus diffeomorphic to a circle in the $u
> 0$ portion of the $(u, p)$ plane. In terms of (1.6), the preceding
implies
that the c = 1 version of (1.6) has a unique $u > 0$ solution with given
initial conditions $(u(0) = u_0 > 0$, \ $u'(0) = p_0$). Moreover,
this solution is periodic.

Now, suppose that $(u, p)$ obey (1.16) and $\frac{4}{9}  p^2 + u^2
 + 1/u = e^2$
for some $e > 0$. (By the way, note that $u^2 + 1/u \geq  3/2^{2/3}$ and so
$e \geq 3^{1/2}/2^{1/3}$ with equality if and only if $u = 1/2^{1/3}$.)
Then ${\underline u}\equiv u/e$ obeys (1.6) with the normalization
$\frac{4}{9} {\underline u}^{'2}+{\underline u}^2+c/{\underline u}=1$
 with $c = e^{-3}$. (The preceding
lower bound on $e$ implies the upper bound $c \leq  2/3^{3/2}$.)
Conversely, if
${\underline u}$ obeys (1.6) with the normalization,
$\frac{4}{9} {\underline u}^{'2}+{\underline u}^2+c/{\underline u}
 = 1$ then $u = e{\underline u}$
obeys (1.16) for $c = 1$ and with the normalization $\frac{4}{9}  p^2 +
u^2 + 1/u
= c^{-2/3}$. With this last fact understood, the existence and uniqueness
theorems for (1.16) directly imply the asserted existence and uniqueness for
solutions to (1.6).

Now consider the assertion that there exists solutions with periods which
are rational multiples of $\pi$. Since the period depends continuously on
the parameter $c$, it is sufficient to prove that the period is not
constant.
The proof of this fact below will also establish the final assertion of the
lemma, namely that the period defines a non-constant, analytic function on
$(0, 2/3^{3/2})$.

To begin, note that the period equals twice the time it takes for $u$ to
travel from its minimum value to its maximum value, and the latter time
 is equal to
\addtocounter{equation}{1}
\begin{equation} %1.17
{\textstyle{\frac 32}}\int^{u_{\max}}_{u_{\min}} (1-x^2-c/x)^{-\frac
12}dx.
\end{equation}
Thus, to estimate the period, estimates for the integral in (1.17) are
required. For this purpose, note that $u_{\min}< u_{\max}$ are the two
 positive
roots of the polynomial $x^3 - x + c$. (For $c \in (0, 2/3^{3/2}]$, this
polynomial has two positive roots and one negative root.) To find these
roots, note that when $c = 2/3^{3/2}$, the polynomial has a double
root at
$1/3^{1/2}$ and a single root at $-2/3^{1/2}$. In particular, the roots
have the
form $\l\pm\d, -2\l$, where $\l, \d > 0$ and $\d <\l$. Moreover,
 when $c$ is
close to $2/3^{3/2}$, then $\l$ is close to $1/3^{1/2}$ and $\d$ is
close to
zero. With the preceding understood, one can rewrite (1.17) by changing
variables to $y = (x - \l)/\d$. The resulting integral is
\begin{equation} %1.18
3^{\frac 12}/2 \int_{[-1,1]} (1-y^2)^{-\frac 12}
(1+\d y/\l)^{\frac 12} (1+\d y/(3\l))^{-\frac 12} dy .
\end{equation}
Note that as $0 < \d/\l < 1$, this integral is an analytic function of the
ratio $\d/\l$. This combination is an analytic function of the roots of the
polynomial $x^3 - x + c$, and thus an analytic function of $c$
except at values
of $c$ where the roots are not distinct. The latter does not occur
in the
interval $(0, 2/3^{3/2})$.

Thus, it remains only to prove that (1.18) is not a constant function of c.
To argue this point, consider this last integral as a function of $\d/\l$. As
such, it can be written as a power series $\d/\l$ when $|\d/\l| << 1$. In
particular, to order $(\d/\l)^2$, the integral in (1.18) is equal to
\begin{equation} %1.19
3^{\frac 12}/2 \int_{[-1,1]} (1-y^2)^{-\frac 12}
[1-\d^2y^2/(6\l^2)] dy .
\end{equation}
This last integral can be evaluated using the trigonometric substitution $y =
\cos(\th)$ and yields
\begin{equation} %1.20
3^{\frac 12}/2\pi [1-\d^2/(12\l^2)]  .
\end{equation}
Note that (1.20) already implies that the period is not a constant function
of $c$ since $\d^2/\l^2$ is not a constant function of $c$. (Note
 that $\d/\l =
0$ at $c = 2/3^{3/2}$ but not otherwise.)

Although this last observation completes the proof of Lemma 1.9, a
digression follows (for the sake of completeness) which uses (1.20) to
estimate the period as a function of c which is accurate to order $(2/3^{3/2} -
c)^{3/2}$. The digression starts with a determination of the constants $\d$ and
$\l$ in terms of $\a$. For this purpose, it is necessary to find a
dependence of
the roots of polynomial $x^3 - x + 2/3^{3/2} - \a$. For small $\a$, it
turns out
that these roots can be expanded as a power series $\a^{1/2}$. In
particular, for
small $\a$, the roots are
\begin{equation} %1.21
1/3^{\frac 12}\pm 3^{-\frac 14}\a^{\frac 12}-\a/6 +
{\cal O}(\a^{\frac 32}) \quad {\mbox{and}} \quad
-2/3^{\frac 12}+\a/3+ {\cal O}(\a^2)  .
\end{equation}
(This last assertion can be checked by substitution.) Note that this is an
expansion in $\sqrt{\a}$ rather than $\a$ because the expansion is being done
around the value $c = 2/3^{3/2}$ where the polynomial $x^3 - x + c$ has a
double
root. In any event, these last expressions identify $\l = 1/3^{1/2} -
 \a/6$ and
$\d = 3^{-1/4} \a^{1/2}$. Then, with these last identifications and (1.20)
understood, the period for small $\a$ is seen to be
\begin{equation} %1.22
3^{\frac 12}\pi [1-\a/(4\sqrt{3})]+{\cal O}(\a^{\frac 32}) .
\end{equation}

\med
{\bf Proof of Lemma 1.10}\qua As remarked previously, because $\o$ is a
homogeneous function of the coordinates $(t, x, y, z)$ and because
$\Sig_{\pm}$
is a cone, it follows that the integral of $\o$ over $\Sig_{\pm}\cap
B_r$ is
finite and has the form $r^3 \d_c$ where $\d_c$ is independent
of $r$. Thus, it
remains only to bound $\d_c$ from below with a positive, $c$--independent
constant. For this purpose, consider the case where $c > 0$. (The $c < 0$
case is
treated almost identically.)  Now write $\o$ as in (1.10). Thus, the
integral of
$\o$ over $\Sig_{\pm}\cap B_1$ is a sum of two integrals, the first being
that of
$dt\wedge  df$ and the second that of $d\phi\wedge  dh$.

The integral of the latter is equal to $2\pi bc$ as the following argument
shows: First, the second and third points of (1.7) assert that $h = c s^3$
 on
$\Sig_+$. Thus, $dh = c s^2 ds$ on $\Sig_+$. With this understood, Stoke's
theorem writes the integral of $d\phi\wedge  dh$ over $\Sig_+ \cap B_1$ as
$c$ times
the integral of $\phi' d\tau$ over $[0, 2\pi b]$. (This last expression is the
pull-back via the $\tau$ parametrization of the restriction of $d\phi$ to
$\Sig_+\cap B_1$). The latter integral is $2\pi b$, a fact which
follows from
the fourth point in (1.7).

Given the preceding paragraph, a positive lower bound for $\d_c$ follows from
the following claim: There exists a constant $\xi \geq 1$ with the property
that $\xi^{-1}$ bounds the integral of $dt\wedge  df$ over $\Sig_+ \cap B_1$
 when
$c < \xi^{-1}$.

To prove this claim, use (1.7) and Stoke's theorem to write the integral of
$dt\wedge  df$ over  $\Sig_+ \cap B_1$ as the integral over the interval $[0,
2\pi b]$ of the function $(3/4) (2u^2 - c/u) u''$. (This last
expression is the pull-back via the t parametrization of the restriction of
$-f \ dt$ to $\Sig_+ \cap B_1$.) Next, use (1.6) to write this last
integral as
\begin{equation} %1.23
\int_{[0,2\pi b]} 3^{-1}u^{-1} (2u^2-c/u)^2d\tau .
\end{equation}
To proceed further, consider the maximum value for $u$ on $[0, 2\pi b]$.
Since $u'  = 0$ at this point, it follows from the normalization $\frac{9}{4}
(u')^2 + u^2 + c/u = 1$ that $u^2 + c/u = 1$ at maximum. Moreover, since
$u''\leq 0$, one can conclude from (1.6) that $u^2\geq  c/2u$ at maximum. Thus,
$u^2\geq  1/3$ at its maximum. Now, let $\tau_0$ be a point at which $u$
achieves
this maximum. Since the aforementioned normalization condition implies that
$|u'|\leq  2/3$, it follows that $u \geq  1/(2\sqrt{3})$ at all $t$ with
$|\tau-\tau_0| \leq  \sqrt{3/4}$. Since $u \leq  1$ (because of the
 normalization
condition), these last points imply that the integral in (1.23) is no
smaller than
\begin{equation} %1.24
(144\sqrt{3})^{-1}(1-12\sqrt{3}c)^2 
\end{equation}
when $c <\sqrt{3/12}$. Thus, when $c < \sqrt{3/24}$, the integral of $dt\wedge
df$
over $\Sig_+ \cap B_1$ is no smaller than $(576\sqrt{3})^{-1}$.

\section{Finite energy and finite local energy}
\setcounter{equation}{0}

The key step in analyzing the behavior near $Z$ of a finite energy,
pseudo-holomorphic subvariety $C$ is the local estimate given in
Proposition
2.1, below, for the amount of ``energy'' which can accumulate in a small
ball about any point in $Z$.

To start, fix a point, $t_0$, in $Z$ which can be taken to be the origin
in $S^1
= \R/\Bbb Z$. Let $\chi$ be a fixed ``bump'' function on $[0,\i)$;
that is, $\chi$ is a function which is non-increasing, equals 1 on [0,1)
 and
equals 0 on $[2, \i )$. Given $s > 0$, promote $\chi$ to a
function $\chi_s$ on
$X$ by the rule
\begin{equation} %2.1
\chi_s=\chi(r/s)
\end{equation}
where $r = (t^2 + \rho^2 + z^2)^{1/2}$.

Now, for each $s > 0$ introduce the ``local'' energy
\begin{equation} %2.2
\s(s)=\int_c\chi_s\o .
\end{equation}
This section proves the following estimate for $\s(s)$:

\med
{\bf Proposition 2.1}\qua{\sl {There exists $\z > 1$ such that for any point $t_0
\in Z$, the function $\s$ is differentiable on $(0, \z^{-1})$.
Furthermore, $\s$
obeys the inequalities:}
\begin{enumerate}
\item $\s(s) \leq z\dt s^3$ 
\item {The area of $C \cap  B(s)$ is bounded by} $\z\dt s^2$.
\item $\frac{d}{ds} \s(s) \geq  \frac{3}{s}\s(s)$\hspace\fill {\rm(2.3)}
\end{enumerate}}

\med
Here is an immediate corollary, a bound for the energy and area of the
intersection of $C$ with a tubular neighborhood of $Z$:
\med

{\bf Proposition 2.2}\qua{\sl {Let $C$ be a finite energy, pseudo-holomorphic
submanifold. Then there exists $\z$ with the following significance:  Fix $s >
0$ and let $C(s)$ denote the intersection of $C$ with the subset of $X$ where}
$(\rho^2+ z^2)^{1/2} < s$. {Then:}

\begin{itemize}
\item[\rm(a)] $\int_{C(s)}\o \leq  \z\dt s^2$

\item[\rm(b)] {The area of $C(s)$ is bounded by} $\z\dt s$.
\end{itemize}}

\med The remainder of this section is occupied with the proof of
Proposition 2.1.

\medskip
{\bf Proof of Proposition 2.1}\qua Consider first the issue of the
differentiability of the function $\s$. In this regard, it proves
useful to
introduce a second bump function: Given $\e > 0$, define the
function $\eta_\e$
on $X$ by the rule $\eta_\e = (1 - \chi((\rho^2 + z^2)^{1/2}/\e))$. Note,
 in
particular, that the support of $\eta_\e$ has closure which is disjoint
from $Z$.

Next, note that for $s$ fixed, the function of $\e$ given by
\addtocounter{equation}{1}
\begin{equation} %2.4
\s(s,\e)=\int_c \eta_\e\dt\chi_s\dt\o
\end{equation}
is a decreasing function of $\e$ with
\begin{equation} %2.5
\lim_{\e\to 0}\s(s,\e)=\s(s) .
\end{equation}
(Remember that the restriction of $\o$ to $C$ is equal to $g$ times the induced
area form on $C$ since $C$ is pseudo-holomorphic.) Note as well that the limit
in question is uniform in $s$. Infact, given $\d > 0$, there exists
$\e(\d)$ such
that for all $\e, \e' < \e(\d)$, and all $s$, one has
\begin{equation} %2.6
|\s(s,\e)-\s(s,\e')| < \d .
\end{equation}
Indeed, (2.6) stems from the following observations: First, because $C$ is
pseudo-holomorphic, the restriction of $\o$ to TC is a positive multiple of the
area form of the induced metric. In fact, this multiple is the function $g$ in
(1.3). Thus, when $\e > \e'$, then the difference in (2.6) is
not greater than the integral of $\o$ over the subset of points in $C$
where
$(\rho^2
+ z^2)^{1/2} \leq  2\dt \e$. And, the latter can be made as small as
desired by
choosing $\e$ small since the integral over $C$ of $\o$ is assumed finite.

By inspection the function $\s(\dt , \e)$ is continuous for fixed $\e$, and
so it
follows from (2.6) that $\s$ is too.

As for the derivative of $\s$, a similar argument finds $\s$ continuously
differentiable with
\begin{equation}  %2.7
\s'(s)=s^{-2}\dt\int_c\chi'_s \dt r\dt\o .
\end{equation}
(Note that the right hand side in (2.7) is no greater than
$\z\dt s^{-1}\dt \s(2s)$ where $\z$ is independent of $s$.)

As for the estimates in (2.3), remark first that the first and second
estimates follow by integration from the third and from the fact that $\omega$
on TC is $g$ times the area form. To consider the estimate in (2.3.3),
introduce
the 1--form
\begin{equation} %2.8
\th=3^{-1}(t\dt df-2f\dt dt-3h\dt d\phi).
\end{equation}
Note that $d\th= \o$. Thus, integration finds
\begin{eqnarray} %2.9
\s(s,\e) & = & -\int_c\eta_\e\dt d\chi_s\wedge\th
 -\int_c\chi_s\dt d\eta_\e\wedge\th\nonumber\\
 & \leq  & s^{-1}\int_c
\chi'_s\dt\kappa \dt |\th|\dt g^{-1}\dt\o +\int_c
|d\eta_\e|\dt |\th|\dt g^{-1}\dt\o.
\end{eqnarray} 
Here, $\k$ is the norm of pull-back to $C$ of the 1--form $|dr|$.  (The
right
most inequality in (2.9) exploits the fact that $\o$ on TC is equal to
$g$ times
the area form of the induced metric.)

To make use of (2.6), remark that a calculation finds,
\begin{equation} %2.10
|\th|\leq 3^{-1}\dt (t^2\dt |df|^2+4\dt |f|^2\dt |dt|^2
+9\dt |h|^2\dt |d\phi|^2)^{12} \leq 3^{-1}\dt r\dt g .
\end{equation}
(The second inequality follows from the identity
\begin{equation} %2.11
t^2\dt |df|^2+4\dt |f|^2\dt |dt|^2
+9\dt |h|^2\dt |d\phi|^2 = (t^2+\rho^2 +z^2)\dt g^2
\end{equation}
which holds on $S^1  \x\R^3$.) Thus, according to (2.9) and (2.10),
\begin{equation} %2.12
\s(s,\e)\leq 3^{-1}\dt s^{-1}\dt\int_c
\chi_s\dt\kappa\dt r\dt\o +\zeta\dt\int_c \eta_{\e/4}\dt
(1-\eta_{4\e})\dt\o .
\end{equation}
As $\e \to  0$, the second term on the right side in (2.12) vanishes,
while the first term on the right side is less than $s\dt \s'(s)/3$ since $\k
\leq  1$ always. (Compare with (2.7).)

\section{Limits}
\setcounter{equation}{0}

Fix a point $t_0 \in Z$ as in the previous section, and again call this
point 0.  For each positive $s$, introduce $B(s)=\{(t, x, y, z): r = (t^2 +
x^2 + y^2 + z^2)^{1/2} < s\}$. Then, for $R \geq 2$ and for
any $s < R/2$,
define $C_s  \subset B(R)$ by the rule
\begin{equation} %3.1
C_s=\{(t,x,y,z)\in B(R): s\dt (t,x,y,z)\in C\cap B(R\dt s)\}.
\end{equation}
Of course, $C_s$ is obtained from $C \cap  B(R\dt s)$ by dilating all of the
coordinates. Since $\o$ is homogeneous under dilation, and the metric
as well,
it follows immediately that $C_s$ is a $J$--pseudo-holomorphic subvariety of
$B(R)$. Furthermore, it follows from Proposition 2.1 that the integral over
this
subvariety of $\o$ is bounded independently of $s$. Also, the area of
$C_s$ is
bounded independently of~$s$.

With the preceding understood, one can consider the limiting behavior as $s
\to  0$ of the family $\{C_s\}$. The next definition describes a
useful notion of convergence. However, a preliminary digression is required
before this definition. The purpose of the digression is to introduce the
notion of an {irreducible component} of a pseudo-holomorphic subvariety.
For this purpose, suppose that $C$ is a finite energy, pseudo-holomorphic
subvariety in a symplect 4--manifold $Y$ with compatible almost complex
structure. By definition, there is a countable set $\L \subset C$ whose
complement is a submanifold of $Y$. Take any such set $\L$. An irreducible
component of $C$ is the closure of a component of $C - \L$. (The set of
irreducible components is insensitive to the choice of $\L$ provided
 that $\L$
is countable and $C - \L$ is a submanifold.)

Note that each irreducible component of a finite energy, pseudo-holomorphic
subvariety is also a finite energy, pseudo-holomorphic subvariety.  (The sum
of the energies of the irreducible components of $C$ equals that of $C$.)

End the digression. Here is the promised convergence definition:

\med
{\bf Definition 3.1}\qua Let $Y$ be a (possibly non-compact) 4--manifold with
symplectic form $\o$ and with a compatible almost complex structure
$J$. Let $\O$
be a countable, decreasing sequence of positive numbers with limit
 zero and let
$\{C_s\}_{s\in\O}\subset  Y$ be a corresponding sequence of finite energy,
pseudo-holomorphic subvarieties. Let $C \subset Y$ be a finite energy,
pseudo-holomorphic subvariety as well. Say that $\{C_s\}$ {\sl converges
geometrically} to $C$ when the following two requirements are met:

$\bullet$\qua Treat each $C_s$ and also each irreducible
 component of $C$ as a
2--dimensional, rectifiable current which assigns to a smooth 2--form with
compact support on $Y$ the integral of the form over the set in question.
 Let
${\cal{C}}$ denote the set of irreducible components of $C$. Then, require
that the
sequence of currents $\{C_s\}$ converges weakly as a current to the current
$\Sig_{C'\in{\cal{C}}} m_{c'} C'$, where $m_C$ is, for each $C \in
{\cal{C}}$, a positive integer.

$\bullet$\qua For each compact subset $K\subset Y$, introduce the number,
$d_K(s)$
which is the supremum over pairs $(x, y) \in (C_s \cap  K)\x (C_0
\cap  K)$
of the sum of the distances from $x$ to $C_0$ and $y$ to $C_s$. Then, for
 each
such $K$, require that the limit as $s \to  0$ of $d_K(s)$
exists and is zero.

\med
With the preceding definition understood, here is the main result of this
section:

\med
{\bf Proposition 3.2}\qua{\sl {Let $\{C_s\}_{s >0}$ be defined as in} (3.1).
{Let $\{s_i\}$ be a countable, decreasing
sequence of
$s$ values with limit zero. This sequence has a subsequence for which the
corresponding set $\{C_s\}$ geometrically converges in each $Y = B(2) -
 Z$ to
a finite energy, pseudo-holomorphic submanifold $C_0 \subset B(2) - Z$.
Furthermore, each connected component of $C_0$ in $B(1) - Z$ is identical
 to the
intersection of $B(1) - Z$ with one of the submanifolds from
Example 1.3, 1.4 or 1.5.}}

\medskip
The proof of this proposition exploits the following local version of
Gromov's compactness theorem \cite{Gr} (see also \cite{Pan}, \cite{PW}, \cite{Ye}, \cite{MS}).

\med
{\bf Proposition 3.3}\qua{\sl Let $Y$ be a smooth 4--manifold with symplectic
form
$\o$ and compatible almost complex structure $J$. Let $\O$ be a countable,
decreasing set of real numbers with limit zero and let
$\{C_s\}_{s\in\O}\subset
Y$ be a corresponding set of finite energy, pseudo-holomorphic subvarieties.
Suppose that there exists $E > 0$ such that $\int_C \o < E$ for all
$C\subset\{C_s\}$. Then there is a finite energy, pseudo-holomorphic subvariety
$C_0 \subset Y$ and an infinite subsequence of $\O$ for which the
 corresponding
subset $\{C_s\}$ converges geometrically to $C_0$.}\medskip

This last proposition is proved at the end of this section. Accept it
momentarily for the proof of Proposition 3.2.

\medskip
{\bf Proof of Proposition 3.2}\qua As remarked above, the sequence $\{C_s\}$
has uniformly bounded energy, so Proposition 3.3 provides a subsequence
which converges geometrically to the limit $C_0$. As $C_0$ has finite
energy, the
results of Section~2 apply to the subvariety $C_0$. For any $a \in (0, 1)$,
consider the value of the function $\s( \ )$ on $C_0$ at the point $a$.
(Take the
point $t_0 = 0 \in Z$ to define this function.) By the convergence
criteria,
$\s(a)$ for $C_0$ is the limit as $s \to  0$ of the values of $m^{-1}\s(a)$ for
$C_s$ where $m$ is a positive integer.
And, by rescaling, this is equal to
\begin{equation} %3.2
m^{-1}\dt (Rs)^{-3}\dt\s(R\dt s\dt a)
\end{equation}
as defined for the curve $C$.

On the other hand, the function $s^{-3}\s(s)$ on $C$ is a monotonically
increasing function of $s$, which means, since it is bounded, that there
is a
unique limit, $\s_0$, as $s \to  0$. Thus, as $s$ tends to zero,
(3.2) converges to
$a^3\dt m^{-1}\dt \s_0$. As a result, one can conclude that the value of
$\s(\dt)$ on $C_0$ at the point $a$ is equal to $a^3\dt m^{-1}\dt \s_0$, with
$\s_0$ being independent of $a$.

Meanwhile, according to Proposition 2.1, the derivative with respect to a of
$a^{-3}\dt \s(a)$ is non-negative. In the case at hand, this derivative
is zero.
According to (2.7)-(2.12), this derivative is zero only if a number of
inequalities are equalities. In particular, the function $\k$ in (2.12) must
equal 1 everywhere on $C_0$ which can only happen when the vector field
\begin{equation} %3.3
s\p_s=t\p_t+\rho\p\rho +z\p_z
\end{equation}
is everywhere tangent to $C_0$. Here, $s^2 = t^2 + \rho^2 + z^2$.
 Then, $J\dt (s
\p_s)$ must also be tangent to $C_0$. It then follows that $C_0$ is the
cone on a
union of integral curves of $J\dt (s \p_s)$. Furthermore, as no two distinct
integral curves intersect, it follows that each component of $C_0$ is
submanifold, and that no two distinct components intersect.

The analysis of the integral curves of $J\dt (s \p_s)$ is facilitated by
writing the latter vector field as
\begin{equation} %3.4
J\dt(s\p_s)=-gt\p_f +2g^{-1}f\p_t+3(g\rho^2)^{-1} h\p_\phi .
\end{equation}
(This equation follows from (1.10) and (1.11) and the identity $s \p_s =
t \p_t + 2f \p_f + 3h \p_h$.) Now, remark that $J\dt (s \p_s)$ annihilates not
only the function $s$, but also the function $h$. Thus, $h =$ constant on an
integral curve of $J\dt (s \p_s)$. If $h = 0$, then the
resulting component of
$C_0$ is one of the cases from Example 1.3 and Example 1.4. At most
finitely many
components of $C_0$ are described in this way because a calculation shows that
each component contributes at least $2 r^3/3$ to the integral of $\o$
 over $C_0
\cap  B_r$. Thus, Proposition 3.2 follows with a demonstration that each $h
\neq
0$ case is described in Example 1.5, and that there are at most
finitely many
of these.

To begin the study of the $h \neq  0$ case, assume, for the
 sake of argument,
$h > 0$ on an integral curve $\g$ of (3.4). (The case $h < 0$ is
argued
analogously.)  Also, normalize the curve so that $s\equiv 1$ on $\g$.

The first observation is that the condition $h =$ constant $> 0$
implies that
the $\R/(2\pi\Bbb Z)$ valued function $\phi$ is increasing along $\g$.
Thus, $\phi$ restricts to $\g$ as a local coordinate and the
corresponding component
of
$C_0$ can be parametrized, at least locally, by setting $t = s {\underline
t}(\phi)$, $f = s^2 {\underline f}(\phi)$ and $h = c s^3$, where $c$ is
a positive
constant and where ${\underline t}$ and ${\underline f}$ are functions
of $\phi$ as indicated.

The condition that $\g$ is an integral curve of $J\dt (s \p_s)$
implies that the
functions ${\underline t}$ and ${\underline f}$  are constrained to obey the
equations:
\begin{eqnarray}
3c\underline{t}_\phi= 2\underline{\rho}^2\underline{f}\nonumber\\
3c\underline{f}_\phi = -\underline{g}^2\underline{\rho}^2\underline{t}
\end{eqnarray}
Here, ${\underline g}$ and ${\underline\rho}$ are the restrictions of the the
functions $g$ and $\rho$ to $\g$. To analyze (3.5), one can exploit the
fact that
both ${\underline t}$ and ${\underline f}$ along $\g$ can be written
implicitly
as functions of the restriction, ${\underline z}$, of $z$ to $\g$. (This is
because $s$ and $h$ are constant on $\g$). Indeed, ${\underline\rho}^2 =
c/{\underline z}$ from the definition of $h$ and therefore $\mu = 2^{-1}
c/{\underline z} - {\underline z}^2$ while $\nu =
\pm  (1 - c/{\underline z} - {\underline z}^2)^{1/2}$. Moreover,
${\underline g}^2
= c/{\underline z} + 4 {\underline z}^2$, so each line in (3.5) gives
an equation
for ${\underline z}_{\phi}$ in terms of ${\underline z}$. These two
equations are
identical and imply that ${\underline z}_{\phi} = 2 {\underline t}/3$.
Differentiate the latter equation once to obtain a second order equation
for the
function ${\underline z}$. This second order equation is (1.6) with the
substitution of $u$ for ${\underline z}$ and $t$ for $\phi$. Then, the
description of the corresponding component of
$C_0$ by (1.7) follows directly from the preceding remarks. Thus, the cases
where $h
\neq   0$ are described by Example 1.5 as required.

To complete the proof of Proposition 3.2, remark that the number of such $h
\neq   0$ cases is finite as there is a $c$--independent lower bound for the
integral of $\o$ over the intersection with $B_1$ of any case from
Example 1.5.
(See Lemma 1.10.)

\medskip
{\bf Proof of Proposition 3.3}\qua
 The first step is to obtain the limit $C_0$ as
a subset of Y. The following monotonicity lemma plays a key role in this
step. (A similar lemma is proved in \cite{PW}.)

\medskip{\bf Lemma 3.4}\qua{\sl 
{Let $Y$ be a smooth 4--manifold with symplectic
form $\o$
and almost complex structure $J$. Let $E$ be a positive number and
let $K\subset
Y$ be a compact set. Then, there exists a constant $c = c(E, K) \geq  1$
with the
following significance: Let $C \subset Y$ be a finite energy,
pseudo-holomorphic
subvariety with $\int_C \o \leq  E$. Suppose that $x \in K \cap  C$. Use
$J$ and
$\o$ to define a Riemannian metric on $Y$, and for $r > 0$, let $B(x,
r)\subset
Y$ denote the ball of radius $r$ with center $x$. Then,}
\begin{equation} %3.6
c^{-1}r^2\leq\int_{B(r,x)\cap C}\o
\leq cr^2
\end{equation}
{if} $r \leq  c^{-1}$.}

\medskip
{\bf Proof of Lemma 3.4}\qua Since $K$ is compact, there exists $r_K$ so that for
each $x \in K$, the ball $B(r_K, x)$ has compact closure in $Y$, and so
that in
this coordinate system, $\o=\o_0+\o_1$, where $\o_0 = dx_1\wedge  dx_2 + dx_3\wedge
dx_4$ and where $|\o_1| \leq |x|$. With the preceding understood, consider
$f(r) \equiv\int_{ B(x,r)\cap C}\o$ and let $f_0(r) = \int_{B(x,r)\cap
C}\o_0$.
Note that both functions are continuous and piecewise differentiable on the
interval $[0, r_K/2]$. Furthermore, as $f(r)$ is the area of the
 intersection
of $C$ with $B(x, r)$, one has
\begin{equation} %3.8
|f(r)-f_0(r)|\leq r\dt f(r),
\end{equation}
as long as $r \leq  r_K/2$.

Meanwhile, $\o_0 = d\th$ for $\th = 2^{-1} (x_1dx_2 - x_2dx_1 +
 x_3dx_4 -
x_4dx_3)$, and so
\begin{equation} %3.7
f_0(r)=\int_{\p B(x,r)\cap C}\th\leq
\int_{\p B(x,r)\cap C}|\th| =2^{-1}r
\int_{\p B(x,r)\cap C} 1=2^{-1}rf'(r) .
\end{equation}
(In deriving the latter, note that the Euclidean ball of radius $r$ agrees
with the Gaussian ball of radius $r$.) Thus, the function $f(r)$
obeys the
differential inequality
\begin{equation} %3.8
f(r)\leq (1+r)\dt 2^{-1}rf'(r) .
\end{equation}
The latter can be integrated to find that when $r \leq  r'\leq r_K/2$.
Then
\begin{equation} %3.9
{r'}^{-2}f(r')\geq 2^{-1}r^{-2}f(r)
\end{equation}
when $r_K \leq 1$. Since $f(r' ) \leq  E$, this last
inequality implies first that $f(r) \leq  r^2 (8 E r_K^{-2})$ for any $r \leq
r_K/2$. This gives the right hand inequality in (3.6). The left-hand
bound
follows by taking $r$ very small on the right hand side of (3.10). In this
case, standard regularity results about the local behavior of
pseudo-holomorphic maps (as in \cite{PW} and \cite{Ye}) imply that for $r > 0$ but
very small, $f(r) \geq  \z^{-1} r^2$, where $\z$ is a fixed constant
which is
independent of all relevant variables. (The point here is that a
pseudo-holomorphic subvariety looks locally like a complex subvariety in
${\Bbb C}^2$, and such a lower bound for $f(r)$ can be readily found for
the latter.)

\med

With Lemma 3.4 understood, the construction of $C_0$ proceeds with the
following lemma:

\medskip
{\bf Lemma 3.5}\qua{\sl{Let $Y$ be a smooth 4--manifold with symplectic form $\o$
and almost complex structure $J$. Let $K\subset  Y$ be a compact
set and let $E$
be a positive number. Then, there exists a constant} $c = c(E, K) \geq  1$
{with
the following significance: Let $C\subset  Y$ be a finite energy,
pseudo-holomorphic subvariety with} $\int_C \o \leq  E$. {Then
given $r > 0$
but less than $c^{-1}$, there is a set of no more that $c r^{-2}$ balls of
 radius
$r$ with center on $C
\cap  K$ which covers $C \cap  K$. Furthermore, this set can be taken
so that the
concentric balls with radius r/4 are pairwise disjoint.}}

\medskip
{\bf Proof of Lemma 3.5}\qua First, choose $r_K$ so that for each $x \in K$, the
ball $B(r_k, x)$ has compact closure in $Y$. Then, for $r$ positive, but
less than
$r_K$, let $\Xi$ be a maximal set of balls of radius $r/4$ with
centers on $K$ which
are pairwise disjoint. Let $c$ be the constant which occurs in Lemma 3.4.
Then, there can be no more that $c E r^{-2}$ such balls because,
according to
Lemma 3.4, each ball contributes at least $c^{-1} r^2$ to the energy.
 Take the
cover to be the set of balls of radius $r$ which are concentric to the balls
in $\Xi$. (If this set did not cover $C \cap  K$, then $\Xi$ would not be
 maximal.)

\medskip
With Lemma 3.5 understood, obtain $C_0$ by the following limiting
argument:
Take an exhaustion of $Y$ by a countable, nested set of compact sets
$\{Y_p\}_{p=1,2,\dots}$. Here, $Y_p\subset  Y_{p+1}$ and $\cup_p Y_p =
Y$. For each
$p$ and integer $N$ sufficiently large, and for each $s \in\O$, choose
a set
$\Xi_{p,N,s}$ of no more than $c_p 4^2 16^{2N}$ balls of radius
 $4^{-1} 16^{-N}$ which
cover $C_s
\cap  Y_p$ and which have centers on $C_s \cap  Y_p$. (Here, take $N$
large enough so
that Lemmas 3.4 and 3.5 can be invoked for the case where $K = Y_p$ and $r
= 4^{-1}
16^{-N}$.) Add extra balls if required so that the number of elements in
$\Xi_{p,N,s}$
is some $s$ independent integer $m_{p,N} \leq  c_p 4^2 16^{2N}$. Label
the centers of
these balls to give a labled set ${\underline x}_{p,n,s}$ of
$m_{p,N}$ points in $Y_p$
which one can think of also a single point in the $m_{p,N}$ fold product
of $Y_p$. For
each fixed $p$ and $N$, the sequence $\{{\underline
x}_{p,N,s}\}_{s\in\O}\subset
\times_{m_{p,N}}Y_p$ has a
convergent subsequence. By taking a diagonal subsequence, one obtains a
subsequence
$\O'\subset\O$ for which each sequence $\{{\underline
x}_{p,N,s}\}_{s\in\O'}$ converges
for each $p$ and $N$. Let ${\underline x}_{p,N,0}$ denote the limit point,
but  thought of as a labeled set of $m_{p,N}$ points in $Y$. Let
$U_{p,N,0}\subset  Y$
denote the set of points with distance $16^{-N}$ or less from some point in
${\underline x}_{p,N,0}$. The claim now is that given $p$, there exists
exists $N(p)$
such that
\begin{equation} %3.11
U_{p',N+1,0}\cap Y_{p-1}\subset U_{p,N,0}
\end{equation}
for all $p'$  and for all $N \geq  N(p)$. Given this claim, take
\begin{equation} %3.12
C_0\cap Y_{p-1}\equiv \bigcap_{N\geq N(p)}
(U_{p,N,0}\cap Y_{p-1}) .
\end{equation}
Note that (3.11) insures that (3.12) is consistent with respect to the
inclusion $Y_p\subset  Y_{p+1}$.

To see (3.11), suppose to the contrary that
$U_{p' ,N+1,0} \cap  Y_{p-1}$ had a point which was not in $U_{p,N,0}$.
Then, for all small $s \in\O'$ , there would be a point $x \in
Y_{p-1}$ with distance $8^{-1} 16^{-N}$ or less from some point in
${\underline x}_{p' ,N+1,s}$ and which had distance $2^{-1} 16^{-N}$ or
more from each
point in ${\underline x}_{p,N,s}$. However, this is not possible if $N$
is large for
the following reason: Let $d_p$ denote the distance between $Y - Y_p$ and
 $Y_{p-1}$.
Let $x'  \in {\underline x}_{p' ,N+1,s}$ denote the point which
is closest to $x$. Then $x'  \in C_s$. Furthermore, as $x \in
Y_{p-1}$, then $x'  \in C_s \cap  Y_p$ if $8^{-1} 16^{-N} < d_p$. Of course,
if $N$ is large, this will be the case. On the other hand, if $x'
\in C_s \cap  Y_p$, then $x'$  has distance $4^{-1} 16^{-N}$ from
each
point in ${\underline x}_{p,N,s}$. Thus, the triangle inequality implies that
\begin{equation} %3.12
{\mbox{dist}}(x,\underline{x}_{p,N,s})\leq 4^{-1}
16^{-N}+8^{-1} 16^{-N} < 2^{-1}16^{-N} .
\end{equation}
With $C_0$ defined by (3.12), the claim next is that:
\begin{itemize} %3.14
\item $C_0$ is closed.
\item $C_0\cap Y_p$ has finite two dimensional Hausdorff measure.

\item For each $p$,
  \begin{equation}\lim_{s\in\O, \ s\to 0} \big(\sup_{y\in C_s\cap Y_p}
  {\mbox{dist}}(y,C_0)+\sup_{y\in C_0\cap Y_p}
  {\mbox{dist}}(y,C_s)\big)=0. 
  \end{equation} 
\end{itemize}

\noindent
Indeed, $C_0$ is closed as it is the intersection of closed sets. The
fact
that $C_0 \cap  Y_p$ has finite two dimensional Hausdorff
measure follows from
the construction and the definition of Hausdorff measure. The second point
in (3.14) follows from the left most inequality in (3.6).

As remarked previously, each $C_s$ can be thought of as a current,
 which is to
say a linear functional on the space of smooth, compactly supported 2--forms
on $Y$. As a current each $C_s$ has norm bounded by $E$. Since the
unit ball in the
space of currents is compact in the weak topology, it follows that there is
a subsequence $\O''\subset\O'$  for which the
corresponding sequence of currents $\{C_s\}_{s\in\O''}$
 converges to a current with norm no greater than $E$. The
latter has support in $C_0$ in the sense that it annihilates any
2--form whose
support avoids $C_0$. Use ${\underline C}'_0$  to denote the limit current.

The final step in the proof of Proposition 3.3 argues the following two
points:
\begin{itemize}
\item There is a complex curve $C'_0$ with a proper,
pseudo-holomorphic map $f'$ into $Y$ which embeds each
component off of a finite set of points and whose image is $C_0$.

\item The current $\underline{C}'_0$  is equal
to the composition of pull-back by $f$ and then integration
over the complex curve. \hfill (3.15)
\end{itemize}

Proposition 3.3 now follows directly from (3.15). (The curve $C_0$ in
Definition 1.1 is obtained by identifying components of $C'_0$
with the same image.)

The assertions in (3.15) will be proved by invoking the analysis in Section
6 of \cite{T3}. In particular, the first point of (3.15) follows from Proposition
6.1 in \cite{T3} with a demonstration that $C_0$ has a {positive cohomology
assignment} as defined in Section~6a of \cite{T3} and reproduced in the
following paragraph. However, with regard to the
application of Proposition 6.1 in \cite{T3}, the reader should take note that the
proof offered in \cite{T3} has errors which occur in Section~6e of the article.
A corrected proof is provided in the revised version of this same article
which appears in \cite{T4}. Also, remark that the compactness requirement for the
ambient manifold in the statement of Proposition 6.1 can be dropped in this
case without compromising the conclusions except to allow the complex curve
to be non-compact as long the associated map $f'$  is proper.

The proof that $C_0$ has a positive cohomology assignment proceeds as
follows:
The point is to assign an integer to each member of a certain subset of
smooth maps from the standard, open unit disk $D\subset\Bbb C$ into $Y$. The
allowed maps are called {admissible  for} $C_0$, and a map $\s$ is
admissible when $\s$ extends as a continuous map from the closure,
 ${\underline D}$,
of $D$ into $Y$ which sends $\p  {\underline D}$ into $Y - C_0$. The
assignment,
$I(\s)$, of an integer to $\s$ must satisfy the following constraints:
\begin{itemize}
\item
$I(\dt)$ must be a homotopy invariant for a homotopy $h$:
$[0,1]\x D\to Y$ which extends as a continuous map from
$[0,1]\x \underline{D}$ into $Y$ sending
$[0,1]\x \p\underline{D}$ into $Y-C_0$.

\item
If $\s$ is admissable and if $\th\co D\to D$ is a proper, degree $k$
map, then $I(\s\th)=kI(\s)$.

\item
If $\s$ is admissable and if $\s^{-1}(C_0)$ is contained in a collection
of disjoint disk $\{D_\nu\}\subset D$, and if each $D_\nu=\th_\nu(D)$
where each $\th_\nu$ is an orientation preserving embedding, then
$I(\s)=\sum_\nu I(\s\dt\th_\nu)$.

\item
If $\s$ is admissable, image$(\s)$ intersects $C_0$ and is a
pseudo-holomorphic embedding, then $I(\s) > 0$.
\hspace{\fill} (3.16)
\end{itemize}

To obtain such an assignment for $C_0$, remark first that each $C_s$ defines a
positive cohomology assignment if the assignment $I(\s)$ for $C_s$ is
 defined to
be the intersection number of $C_s$ with any sufficiently small perturbation of
$\s$ which makes the latter an immersion near $C_s$. (The last point in (3.15)
follows from the fact that pseudo-holomorphic submanifolds have strictly
positive local intersection numbers. See  \cite{Mc}.) This cohomology assignment
for $C_s$ will be denoted below by $I_s(\dt )$ to avoid confusion.

Now, with $\s$ fixed and admissable for $C_0$, then $\s$ will also
be admissable for
$C_s$ if $\s$ is sufficiently large. This follows from the last point of
(3.14).
Furthermore, as $\{C_s\}_{s\in\O''}$ converges to
$C_0$ as a current when $s \to  0$, the values of $I_s(\s)$ for $s
\in\O''$
become independent of $s$ as $s \to  0$.
(The reason for this is as follows: One can first perturb $\s$ so that the
latter is immersed near $C_0$. This will not affect the value of $I_s(\s)$ for
small s by the first point in (3.16).  Then, the value of $I_s(\s)$ can be
computed by integrating a smooth form over $C_s$ which has support in some
small tubular neighborhood of the image of $\s$.)

The fact that  $I_s(\s)$  becomes s independent for small $s$ in $\O''$
allows for the following definition of $I(\s)$ for $C_0: I(\s) =
\lim_{s\in\O''  , s\to 0}  I_s(\s)$ . The fact
that $I(\dt )$ obeys the points in (3.16) follows from the fact that
these
same points are obeyed by each $I_s$. In this regard, the only subtle point in
(3.16) is the final one. And, the final claim in (3.16) follows directly
from Lemma 3.6, below. Indeed, Lemma 3.6 implies that when $s$ is
 large, then $C_s$
intersects the image of a perturbation, $\s'$, of $\s$ which is
homotopic to $\s$ by a
homotopy which is $C_s$ admissible, and which is pseudo-holomorphic near
 $\s^{\,'
-1}(C_s)$. This makes $I_s(\s' ) > 0$ and hence  $I_s(\s)  > 0$ too. (Of
course, the
conclusion that  $I_s(\s)$  is positive follows with apriori knowledge that
$C_s$ 
intersects $\s(D)$ for large $s$.)

\medskip
{\bf Lemma 3.6}\qua{\sl {Let} $\s\co D \to  Y$ {be a pseudo-holomorphic map
and let
$\e > 0$ be given and let $D'\subset  D$ be the concentric disk of
radius
$1 - \e$. Then, there is a neighborhood $U\subset Y$ of $\s(D' )$ and for
each point $y \in U$, there is a smooth map $\s_y\co D \to  Y$ with the
following properties:}

\begin{itemize}

\item ${\mbox{dist}}(\s_y(z), \s(z)) < \e$ \ {\  for all $z
\in D$, and}
$\s_y = \s$ \  {near} \ $\p {\underline D}$,

\item $\s_y|_{D'}$ \  {is pseudo-holomorphic,}

\item $y \in\s_y(D' )$.

\end{itemize}
}

\medskip
The proof of Lemma 3.6 is deferred for a moment to complete the discussion
of (3.15). In particular, remark that with a positive cohomology assignment
for $C_0$ given, Proposition 6.1 in \cite{T3} (with the corrected proof in
\cite{T4})
proves the first point in (3.15).  The argument for the second point in
(3.15) is straightforward, but will be omitted save for the following
outline: First, as noted, a pseudo-holomorphic subvariety has a canonical,
positive cohomology assignment which is defined by taking the intersection
number with maps from the disk which avoid the singular points and have
transversal intersections with the remainder. Second, the current defined by
a pseudo-holomorphic subvariety via integration can be computed from the
values of the corresponding, canonical positive cohomology assignment.
Third, according to Proposition 6.1 in \cite{T3}, the cohomology assignment
given by intersection with an irreducible component of $C_0$ as a
pseudo-holomorphic subvariety is some positive multiple of that given above
for that component as the $s \to  0$ limit of $\{I_s(\dt
)\}_{s\in\O''}$. Fourth, the
limits of $I_s$ were defined using the fact that the set $\O''$  was
 chosen so
that the currents $\{C_s\}_{s\in\O''}$  converge
as $s$ tends to zero.

With the preceding understood, the proof of Proposition 3.4 is completed
with the proof of Lemma 3.6.

\medskip
{\bf Proof of Lemma 3.6}\qua As might be expected, the argument is perturbative
in nature. To start, introduce the bundle $\s^*TX \to  D$. The latter
admits a smooth map, $\phi$, into $X$ which is obtained by
composing the
tautological map to $TX$ with the Riemannian metric's exponential
 map. Thus, a
section $\eta$ of $\s^*TX$ defines, via composition with $\phi$, a
perturbation of
$\s$, with the zero section giving $\s$ itself.  Note that this
composition sends any
neighborhood of the zero section onto a neighborhood of $\s(D)$.

The condition that the resulting map, $4\phi\dt \eta$, is
pseudo-holomorphic defines
an equation for $\eta$ which is elliptic, at least for small $\eta$.
This equation has
the form
\addtocounter{equation}{2}
\begin{equation}
{\cal D}\eta+{\cal R}(\eta,\nabla\eta)=0
\end{equation}
where ${\cal D}$ is an elliptic, first order operator, and where ${\cal R}$
is the remainder, and is small in the sense that $|{\cal R}(\l\eta,
\l\nabla\eta)|
\leq  \z_{\eta}|\l|^2$ for $\l \in\Bbb C$. (See, eg Chapters 1 and 2 of
\cite{MS}.) Regarding ${\cal D}$, note that $J$ endows $\s^*TX$ with a complex
structure, and with this understood, ${\cal D}$ is a zero'th order,
$\R$--linear (not
${\Bbb C}$ linear) perturbation of the standard d-bar operator.

Given the preceding, the Lemma 3.6 follows from

\medskip
{\bf Lemma 3.7}\qua{\sl {There exist constants $\z > 1$ and $\d > 0$ 
 such that
for each point $z \in D'$ and} $\eta_z \in \s^*TX|_z$ with $|\eta_z|
<\d$, {there exists a solution $\eta$ to} (3.17) {over the
concentric disk
with radius $1 - \e/2$ which obeys} $\eta(z) = \eta_z$ {and} $|\eta|
\leq  \z
|\eta_z|$.}

\medskip
{\bf Proof of Lemma 3.7}\qua There are three steps to the proof.

\med
{\bf Step 1}\qua Introduce the disks $D_{-1}, D_0$ and
$D_1$ which are
concentric to $D$ and have respective radii $1 - \e/16$, $1 - \e/8$ and $1
- \e/4$.
Let $\k\co D
\to  [0, 1]$ be a smooth function with compact support on $D_1$ which
equals 1 on the radius $1 - \e/2$ disk. Given $\k$, extend (3.17) to an
equation
over $D_0$:
\begin{equation}
{\cal D}\eta+\kappa {\cal R}(\eta,\nabla\eta)=0
\end{equation}
The point is that the latter equation is linear on $D_0 - D_1$. A
solution to
(3.18) will be found having the form $\eta = \eta^0 + {\cal D}^{\dagger}u$,
where
$\eta^0$ is a smooth section of $\s^*TX$ over $D_0$ which is annihilated
by ${\cal
D}$, and where ${\cal D}^{\dagger}$ is the formal $L^2$ adjoint of ${\cal
D}$. Also,
$u$ is a section of $\s^*TX {\otimes} T^{0,1}D$ which is
constrained to vanish on the
boundary of the closure of $D_0$. With $\eta^0$ fixed and with small norm
(suitably
measured), $u$ can be found as a fixed point of the map ${\cal T}$
 which sends $u$ to
\begin{equation}
{\cal T}(u)=-({\cal D}{\cal D}^{\dagger})^{-1}(\kappa {\cal R}(\eta'+
{\cal D}^{\dagger} u,\nabla(\eta^0 +{\cal D}^{\dagger} u)).
\end{equation}
Here, $({\cal D}{\cal D}^{\dagger})^{-1}$ signifies the Dirichlet Green's
function for
${\cal D}{\cal D}^{\dagger}$. Standard techniques can be used to prove that
(3.19) has a
unique small solution $u$ provided that $\eta^0$ is small in an
appropriate sense. In
particular it is sufficient to measure $\eta^0$ with the norm $\| \
\dt \ \|$ which
is the sum of the $L^2_1$ norm over $D_0$ and the $C^{0,1/2}$ Holder norm
over $D_1$.
Indeed, a simple contraction mapping argument (as in the proof of
 Lemma 5.5 in \cite{T3}
or \cite{T4}) finds $\d > 0$ and $\xi \geq  1$ such that when $\|\eta\|<\d$,
 then (3.19) has a unique solution $u$ which obeys the
constraint $\|u\|+\|\nabla u\|
\leq  \xi^{-1} \d$. Infact, the solution obeys
\begin{equation}
\|u\| +\|\nabla u\|\leq\xi \|\eta^0\|^2 .
\end{equation}
Moreover, the solution $u$ varies smoothly with varying $\eta^0$.
(The estimates
which prove that ${\cal T}$ is a contraction on a suitable ball are
facilitated by the fact that the Green's function $({\cal D}{\cal
D}^{\dagger})^{-1}$ has
the same singular behavior as the Laplacian on the disk. An essentially
identical contraction mapping argument with similar estimates can be found
in the proof of Lemma 5.5 in \cite{T3} and \cite{T4}.)

With $u$ and (3.20) understood, Lemma 3.7 follows from the assertion that
there exists $\z \geq  1$ such that for each $z \in D'$  and
each $\eta_z \in\s^*TX|_z$, there is a section $\eta^0$ of $\s^*TX$
over $D_0$ which
restricts to $z$ as $\eta_z$, is annihilated by ${\cal D}$, and obeys
$\|\eta^0\|
\leq  \z |\eta_z|$. The next two steps prove this assertion.

\medskip
{\bf Step 2}\qua The first part of the proof of the
preceding assertion
gives an affirmative answer to the following question: Given $z \in D'$
 and a point $\eta_z \in \s^*TX|_z$, is there an element $\eta^0$ in the
kernel of ${\cal D}$ over $D_0$ which obeys $\eta^0(z) = \eta_z$?

To find such a section $\eta^0$, consider fixing a point $y \in D_{-1} -
D_0$ and a
point $\eta_y$ in $\s^*TX|_y {\otimes}T^{0,1}D|_y$. Then, let $v$ denote the
section of
$\s^*TX$ over $D_0$ which is given by ${\cal D}^{\dagger}({\cal D}{\cal
D}^{\dagger})^{-1}
( \dt  , y)\eta_y$, where, in this step, $({\cal D}{\cal
D}^{\dagger})^{-1}$ denote the Greens with
Dirichlet boundary conditions on the boundary of the closure of $D_{-1}$.
Note that
${\cal D} v = 0$ on $D_0$.

Now, consider varying $y$ over a small ball in $D_{-1} - D_0$ and
varying the
corresponding $\eta_y$. Restricting the corresponding $v$'s to the point
$z$ generates
a vector subspace $V \subset\s^*TX|_z$. The claim is that $V = \s^*TX|_z$.
Accept this
claim and it follows that there is an element in the kernel of ${\cal D}$
over $D_0$ which restricts to $z$ as $\eta_z$.

Here is the proof of the claim: Were the claim false, there would exist a
non-zero element $\eta'  \in\s^*TX|_z$ with the property that
$\eta^{'\dagger} ({\cal D}{\cal D}^{\dagger})^{-1}(z, y) = 0$ at all
points $y$ from
some nonempty, open subset of $D_{-1} - D_0$. Since $z \in D_1$,
the section
$\eta^{'\dagger} ({\cal D}{\cal D}^{\dagger})^{-1}(z, \dt  )$ of
 $\s^*TX {\otimes}
T^{0,1}D$ is annihilated by ${\cal D}^{\dagger}$ over $D_{-1} - D_0$. If
it vanishes on
a non-empty open set, Aronszajn's unique continuation theorem \cite{Ar}
forces it
to vanish everywhere. Such is not the case. (Consider a point very close to
z where the small distance asymptotics of $({\cal D}{\cal D}^{\dagger})^{-1}(z,
\dt  )$
can be invoked.)

\medskip
{\bf Step 3}\qua The previous step established the
following: Given $z
\in D'$  there exists a set $\{(z_j, \eta_j)\}_{1\leq j\leq 4}$
of four pairs consisting of points $\{z_j\}\subset  D_{-1} - D_0$ and
for each $z_j$,
a corresponding non-zero element $\eta_j$ in the fiber of $\s^*TX \otimes
T^{0,1}D$ at
$z_j$. This set has the following property:  The set $\{{\cal
D}^{\dagger}
({\cal D}{\cal D}^t)^{-1}(z, z_j)\eta_j\}_{1\leq j\leq 4}$ of vectors in $\s^*TX|_z$ spans
 $\s^*TX|_z$.
Now, because the spanning condition is an open condition, one can
additionally conclude that $\{{\cal D}^{\dagger}
({\cal D}{\cal D}^t)^{-1}(z', z_j)\eta_j\}_{1\leq j\leq 4}$ also spans
$\s^*TX|_{z'}$  for
all points $z'$  in some open neighborhood of $D$.

With the preceding understood, the fact that the closure of $D'$
is compact implies that there exists a finite $N (\geq  4)$ and
 a set $\{
(z_j, \eta_j)\}_{1\leq j\leq N}$ with $\{z_j\}\subset  D_{-1} - D_0$ and
with $\eta_j$
in the fiber of $\s^*TX \otimes  T^{0,1}D$ at $z_j$ and such that
$\{{\cal D}^{\dagger}
({\cal D}{\cal D}^t)^{-1}(z, z_j)\eta_j\}_{1\leq j\leq N}$ spans $\s^*TX|_z$ at
each point
of $D'$.

The fact that $D'$  has compact closure in $D_0$, and the fact that
the integer $N$ here is finite imply the following: There exists $\z
\geq  1$
such that for each point $z \in D'$  and $\eta_z \in \s^*TX|_z$,
there exists a section, $\eta^0$, of $\s^*TX$ over $D_0$ which obeys
 $\eta^0(z) =
\eta_z$, ${\cal D}\eta^0 = 0$ and $\|\eta\|\leq\z|\eta_z|$. Given this last
remark, the conclusions of Step 1 above can be brought to bear to
finish the proof of
Lemma 3.7.

Before ending this section, it is pertinent to remark here that the proof of
Proposition 3.3 just offered proves, with minimal modifications, a somewhat
stronger local compactness theorem.  This stronger version of Proposition
3.3 is not necessary for this article, but is stated below for possible
future reference:

\medskip
{\bf Proposition 3.8}\qua {\sl Let $Y$ be a smooth 4--manifold with symplectic
form $\o$
and a compatible almost complex structure $J$. Let $\{Y_s\}_{s\to 0}$ be
a countable exhaustion of $Y$ by a nested sequence of open sets (so
$Y_s\subset
Y_{s'}$  with compact closure when $s > s'$). In
addition, suppose that $\{(\o_s, J_s)\}_{s\to 0}$ is a corresponding
sequence whose typical element is a pair consisting of a symplectic form,
$\o_s$, on $Y_s$, and a compatible almost complex structure $J_s$. Require that
$\{\o_s,
J_s)\}_{s\to 0}$ converges to $(\o, J)$ in the $C^{\i}$
topology on compact subsets of $Y$. Next, let $E \geq  1$ be given with a
sequence $\{C_s\subset  Y_s\}_{s\to 0}$ where $C_s$ is a finite energy,
$J_s$--pseudo-holomorphic subvariety with $\int_{C_s}\o_{s} \leq E$.
Then there is a $J$--pseudo-holomorphic subvariety $C_0\subset Y$ and a
subsequence of $\{C_s\}$ (hence renumbered consecutively from 1) with the
following properties:

\begin{itemize}
\item {The integral of $\o$ over $C_0$ is bounded also by}
$E$.

\item $\{C_s\}$ {converges geometrically to} $C$.
\end{itemize}}

\medskip
{\bf Proof of Proposition 3.8}\qua The only essential modification of the proof
of Proposition 3.3 occurs in (3.18) and (3.19) with the addition of an
$s$--dependent term ${\cal R}_s$ to ${\cal R}$ which is $\R$--linear in $\eta$
and accounts for the difference between $J$ and $J_s$. However, as $J - J_s$
converges to zero on compact subsets of $Y$, this term has no substantive
affect on the subsequent arguments.

\section{The multiplicity of the limit}
\setcounter{equation}{0}

Fix $t_0 \in Z$. Then the constructions of the previous section find a
convergent subsequence in the ball $B(1)$ from the sequence
$\{C_s\}_{s>0}$ of
rescaled versions of the original subvariety $C$. As currents, the
subsequence
of $\{C_s\}_{s>0}$ converges to a current in $B(1)$. This current can be
written as a sum of currents with integer multiplicities, where each term
corresponds to one of the cases from Examples 1.3, 1.4 and 1.5. The
following proposition makes a somewhat stronger statement; for it asserts,
in part, that all convergent subsequences of $\{C_s\}_{s>0}$ have the same
limit current.

\medskip
{\bf Proposition 4.1}\qua{\sl {Let $C$ be a finite energy, pseudo-holomorphic
subvariety. Fix a point $t_0 \in Z$. Let B denote the ball of radius
1 in $S^1\x\R^3$
 with center $t_0$. For $s > 0$, but small, introduce the sequence
$\{C_s\}$ of pseudo-holomorphic subvarieties of $B$ as defined in (3.1).
Then, as a sequence of currents, every limit of $\{C_s\}_{s>0}$ has the form}
\begin{equation}
\sum_{1\leq i\leq p} C^{(i)}_\rho +q_+\dt
C_{z+} +q_-\dt C_{z-} + \sum_{1\leq i\leq n_+} C^{(i)}_{h+} +
\sum_{1\leq i\leq n_-} C^{(i)}_{h-} .
\end{equation}
{Here, $\{C^{(i)}_{\rho}\}$ is a finite set of (not necessarily distinct)
pseudo-holomorphic submanifolds which are described in Example 1.3 (with $\nu$
determined by the index $i$); $C_{z\pm}$  are described by Example 1.4 (the
$\pm$
signifies when $\nu = (0, 0, \pm 1))$ while $q_{\pm}$  are
non-negative integers;
and $\{C^{(i)}_{h\pm}\}$ is a finite set of (not necessarily distinct)
pseudo-holomorphic submanifolds with are described in Example 1.5 where the
constant c is determined by the index $i$ while the $\pm$  correlates
with the
signs used in (1.7). Furthermore, the integers $p, q_+, q_-, n_+, n_-$ are
identical for any limit; and the corresponding sets $\{c^{(i)}\}\subset  (0,
2/3^{3/2}]$ which determine (in part) the Example 1.5 data
$\{C^{(i)}_{h\pm}\}$
are identical for any limit up to permutations.}}

\medskip
This proposition implies that the limiting data in (4.1) is unique except
possibly up to separate rotations $(\phi\to\phi+{\rm constant}$) of each
element in $\{C^{(i)}_{\rho}\}$ and $\{C^{(i)}_{h\pm}\}$.

The remainder of this section is occupied with the proof of Proposition 4.1.

\medskip
{\bf Proof of Proposition 4.1}\qua The fact that all possible limit currents
are given by a formal sum of the form in (4.1) (with all coefficients
non-negative) is a restatement of Proposition 3.2. At issue here is the
assertion that there is a unique limit.  The proof of this assertion is
broken into four steps.

{\bf Step  1}\qua This step verifies the assertion
that $p - q$ is
uniquely defined and independent of $t_0$. Indeed, this assertion
follows from the
observation that $p - q$ is the intersection number between $C$ and a
2--sphere
of the form $t =$ constant, $\rho^2 + z^2 = \d > 0$.  This can be
seen from the
following argument: Let $S$ be the sphere $t = t_0$ and $\rho^2 + z^2
= {1/2}$. This
sphere has trivial normal bundle, so one can identify a tubular neighborhood
of $S$ with the product $S \x  D$, where $D$ is small radius disk in
the plane.
(Here, the sphere $S$ is identified with $S \x  0\subset  S \x  D$.)
 Next, fix a
compactly supported, non-negative 2--form on $D$ with integral equal to 1.
 This
form defines a closed 2--form, $\mu$, on $S \x  D$ which pairs with the
 current in
(4.1) to give $p - q$. (All of the cases of Example 1.5 have
intersection
number zero with a $t =$ constant 2--sphere.) Thus, with $\e >
0$ fixed, and $s$
large, this $\mu$ pairs with $C_s$ from the subsequence to give a number,
say $y_s$,
which is within $\e$ of $p - q$.

However, $\mu$ is closed, and thus $y_s$ is equal to the pairing between
$C_s$ and any
2--form $\mu'=\mu + d\nu$, where $\nu$ is a 1--form with compact support
in $S\x D$. In particular, fix $s$ large, and perturb $S$ slightly (if
necessary) so
that its intersection with $C_s$ is transversal. Then one can find a 2--form
$\mu'$ with
support in any tubular neighborhood of this perturbation (and thus in $S\x
D$) whose
integral over
$C_s$ is equal to the intersection number of $C_s$ with the perturbed
version of $S$.

With the preceding understood, it follows that this intersection number is
within $\e$ of $p - q$, and thus equal to $p - q$ as both numbers are
integers.
Finally, rescaling shows that this intersection is the same as that between
$C$ and a perturbation of $t = t_0$, and $\rho^2 + z^2 = s^2/2$
sphere. Finally,
remark that the intersection between $C$ and such a sphere is the
same as that
between $C$ and any tiny perturbation of the sphere $t = t_0$, $\rho^2 +
z^2 = 1/2$,
as these two spheres are homologous in the complement of $Z$. In
particular, this
shows that the difference, $p - q$, is independent of the chosen
subsequence
of $\{C_s\}$.

The same invariance under perturbation argument shows that $p - q$ is also
independent of $t_0$.

\medskip
{\bf Step  2}\qua  Let $\O\subset  (0, 1]$ be a
countable, decreasing
set with limit zero and such that the sequence $\{C_s\}_{s\in\O}$ converges as
described in Proposition 3.2. The limit in (4.1)  determines, among other
things, the
integers $p$ and $q_{\pm}$ . This step identifies these integers as
intersection
numbers with the following claim: There exists $d_0 > 0$ such that
when $\d$ is
positive and less than $\d_0$, then:
\begin{itemize}
\item $p$ is the intersection number in $B(1)$
between $C_s$ for small $s$ and the pseudo-holomorphic cylinder
$t=t_0$, \ $f=1/100$ and $|z| <\d$.
\item $q_{\pm}$ is the intersection number in $B(1)$
between $C_s$ for small $s$ and the pseudo-holomorphic disk
$t=t_0$, \ $f=-1/100$, $\rho <\d$ and $\pm z >0$.
\hspace{\fill} (4.2)
\end{itemize} 

The remainder of this step offers a proof of the claim for $q_+$.
The proofs
for $p$ and $q_-$ are completely analogous and will be omitted.

To begin the argument, consider that $\rho$ is bounded from below on a
case from
Example 1.5 by $s\dt c$. Thus, since only finitely many of the
Example 1.5
cases appear in (4.1), there exists $\d_0 > 0$ such that each
$C^{(i)}_{h+}$ which
appears in (4.1) has empty intersection in $B(1)$ with the disk $t
 = t_0$, $f =
-1/100$, $\rho < \d_0$ and $z > 0$. With the preceding understood,
take $\d < \d_0/4$
and let $\Sig$ denote the $t = t_0$, $f = -1/100$, $\rho < \d$ and $z > 0$
disk. Let
$D\subset {\Bbb C}$ denote the standard disk, and fix an
 identification of a tubular
neighborhood of $\Sig$ in $S^1\x B^3$ with $\Sig\x D$ so that $\Sig$
in the former is
identified with $\Sig\x\{0\}$ in the latter. Also, arrange the
identification so that the functions $\phi$ and $h$ are constant on each
$\{$point$\}\x D\subset\Sig\x D$. Finally, require $f < 0$ on this tubular
neighborhood. Now, given $\e > 0$, let $\mu_\e$ denote a smooth 2--form on $D$
with integral
1 and support in the concentric disk $D_\e \subset D$ of radius
$\e$. Modify $\mu_\e$
(via multiplication by an $\e$--independent, compactly supported function on
$\Sig$) so
that the result, $\mu'_\e$, has compact support in $S^1  \x  (B^3 -
\{0\}$) and also
agrees with $\mu_\e$ on some $\e$--independent neighborhood the $h =
 0$ disk in
$\Sig\x  D$. This last condition insures that $\mu'_\e$  is closed near the
$h = 0$
disk in $\Sig\x D$.

By definition, the current in (4.1) pairs with  $\mu'_\e$  to give the
integer $q_+$. It then follows that the current $C_s$ for $s \in\O$ and
small will
pair with each  $\mu'_\e$   to give $q_+$ also. Since $\e$ can be
arbitrarily
small here (though positive), and  $\mu'_\e$   is closed near the $h = 0$
disk in $\Sig\x  D$, the perturbation argument from Step 1 can be
applied here to
prove that $C_s$ has intersection number $q_+$ with $\Sig$ for small $s$.

\medskip
{\bf Step  3}\qua
This step argues that $p, q_+$ and $q_-$ are the same for any
limit of the sequence $\{C_s\}$. Again, the details will be given for the
case of $q_+$ with the analogous arguments for $q_-$ and $p$ left to
the reader. (In
fact, once $q_{\pm}$  are known to be unique, then the uniqueness of
$p$ follows
from Step 1.)

To begin the discussion, suppose, for argument's sake, that there exist
decreasing sequences $\O,\O'\subset [0, 1]$, both having limit zero and
such that the corresponding sequences $\{C_s\}_{s\in\O}$ and
$\{C_s\}_{s\in\O'}$
converge geometrically to limits given by (4.1) with
the corresponding $q_+$ and $q'_+$  not equal. A contradiction will
follow from this assumption.

To obtain the contradiction, reintroduce, for $\d > 0$ and very
small, the
pseudo-holomorphic disk $\Sig$ where $t = t_0$, $f = -1/100$, $\rho
< \d$ and $z >
0$. Then, use (4.2) to conclude that when $\d$ is small, then
$C_s$ for all small $s
\in\O$ has intersection $q_+$ with this $\Sig$. Likewise, one can choose
$\d$ still
smaller if necessary so that $C_s$ for all small $s'  \in\O$ has
intersection number $q'_+$  with $\Sig$. Moreover, for sufficiently
small $s$, these intersection points in either case will lie where $\rho$
is less
than any given positive constant. In particular, one can find a small number
$c > 0$ such that the following is true:
\begin{itemize}
\item The period of the corresponding solution to (1.6)
is an irrational multiple of $\pi$.
\item The level set $h/(\rho^2+z^2)^{3/2}=c$
intersects $\Sigma$ as a circle.
\item For sufficiently small $s$ in either $\O$ or $\O'$,
all intersection points of $C_s$ with $\Sigma$ lie where
$h/(\rho^2+z^2)^{3/2}<c$.\hspace{\fill} (4.3)
\end{itemize} 
(See Lemma 1.9 about the first point.)

Now remember that $C_s$ is obtained from $s$ by rescaling, so
that the
intersections between $C_s$ and $\Sig$ are obtained by rescaling those
between $C$ and
$\Sig_{1 /s}$, where $\Sig_{1 /s}$ has $t = t_0$, $f = -s^2/100$, $\rho <
 s\d$ and $z >
0$.

With the preceding understood, consider the following construction:  To any
given $s \in\O$, let $s' (s) \in\O'$ denote the
largest element which is less than $s$. For $\s \in [s' (s), s]$,
consider the intersection of $C$ with $\Sig_{1 /s}$. For $\s$ near $s$,
 the
intersection number is $q_+$ and all intersection points have $h/(\rho^2 +
z^2)^{3/2}
< c$. Likewise, for $\s$ near $s'$, the corresponding intersection
number
is $q'_+$  and all intersection points also have $h/(\rho^2 + z^2)^{3/2}
<c$.  Since
intersection number is a deformation invariant (for compactly supported
homotopies),
it follows that there exists $s'' \in (s' (s), s)$ where $C$ has an
intersection
with $\Sig_{1 /s''}$ which lies where $h/(\rho^2 + z^2)^{3/2} = c$.

This last conclusion implies that there is a sequence
$\O''\subset [0, 1]$ which decreases to zero and which has the property
that each $C_s$ for $s \in\O''$  intersects the
 $h/(\rho^2 + z^2)^{3/2} =c$ set. According to Proposition 3.2, this
sequence has a
convergent subsequence with a geometric limit which is a union of cases from
Examples 1.3, 1.4 and 1.5.  This limit must include a case from Example 1.5
which is determined by the given value of $c$. However, as this value of $c$
corresponds to a solution to (1.6) with whose period is not in $\pi\Bbb Q$,
the corresponding pseudo-holomorphic submanifold will not have
finite energy. This is the contradiction which establishes that $q_+ =q'_+$.

\medskip
{\bf Step  4}\qua
 Each limit as in (4.1) determines a set $\{C^{(i)}_{h\pm}\}$ of
pseudo-holomorphic
submanifolds which are described by Example 1.5. The latter determines a
corresponding set of constants $\{c^{(i\pm)}\}\subset (0, 2/3^{3/2}]$.
(Note that this
set need not consist of distinct elements.) The purpose of this step is to
establish
that the set $\{c^{(i\pm)}\}$ does not depend on chosen subsequence. This
is to say
that all limits of $\{C_s\}$ determine the same set $\{c^{(i\pm)}\}$ up to
permutations.

What follows is the argument for $\{c^{(i+)}\}$. The argument for
$\{c^{(i-)}\}$
is completely analogous and will be omitted. To start the argument
suppose that $c \in (0, 2/3^{3/2}]$ is such that the corresponding
 solution to
(1.6) has period which is a rational multiple of $\pi$ . This $c$ determines
a least positive integer, $b(c)$, for which the corresponding solution to
(1.6) is periodic on $[0, 2\pi b]$. Extend this assignment $b(\dt )$ to a
function on $(0, 2/3^{3/2}]$ by declaring that $b(c) = \i$  when the
period in question is an irrational multiple of $\pi$ . With the preceding
understood, the assertion that the set $\{c^{(i+)}\}$ is limit
independent is
equivalent to the assertion that
$\sum_{c^{(i+)} >c} b(c^{(i+)})$ is limit independent for each $c \in (0,
2/3^{3/2}]$ with $b(c) =\i$ .

To establish the latter claim, fix $c$ with
$b(c) = \i$  and argue as in Steps 1 and 2 above to prove that
\begin{itemize}
\item $\sum_{c^{(i+)} >c} b(c^{(i+)})$ is the intersection number in the
 ball $B(1)$
between $C_s$ for $s \in\O$
sufficiently small and the pseudo-holomorphic submanifold where $\phi = 0$
and $h = c$. \hspace{\fill} (4.4)
\end{itemize} 
Then, by rescaling, (4.4) is the same as the assertion that
$\sum_{c^{(i+)}  >c} b(c^{(i+)})$ is the intersection number in the
 ball $B(s)$
(for $s \in\O$) between $C$ and the pseudo-holomorphic submanifold where
$\phi = 0$
and $h = c s^3$.

Given the preceding, suppose that $\O$ and $\O'$  are a pair of
countable, decreasing subsets of [0, 1] with limit zero and such that the
corresponding sequences $\{C_s\}_{s\in\O}$ and $\{C_s\}_{s\in\O'}$
converge as described in (4.1) and to yield the data
$\{c^{(i+)}\}$ and $\{{c'}^{(i+)}\}$, respectively. Suppose that the
sum $\sum_{c^{(i+)}  >c} b(c^{(i+)})$ differs from its primed analog. This
supposition leads to a contradiction as follows: By an argument which is
completely analogous to that used in Step 3 above, one can construct, for
each $s \in\O$, a number $s''\equiv s''(s)\in (0,s)$
with the property that $C$
intersects the pseudo-holomorphic submanifold where $\phi = 0$ and $h
 = c
s^{'' 3}$ on the boundary of the ball
$B(s'')$. The corresponding sequence
$\O''\equiv \{s''(s):s\in\O\}$
would then converge geometrically as described in
Proposition 3.2 to a limit which is given by (4.1). However, for this limit,
the corresponding set $\{{c''}^{(i+)}_h\}$ would
contain an element which has $h = c$ on the boundary of the ball $B(1)$. The
corresponding $c \in (0, 2/3^{3/2}]$ for such an element has $b(c) = \i$
and thus, the corresponding solutions to (1.6) do not have periods in
$\pi\Bbb Q$. As
these values for $c$ are not allowed in Example 1.5 (since
 they yield infinite energy
pseudo-holomorphic varieties), the desired contradiction has been obtained.

\section{The integers $n_{\pm}$ }
\setcounter{equation}{0}

The previous section associates to each $t_0 \in Z$ the ordered set $(p, q_+,
q_-, n_+, n_-$) of non-negative integers.  The purpose of this
section is to
prove:

\medskip
{\bf Proposition 5.1}\qua{\sl {There are only finitely many} $t_0 \in Z$ {\it
where}
$n =n_+ + n_- \neq   0$.}

\medskip
This proposition is a corollary to a second proposition of some independent
interest. The statement of this second proposition requires a brief
digression. To start the digression, fix $\d\in [1/16, 1/8]$ and, for small
$s > 0$, consider the compact set $\O(s)$ in $B^3$ where
\begin{equation}
|h|\leq s^3, \quad \rho\leq\sqrt{\d} \quad {\mbox{and}}\quad
|z|\leq\d.
\end{equation}
Thus, the boundary of $\O(s)$ consists of the parts:
\begin{enumerate}
\item $h=\pm s^3$ \ and \ $\rho\leq\sqrt{\d}$ \ and \
$|z|\leq\d$
\item $\rho=\sqrt{\d}$ \ and \ $|z|\leq s^3/\d$
\item $|z|=\d$ and  $\rho\leq s^{3/2}/\sqrt{\d}$ \itemnum{5.2}
\end{enumerate}
Use $C(s)$ to denote $C \cap  (S^1 \x\O(s))$. With $\d$ fixed, introduce the
following function of $s$:
\begin{equation}
\mu(s)=\int_{C(s)}d\phi\wedge dh
\end{equation}
Now, end the digression and consider the following proposition:

\medskip
{\bf Proposition 5.2}\qua {\sl There is a constant $\z > 0$ such
 that $\mu(s) \leq \z
s^3$ for all $s$ sufficiently small.}

\medskip
Note that this proposition has corollary which bounds the integral of
$d\phi\wedge dh$
over the set $C' (r) = C \cap  (S^1 \x\{\rho^2 + z^2 \leq  r\})$
by an $r$--dependent multiple of $r^3$. That is, there exists $\z > 0$
 which is
independent of $r$ and is such that
\begin{equation}
\int_{C'(r)}d\phi\wedge dh \leq \zeta r^3
\end{equation}
for all $r < 1$. (This bound follows from Proposition 5.2 because the ball
of radius $r$ is contained in $\O(s = 2^{1/2}/3^{1/2} r)$.)

\medskip
{\bf Proof of Proposition 5.1 assuming Proposition 5.2}\qua The first point to
make is that the form  $d\phi\wedge dh$ restricts to $C$ as a non-negative
multiple of
$\o$, and thus $\mu(s)$ is no smaller than the integral of  $d\phi\wedge dh$
 over any
subset of $C$ in $C(s)$. With this understood, the plan is bound
$\mu(s)$ from below
by considering the restriction of the integral in question to the
intersection
of $C$ with small balls about the points in $Z_0$ where $n \neq   0$.

However, for this purpose, it is necessary to digress with a return to
Example 1.5. To start the digression, suppose that $\Sig$ is given by
Example
1.5. Let $b$ denote the positive integer supplied in the example. (This
integer is the smallest $b'  \in\{1, 2, \dots\}$ such that
the solutions to (1.6) are periodic on $[0, 2\pi b' ]$.) With
$s > 0$ given, note that set of points on $\Sig$ where $|h| \leq  s^3$ is
contained in the ball of radius $s/c^{1/3}$ and center at 0.  (Here, $c$
is the
constant in $(0, 2/3^{3/2}]$ which determines $\Sig$ up to
rotations in the $x-y$
plane.) In particular, the set where $|h| \leq  s^3$ on $\Sig$ is
 compact and
Stoke's theorem finds
\begin{equation}
\int_{\Sigma\cap \{|h|\leq s\}}d\phi\wedge dh=
2\pi bs^3\geq 2\pi s^3 .
\end{equation}

End the digression. With the preceding understood, suppose that the set of
points in $Z$ where $n \neq  0$ contains at least $N$ points, where $N$
is some given
integer.  Then, it follows from (5.3) and Propositions 3.2 and 4.1 that $\mu(s)
\geq 2\pi  N s^3$ when $s$ is small. This last lower bound for
$\mu$ plus
Proposition 5.2 bounds $N$.

\medskip
{\bf Proof of Proposition 5.2}\qua  The existence of a bound $s^{-3}\dt \mu(s)
\leq
\z$ is argued as follows: First of all, one can remove from $C$ all
 irreducible
components which are described by cases in Examples 1.3 and 1.4 since these
make no contribution to $\mu(s)$.

Next, note that the intersection of $C$ with the surfaces $\rho =
\sqrt{\d}$ and $|z|
 = \d$ defines a smooth curve for all but finitely many $\d$ in
$[1/16, 1/8]$.
This follows from Sard's theorem. Likewise, the intersection of $C$ with the
surface $|h| = s^3$ is a smooth curve for all but countably many
values of
$s$ in $(0, 1]$. Thus, with $\d$ chosen appropriately, and with $s$ chosen
from the
complement of a countable set, the number $\mu(s)$ can be computed via Stokes'
theorem as
\begin{eqnarray}
\mu(s)=-s^3\dt\int_{C(s)\cap \{h=s\}}d\phi+
s^3\dt&&\kern-1.5em\int_{C'(s)\cap \{h=-s\}}d\phi \\
&&-\int_{C'(s)\cap \{\rho=\sqrt\d\}}h\dt d\phi-
\int_{C'(s)\cap \{|z|=\d\}} h\dt d\phi . \nonumber 
\end{eqnarray}

(The integration by parts formula here can be justified by employing the
functions $\{\eta_\e\}$ as in (2.4). This step is straightforward and
is left
to the reader.)

Consider now the various terms in (5.6). As $C$ does not coincide
 with the
product of $S^1$  and either the positive or negative $z$--axis, then for
sufficiently small $s$, the right most term in (5.6) will be zero for
all but
finitely many choices of $\d$.

To understand the contribution to $\mu(s)$ from the term which is
second from
the right in (5.6), remark first that $|d\phi|=\rho^{-1}$, so the
restriction of
$d\phi$ to $C \cap \{\rho = \d^{1/2}, |z| \leq  s/\d\}$ has norm less than or
equal to $\d^{-1/2}$. Meanwhile, on this part of $C$, one has $h =
z\dt \d$. Finally,
this part of $C$ is uniformly far from $Z$ (since $\d \geq
 1/4$), and so the
length of this part of $C$ is bounded independently of $s$. With the
preceding
understood, it follows that the term which is second from the right in (5.6)
has absolute value no greater than $\z\dt s^3$, where $\z$ is a constant
which
depends on $C$, but not on $s$.

To analyze the left most term in (5.6), consider completing the integration
path $\nu = C \cap \{h = s^3, \ \rho \leq  \sqrt{\d}, \ z \leq
 \d\}$ to make a
closed path in $C$ which avoids where $\rho= 0$ and which stays
uniformly away
from the set $Z$. For example, one could add to $\nu$ the following:
\begin{equation}
\nu'= C\cap\{\rho^4+z^2=\d^2+s^3/\d^2, \ s^{3/2}/\d\leq\rho, \
s^3/\d\leq z\} 
\end{equation}
or any slight perturbation thereof. For generic $s$, this last is a smooth
curve in $C$ whose distance from $Z$ is at least $\sqrt{\d}$. Thus, the
integral of
$d\phi$ over $\nu'$  is uniformly bounded, independent of $s$.
Furthermore, the curve $\nu'$ is homologous in $C - C \cap  \rho^{-1}(0)$
to $\nu$ since both are obtained by slicing $C$ with functions which
are defined
on $B^3 - \rho^{-1}(0)$. (Consider restricting the family of
functions $\{w_t =
t\dt (\rho^4 + z^2 - \d^2 - s^6/\d^2) + (1 - t)\dt (\rho^2z -
s^3)\}_{t\in [0,1]}$ to
the set to $C \cap \{s^{3/2}/\sqrt{\d} \leq  \rho, \  s^3/\d \leq  z\}$.
Then, the
family $\{w_t^{-1}(0)\}_{t\in [0,1]}$ defines a homology rel endpoints
between $\nu$
and $\nu'$  in $C - C \cap  \rho^{-1}(0)$.)

Since $\nu$ and $\nu'$  are homologous rel end points in $C \cap
\rho^{-1}(0)$, the integral of $d\phi$ over $\nu$ is also uniformly
bounded,
independent of $s$. Thus, the left-most term in (5.6) has absolute value
 no greater
than $\z\dt s^3$, where $\z$ depends on $C$, but not on $s$ when the
latter is small.

The term which is second to the left in (5.6) is handled in the same way as
the previous term.

It follows from the preceding that $\mu(s) \leq \z\dt s^3$, where
 $\z$ depends on
$C$, but is independent of $s$, as was required.  Note that this
estimate is
independent of the choice of $\d$ as well since with s fixed, $\mu(s)$ is a
continuous, increasing function of $\d$.

\section{Semi-continuity of $p$ and $q$ as functions on~$Z$}
\setcounter{equation}{0}

For each point $t_0 \in Z$, Proposition 4.1 assigns the ordered pair of
non-negative integers $(p, q = q_+ + q_-)$. That is, $(p, q)$
 from Proposition 4.1
can be thought of as a pair of functions on $Z$. The purpose of this section
is to investigate $(p, q)$ as functions on $Z$.

Before starting, remark that there are only a finite number of points in $Z$
where the integer $n$ in (4.1) is non-zero. It follows from the
convergence
properties as described in Proposition 3.2 that there exists $r_0 >
0$ such
that $C \cap \{\rho^2 + z^2 \leq r_0^2\}$ can be decomposed
as the disjoint
union $C'  \cup C''$  of
pseudo-holomorphic subvarieties with the following properties: All points in
$Z$ have $n \equiv n_++n_- = 0$ when the latter is computed using
$C'$ instead of $C$. Conversely, all points in $Z$ have $p = q = 0$
when these integers are computed using $C''$
instead of $C$. With this understood, it follows that no generality
is lost by
considering solely the case where $C = C'$. This condition will
henceforth be assumed.

In addition, one can assume with out loss of generality that $C \neq   C_1 \cup
C_2$ where $C_2$ is described in either Example 1.3 or 1.4.

Consider now a sequence $\{t_i\}\subset  Z$ which converges to some $t_0
\in Z$. As
there is a uniform bound on the values for $p$ and $q$ over the whole
of $Z$ (from
Proposition 2.1), one can assume, without loss of generality, that $p$ and $q$
are constant on the sequence $\{t_i\}$. Use $(p_1, q_1)$ to denote this
constant value. The goal in this section is to prove the following 
proposition:

\medskip
{\bf Proposition 6.1}\qua {\sl {Define $(p_0, q_0)$ and $(p_1, q_1)$ as
above.  Then}
$p_0 \geq p_1$ {and} $q_0 \geq  q_1$.}

\medskip
(Remember from the proof of Proposition 4.1 that $p_0 - q_0 = p_1 - q_1$.)

The proof of Proposition 6.1 depends on Lemma 6.2, below. The statement of
the latter requires the introduction, for $t'  \in Z$ and $s > 0$,
of the pseudo-holomorphic submanifold $\Sig(t' , s)\subset  X$ which is
defined by the equations $t = t'$  and $f = s^2$.

\medskip
{\bf Lemma 6.2}\qua {\sl {There exists $\d, \e > 0$ such that the following is
true: If
$s > 0$ but small, and if $t'  \in Z$ obeys
$|t_0-t'|<\e$, then:}
\begin{enumerate}
\item
{The intersection of $\Sig(t' , s)$ with $C \cap  B(\d)$ is contained
in} $C \cap  B(\d/2)$.
\item {The integer $p_0$ is equal to the algebraic intersection number
 between
$C \cap B(\d)$
and any sufficiently small perturbation of} $\Sig(t' , s)$.
\end{enumerate}}

\medskip
Remark that when $C$ and $\Sig(t' , s)$ are transverse, the algebraic
and geometric intersection numbers agree since both $C$ and $\Sig$ are
pseudo-holomorphic.

\medskip
{\bf Proof of Proposition 6.1}\qua Accept Lemma 6.2 for the moment.  Fix the
index $i$ large, and consider the interesection number between $C
\cap  B(\d)$
and $S(t_i, s)$ as $s$ tends to zero. These intersections will
still be in the
ball $B(\d/2)$, so the algebraic intersection number here is still
 $p_0$. To be
precise, if $\Sig(t_i, s)$ is not transversal to $C$ in $B(\d)$, then make
any arbitrarily small perturbation of the former to make a submanifold
$\Sig'$  which is transversal to $C$ in $B(\d)$. Then, the intersection
number of $\Sig'$  with $C \cap  B(\d)$ will equal $p_0$.
Furthermore,
using the analysis of Lemma 5.5 in \cite{T3}, one can find such small
perturbations for which $\Sig'$  is pseudo-holomorpic near its
intersection points with $C$. In this case, $\Sig'$  will intersect $C$
in precisely $p_0$ points because all points of intersection must have
positive
local sign.

Now, remark that for $s$ tiny, the algebraic intersection number of
$\Sig(t_i, s)$
with $C$ in a very small ball about $t_i$ must equal $p_1$. Indeed, this
is simply
Lemma 6.2 but with $t_i$ replacing $t_0$. In fact, since $C$ is
pseudo-holomorphic,
for s sufficiently small, any perturbation of $\Sig(t_i, s)$ which is
pseudo-holomorphic and transversal to $C$ in $B(\d)$ and which
is sufficiently
close to $\Sig(t_j, s)$ will intersect $C$ in exactly $p_1$ points
 (all with positive
sign) in some very small ball about $t_j$. Since all other points of
intersection with $C$ in $B(\d)$ will also have positive sign, this forces the
inequality $p_1 \leq  p_0$.

\medskip
{\bf Proof of Lemma 6.2}\qua Construct the sequence $\{C_s\}$ for the point
$t_0$. For all $s$ sufficiently small, the current $C_s$ will be close to the
current in (4.1) with $n_+ = n_- = 0$, $p = p_0$ and $q_+ + q_- =
 q_0$. Furthermore,
all points in $C_s$ will be very close to either the $z$--axis or the
$x-y$ plane. Fix
$t'  \in Z$ with $|t_0 - t'| < 1/100$. Thus,
for any value of $\e \in (0, 1/100)$, all intersection points of $C_s$ with
 the
submanifold $\Sig(t' , \e)$ are in $B(1/2)$ when $s$ is very large. (This
is because $\Sig$ exits $B(1/2)$ where $|\rho|^2 + |z|^2 = 1/2$, and
$2^{-1}\dt \rho^2
- z^2 = \e$, which is far from where either $z$ or $\rho$ is zero.)
Moreover,
according to (4.2), $p_0$ is the intersection number between $C_s$ and
$\Sig(0, \e)$
when $s$ is large. Now, rescale the discussion back to $C$ to
 conclude that there
exists $\d > 0$ such that when $t'  \in Z$ obeys $|t_0 - t'| <
\d/100$, and when $s$
is sufficiently small, then the intersection between $C \cap
 B(\d)$ and $\Sig(t' ,
s)$ is contained in the ball $B(\d/2)$ and the intersection number
equals $p_0$.

\section{Regular points}
\setcounter{equation}{0}

In this section, suppose that $C$ is a pseudo-holomorphic subvariety
with the
property that all points in $Z$ have $n_{\pm}  = 0$. With the preceding
understood, call a point $t \in Z$ regular when $t$ is a local minimum of
 the
function $p(\dt )$. This is to say that $p(t) \leq  p(t' )$ for
all $t'$  sufficiently close to $t$. Let $Z_R$ denote the set of
regular points in $Z$. The complement of $Z_R$ will be denoted by $Z_S$ and
a point
in $Z_S$ will be called {singular}.

The first observation here is that the set of regular points is open and
dense in $Z$. The fact that $Z_R$ is open follows from
Proposition 6.1, and the
fact that $Z_R$ is dense follows because $p(\dt )$, being
integer valued and
bounded from below, takes on a local minimum on every open set. Note that
$p(\dt )$ is locally constant on $Z_R$ when considered as
integer valued
functions on $Z_R$.

The following proposition summarizes the behavior of $C$ near a regular point.

\medskip
{\bf Proposition 7.1}\qua{\sl {The functions $p$ and $q_{\pm}$  (from} (4.1))
{are constant
on $Z_R$, and the latter take values in} $\{0, 1\}$. {Furthermore,
 given $t_0
\in Z_R$, there exists $s_0 > 0$ with the following significance:
Let $C_*$
denote a subset of $C$ where} $t \in (t_0 - s_0, t_0 + s_0)$ {and}
 $(\rho^2 +
z^2)^{1/2} < s_0$. {Then $C_*$ decomposes as the disjoint union}
$C_* = C_{\rho}
\cup C_{z+}\cup C_{z-}$ {where:}
\begin{enumerate}
\item If $q_{\pm} = 0$, {then} $C_{z\pm}$  {is empty;
otherwise} $C_{z\pm}$
{is the restriction of} $\{(t, x, y, z): x = y = 0$ {and} $\pm z >
0\}$
{from Example} (1.4).
\item $C_{\rho}$ {consists of a finite set of components, and $|z| <
s_0/100$ on
each.}
\item {A component of $C'$  of $C_{\rho}$ has the following properties:}

\leftskip 20pt

\item[\rm(a)] {Where} $\rho = (x^2 + y^2)^{1/2} < s_0$,
{the pair of functions $(t, \rho)$ give a proper
embedding of $C'$  into}  $(t_0 - s_0, t_0 + s_0) \x  (0, s_0)$.

\item[\rm(b)] {The restrictions of the functions $z$ and $\phi = \arctan(y, x)$ to
$C'$  extend from where
$t \in (t_0 - s_0, t_0 + s_0)$ and $\rho\in (0, s_0)$ to where $\rho = 0$
as real
analytic functions.}

\item[\rm(c)] {The Taylor's expansions of $\phi$ and $z$ off of $\rho = 0$ are}
\begin{eqnarray}
\phi &=& \phi_0(t)-8^{-1}\phi''_0(t)\dt\rho^2 +{\cal O}(\rho^4),\nonumber\\
z &=& 4^{-1}\phi''_0(t)\dt\rho^2 +{\cal O}(\rho^4) .
\end{eqnarray}

\leftskip 0pt

\end{enumerate}}

\medskip
The remainder of this section is occupied with the statement and proof of
Proposition 7.2, below; the latter implies Assertions 1 and 2 of Proposition
7.1. The proof of Assertion 3 of Proposition 7.1 is completed in Section~8
and in the Appendix.

\medskip
{\bf Proposition 7.2}\qua {\sl{Let $t_0 \in Z_R$. Let $p = p(t_0)$.
 There is an open
neighborhood $U\subset Z_R$ of $t_0$ with the following significance: Given
 $\e$,
there exists $s > 0$, such that the set of points in $C$ where $t \in U$
and $\rho^2 +
z^2 < s^2$ is a disjoint union $C_{\rho} \cup C_{z+} \cup C_{z-}$ which
 have the
following properties:}
\begin{itemize}
\item $C_{\rho}$ {has at most $p$ components and $|z|/r <
\e$ on each.}
\item $C_{z\pm}$  {is either empty or else the
restriction of}
 $\{(t,x, y, z): x = y = 0$ {and} $\pm z > 0\}$ {to where $t \in
U$ and}
$\rho^2 + z^2\leq  s^2$.
\end{itemize}}

\medskip
Aside from giving Assertions 1 and 2 of Proposition 7.1, this last
proposition also implies the initial assertion that $p$ and $q_{\pm}$  are
constant on $Z_R$. Indeed, the constancy of $p$ follows from
 that of $q_+$ and $q_-$
since $p - (q_+ + q_-)$ is the intersection number of
$C$ with a linking 2--sphere
around $Z$. The constancy of $q_{\pm}$  follows from the
following claim: If $C$
contains a non-empty open set which is identical to an open set in some
connected, pseudo-holomorphic subvariety $C'$ , then $C$ contains
the whole of $C'$. Because $C$ is assumed closed in $S^1  \x  (B^3 -
0)$, the latter claim follows from the standard regularity theorems
about
pseudo-holomorphic maps (see, eg~\cite{PW} or \cite{Ye}.  Alternately, one can
invoke Aronszajn's unique continuation theorem in \cite{Ar}.)

The proof of Proposition 7.2 is broken into two parts corresponding to the
two points of the proposition.

\subsection*{a)\qua Proof of the first assertion of Proposition 7.2}

Although somewhat lengthy, the first assertion is essentially a corollary of
the monotonicity result from Proposition 2.1.  In any event, the proof is
conveniently divided into four steps.

{\bf Step 1}\qua The first observation concerns $\s$ from
(2.2), but
considered as a function of two variables, $t$ and $s$.

\medskip
{\bf Lemma 7.3}\qua {\sl {Let $t_0 \in Z_R$. Let $p = p(t_0)$ and
$q = q(t_0)$.
Then there is an open neighborhood  $U\subset Z_R$ of $t_0$ with the
following significance: For each $t \in U$, introduce the function $\s$ of}
(2.2) {and thus consider the latter as a function, $\s(s, t)$ of variables
$s > 0$ and $t \in U$.}  {Given
$\e > 0$,
there exists $s_\e > 0$ such that for all $s < s_\e$ and} $t \in U$,
\begin{equation}
s^{-3}\s(s,t)-(p+2q)\dt {\textstyle{2\over3}} < \e .
\end{equation}}

\med
{\bf Proof of Lemma 7.3}\qua Choose the set $U$ to be connected and to have
compact closure in $Z_R$. This implies that the function $p(\dt )$ is constant
on the closure of $U$. Now, supppose that there exists $\e > 0$, but
that no $s_\e
> 0$ exists which makes (7.2) true for all $t \in U$. Then, there would
exist an infinite sequence $\{t_m\}_{m=1,\dots}\subset U$ such that for
each $m$, the
inequality in (7.2) is violated using $s = 1/m$. Let $t \in Z$
denote a limit
point of the sequence $\{t_m\}$. Since $t \in Z_R$ and $p(t) = p$, there
exists $s' > 0$ such that for all $s < s'$, the
inequality in (7.2) holds with $\e$ replaced by $\e/2$. Fix such $s$. It is an
exercise to verify that $\s(s, \dt  )$ is a continuous function on
$Z$; and this
implies that the inequality in (7.2) holds for those $t_m$ which are
sufficiently close to $t$.  This contradicts the initial assumptions and thus
proves the lemma.

\medskip
{\bf Step 2}\qua A straightforward modification of the preceding
argument proves:

\medskip
{\bf Lemma 7.4}\qua{\sl{ Let $t_0 \in Z_R$ and set $p = p(t_0)$ and $q =
q(t_0)$. There
is an open neighborhood $U\subset Z_R$ of $t_0$ with the following
significance: For
each $t\in U$ and $s > 0$, introduce the function}
\begin{equation}
\g(s,t)=\int_c\chi_s\dt d\phi\wedge dh .
\end{equation}
{Then, given $\e > 0$, there exists $s_\e > 0$ such that for all $t \in U$
and}
$s < s_\e$,
\begin{equation}
s^{-3}\dt\g (t,s) < \e .
\end{equation}}

\medskip
(Note that (1.10) and (1.11) imply that $d\phi\wedge  dh$ restricts to $C$ as an
everywhere non-negative 2--form.)

\medskip
{\bf Step 3}\qua This step is concerned with the
following claim:  There
exists a neighborhood $U$ of $t_0$ and, given $\e > 0$, there
 exists $s_\e > 0$ such
that either
\begin{equation}
|z|/\rho < \e  \quad {\mbox{or}}\quad
\rho/|z| < \e
\end{equation}
at points of $C$ where $t \in U$ and $\rho^2 + z^2 \leq  s_\e^2$.

To prove this claim, consider an open interval $U\subset Z_R$ which
contains $t_0$, is
connected and which has compact closure in $Z_R$. Now, suppose that there
exists $\e > 0$ and no $s_\e > 0$ which makes (7.5) true at points of
$C$ where $t
\in U$ and $\rho^2 + z^2 \leq  s_\e^2$. There would then be an infinite
sequence
$\{t_m\}\in U$ with the following properties: There is a point $w_m$ on
$C$ with
$t$ coordinate $t_m$, with $r = (\rho^2 + z^2)^{1/2}$ coordinate $r_m
= 1/m$ and with
$r_m^{-3}\dt h(w_m) \geq  \e^3$. No generality is lost by assuming that
$\{t_m\}$
converges to a point $t_{\i}$  in the closure of $U$.

With the preceding understood, for each $m$, rescale so that the ball of
radius $4r_m$ centered at $t_m \in Z$  becomes the ball of radius 4. Then,
translate along $Z$ so that $t_m$ is moved to the origin. The result is a
sequence $\{C_m\}$ of pseudo-holomophic submanifolds on the set where $t^2 +
\rho^2 + z^2 < m/\z$ and where $\rho^2 + z^2 > 1/16$. Here, $\z \geq  1$
is a fixed
constant. Note that $w_m$ is sent to a point with $t$ coordinate zero and
with $r
= 1$.

The sequence $\{C_m\}$ has a uniform bound for $\int\o$ on any fixed ball
about
the origin. Thus, Proposition 3.3 finds a subsequence which converges
geometrically (as described in Definition 3.1) where $\rho^2 + z^2 > 1/16$ to a
pseudo-holomorphic subvariety, $C_{\i}$. Furthermore it follows from
Lemma 7.4 that $d\phi = dh = 0$ on $C_{\i}$. This implies that $C_{\i}$
is a union of subvarieties, where each component has the form $\phi =$
constant and $h =$ constant from Example 1.6. Note in particular,
 that there
must be at least one component of $C_{\i}$  where $h \geq  \e^3$.

To obtain the final contradiction, consider now fixing some $t'
\in Z$ very near to $t_{\i}$  and counting the intersections of $C$
 with
a pseudo-holomorphic disk which obeys the conditions $t = t'$, $f=$
constant $> 0$ and $|z|/\rho < 1/4$. (This is a subdisk of one of the
cases in Example 1.6.) Consider, in particular, the case where $f = c$ where $c
= |t_{\i} - t'|$. It then follows from (4.1) (using
$t_{\i}$  in (4.1) instead of $t_0$) that this intersection number is
equal to $p$ when $t'$  is sufficiently close to $t_{\i}$.
(Remember that all intersection points count positively to this intersection
number.) With the preceding understood, fix $t'$  very close
$t_{\i}$  and consider decreasing $c$. For $c$ sufficiently small, it follows
from (4.1) (using $t'$  for $t_0$ in (4.1) this time) that the
intersection number is again, equal to $p$. On the other hand, if
$t'  = t_m$ and $m$ is large, it follows from the behavior of $C_m$ at
large $m$ that the intersection number must be less than $p$ when $c$
 decreases
much past $r_m^2$.

It follows from the preceding step that there is an open neighborhood $U$
 of
$t_0$ and, given $\e$, there exists $s_\e > 0$ such that where $t \in
U$ and $\rho^2 +
z^2 \leq s_\e^2$, $C$ can be written as $C = C_\rho \cup C_z$,
where $|z|/\rho< \e$ on
$C_r$ and $\rho/|z|< \e$ on $C_z$. That is, (7.5) holds.

\medskip
{\bf Step 4}\qua With (7.5) understood, the proof of the
first
assertion of Proposition 7.2 is complete with the establishment of a
bound by $p$ on
the number of components of $C_r$. To obtain such a bound, agree first
 to
re-introduce the family of pseudo-holomorphic cylinders given by $t =
t'$, $f =$ constant $> 0$ and $|z|/\rho < 1/2$ where
$t'$  now ranges over some small interval $U\subset Z_R$ which contains
$t_0$. If $c > 0$, let $\Sig(t' , c)$ denote the $t = t'$
and $f = c$ case.

Now, let $C'$  be a component of the intersection of $C_\rho$ with some
small neighborhood of $t_0$ in $S^1 \x  B^3$. Note that when $t'$  is
near to $t_0$ and $c$ is small, $\Sig(t' , c)$ has positive intersection
number with $C'$. (This follows from Proposition 3.2.) Because
of (7.5), this intersection number will be independent of both
$t'$  and $c$ as long as $|t'  - t_0|$ and $c$ are both small.

On the other hand, the intersection number of $\Sig(t_0, c)$ with $C_\rho$
is, for
small
$c$, equal to $p$. (This follows by rescaling (4.2).)  For such $c$, the
intersection number of $\Sig(t' , c)$ with $C_\rho$ must equal $p$ for all
$t'$  near to $t_0$. Furthermore, because of (7.5), this
intersection number will also be independent of $t'$  and $c$ as
long as both $|t'  - t_0|$ and $c$ are small. Thus, a counting
argument shows that there is a small neighborhood of $t_0$ in $S^1  \x
B^3$ whose
intersection with $C_\rho$ has at most $p$ components.

\subsection*{b)\qua Proof of the second assertion of Proposition 7.2}

The task here is to verify that $C_{z\pm}$  have the indicated form.
To begin,
consider a component $C'$  of $C_z$. It follows from (7.5) that one
can assume, without loss of generality that $z > 0$ on this component.
Concerning $(t, u)$, consider:

\medskip
{\bf Lemma 7.5}\qua {\sl {The functions $t$ and $u = (-f)^{1/2}$ restrict to a
neighborhood of all but a set of at most countably many points in $C'$  as
local coordinates. The countable set, $\L$, here consists of the singular
points of $C'$  and the critical points of $dt$, (or, equivalently,
$du$). Note that this set $\L$ has no accumulation points in} $C'$.}

\medskip
{\bf Proof of Lemma 7.5}\qua By assumption, $C'$  has at most a
countable number of singular points. With the preceding understood, let $o
\in C'$  be a smooth point. Then, by assumption, there is a
pseudo-holomorphic embedding $f$ of the standard disk in $\Bbb C$ into
$S^1  \x  (B^3 - \{0\}$) which sends 0 to $o$ and the disk to a
neighborhood of
$o$ in $C'$. Furthermore, it is an exercise to show that there are
coordinates $(\s^1 , \s^2)$ on the disk with respect to which the
pull-backs of the
functions $t$ and $u$ obey the equations
\begin{equation}
t_1\kappa =u_2 \quad {\mbox{and}} \quad
t_2\kappa =-u_1.
\end{equation}
Here, $\k = (\rho^2 + 4 z^2)^{1/2}/2u = 1 + {\cal O}(r^2/z^2)$.
(In (7.6), the
subscripts ``${}_1$'' and ``${}_2$'' denote the partial derivatives with
respect to
$\s_1$  and $\s_2$.) Equation (7.6) implies, first of all, that the
 critical points of
$t$ and $u$ agree on $C'$.

Now, suppose that $o$ were also an accumulation point of critical
 points of
$dt$. The argument below deduces a contradiction from this
supposition. To
begin the argument, note first that (7.6) implies a homogeneous Laplace
equation for either $t$ or $u$ on the disk. The equation for $t$ reads:
$(\k t_1)_1 +
(\k t_2)_2 = 0$. This equation implies, via Aronszajn's unique
continuation
theorem \cite{Ar}, that $t - t(0)$ does not vanish to infinite order at 0.

On the other hand, if 0 were an accumulation point of the critical set of $t$,
then Taylor's theorem with remainder and the Laplace equation would force $t
- t(0)$ to vanish to infinite order at 0. The argument for this last
assertion is as follows: If $t$ has a non-zero Taylor's coefficient, then
 it
has one of smallest order, say $n$. With $n$ understood, use Taylor's
theorem to
write $t = t(0) + \sum_{0\leq k\leq n} a_k \s^k_1 \s_2^{n-k} +
{\cal O}(|\s|^{n+1})$.
The definition of n requires some $a_k \neq   0$. If $dt$ is to
vanish arbitrarily
close to 0, then there must exist $(\s_1,\s_2) \neq   0$, with the property
that
$\sum_{0\leq k\leq n} a_k k \s_1^{k-1} \s_2^{n-k} = 0$ and
 $\sum_{0\leq k\leq n}
a_k (n - k) \s_1^k \s_2^{n-k-1} = 0$. By rotating the coordinate
system if necessary,
one can assume, with out loss of generality, that this occurs where
$\s_2 = 0$. Thus,
the vanishing of $dt$ at points arbitrarily close to 0 implies that
$a_{n-1} = a_n =
0$. Meanwhile, the aforementioned Laplace equation for $t$ requires that
the sum\newline
$\sum_{0\leq k\leq n} a_k (k(k-1) \s_1^{k-2} \s_2^{n-k} + (n-k)(n-k-1) \s_1^k
\s_2^{n-k-2})$ vanish at all $\s$. This last fact implies, via an induction
on $n -
k$, that $a_k = 0$ for all $k$, thus contradicting the defining
property of $n$.

With the preceding lemma understood, remark that a neighborhood of any point
in $C' - \L$ where $\rho^2 + z^2$ is small is given as the set of points
$(t, z, x, y)$ where
\begin{itemize}
\item\quad $x=x(t,u), \ \ y=y(t,u)$
\item\quad $z^2=u^2+2^{-1}\dt (x^2+y^2)$
\item\quad $u\in (0,s_1) \ \ {\mbox{and}} \ \ t\in U$.\itemnum{7.7}
\end{itemize}
Here, $s_1$  is positive and much less than $s_0$. In particular, given
$\e > 0$,
one can choose $s_1$  so that
\begin{equation}
(|x|+|y|)/u < \e .
\end{equation}
Here, $x$ and $y$ are implicitly restricted to a neighborhood of the
given point
in $C'$  where they are considered as functions of the variables
$t$ and $u$.

Now, let $\eta = x + i\dt y$ and introduce the complex coordinate $w = t +
i\dt u$ and the associated $\p   = 2^{-1}\dt (\p_t - i\cdot \p_u)$. Then,
$\eta$ obeys
an equation of the form
\begin{equation}
\p\eta +i\dt\eta/4z=a\dt\nabla\eta +b\dt\nabla\bar\eta  ,
\end{equation}
where $a$ and $b$ are functions of $u$ and $\eta$ which obey $|a|+|b|\leq
\z\dt |\eta|^2/u^2$. Note that (7.9) is augmented with the knowledge
that  $|\eta|/u$ is uniformly bounded by $\e$. Furthermore,
\begin{equation}
\lim_{u\to 0} |\eta|/u=0 .
\end{equation}
Now, it is important to realize that the set $\L\subset C'$  need not be
empty, and may, as far as has been established as yet, have infinitely many
members. In particular, for the sake of argument, one can assume that $\L$
 has
an infinite number of elements, for othewise, one could decrease $s_0$ and
 thus
avoid them entirely in the subsequent discussion.

In any event, the complex function $\eta$ obviously extends (by definition)
as a
continuous function on $C'$, but not so as a global function of
the complex coordinate $w$.  As a function of the latter, $\eta$ is
multivalued,
where the different values correspond to different sheets of $C'$
over the $w$ plane. The points of $\L$ then correspond to places
where the sheets
coincide. (In any event, there are at most $q_+$ sheets because each sheet has
intersection number 1 with a $t =$ constant, $f =$ constant $< 0$
psuedo-holomorphic submanifold.)

The structure of $C'$  as viewed using the complex coordinate $w$
can be described slightly differently as follows: The assignment of the
complex number $t + i u$ to a point in $C'$ defines a continuous
map, $\pi$ , from $C'$ to the $(t, u)$ plane. In a neighborhood
of $w \not\in\pi (\L)$, the sheets of $C'$ are parametrized by a
set $\{\eta_1,\dots, \eta_{q'}\}$ of complex valued functions with
distinct values at each point. Here, $q' \leq  q_+$ is fixed.
Each $\eta=\eta_j$ obeys (7.9) and (7.10). Meanwhile, the points in
$\pi(\L)$ are
characterized by the fact that the values of some pair from $\{\eta_1, \dots,
\eta_{q'}\}$ coincide. Furthermore, as $w$ traverses a circle around
some point in $\pi(\L)$, the set $\{\eta_1, \dots, \eta_{q'}\}$ can
be permuted by some non-trivial monodromy.

With the preceding understood, consider now the complex valued function
\begin{equation}
\th =\prod_{1\leq j\leq q'} \eta_j .
\end{equation}
Note that $\th$ is a globally defined function on a neighborhood of 0
in the
upper half of the $w$ plane which obeys
\begin{equation}
\lim_{u\to 0} \th/u^{q'}=0 .
\end{equation}
Furthermore, by virtue of (7.9), this complex valued function obeys an
equation of the form
\begin{equation}
\p\th +i\dt\th\dt q'\dt (1+m)/4u=0  ,
\end{equation}
where $m$ is a smooth complex valued function away from $\pi(\L)$
about which
more will be said below.

Here is the derivation of (7.13): First, observe that both the Euclidean
metric and $\omega$ are invariant under the group of rotations in the (x, y)
plane. This implies that if $\eta = \eta_j$ is a solution to (7.9), then so is
$\l\dt\eta_j$ where $\l$ is any unit length complex number. This fact implies
that the functions $a$ and $b$ in (7.9) must have the form
\begin{itemize}
\item\quad $a=\a\dt |\eta_j|^2/u^2 $
\item\quad $b=\b\dt\eta^2_j/u^2$.\itemnum{7.14}
\end{itemize} 
Here, $\a$ and $\b$ are functions of $u$ and $\eta_j$ which obey $|\a|+|\b|
\leq \z$. With (7.14) understood, then (7.9) for $\eta=\eta_j$ can be
 written as $\p
\eta_j + i\dt \eta_j\dt (1 + m_j)/4u = 0$, where
\begin{equation}
m_j=(u/z-1)-4\dt(\a\dt\bar\eta_j\dt\nabla\eta_j+\b\dt\eta_j\dt\nabla
\bar\eta_j)/u .
\end{equation}
With $m_j$ understood, then (7.13) for $\th=\prod_j\eta_j$ follows
immediately with $m= \sum_j m_j$.

The subsequent manipulations of (7.13) require an estimate for the size of
$m$. Here, remark that $m$ is smooth away from the points in
$\pi(\L)$. The
following lemma describes $|m|$ in more detail.

\medskip
{\bf Lemma 7.6}\qua{\sl {Although $|m|$ may be unbounded near the point
of $\L$, this
function is none the less locally square integrable.  Moreover, there exists
$c \geq  (0, 1)$ and, given $\e > 0$, there exists $s_\e > 0$ such
 that the
following is true: Suppose that $s, r \in (0, s_\e)$. Let $A$ denote the
intersection of the $u \geq  0$ portion of the $(t, u)$ plane with an annulus
with center where $|t| < c$ and $u = s$, and whose inner and
outer radii
are $r/2$ and $r$, respectively. Then:}

\begin{itemize}
\item {If $r \leq  s/4$, then} $\int_A |m|^2 dt du
 \leq \e^2 (r/s)^c
s^2$.

\item {If $r \geq  s/4$, then} $\int_A |m|^2 dt du
 \leq \e^2 r^2$ .
\end{itemize}}

\medskip
This lemma is proved below. Accept it for now.

With (7.13) understood, then the second assertion of Proposition 7.2 follows
with a proof that the only solution to (7.12) with $\lim_{u\to 0} \th = 0$
is the function $\th\equiv  0$. (Remember that $C'$ is assumed to
be connected.)

Progress on this task requires

\medskip
{\bf Lemma 7.7}\qua  {\sl {Given $c > 0$, there exists a constant $\z > 0$ with the
following significance: Let $m$ be a square integrable, complex valued function
on the
$(t, u)$ upper half plane which obeys the conclusions of Lemma 7.6
for some $\e
> 0$ and for the given constant $c$. Suppose in addition that there exists
$R\in (0, c)$ such that $m = 0$ where $(t^2 + u^2)^{1/2} \geq R$.
Then there exists
a real valued function, $\k$, on the $(t, u)$ upper half plane which, where
$t^2 +
u^2 \leq  R/2$, obeys}
\begin{itemize}
\item \quad $\p  \k + i\dt \k\dt (1 + m)/4u = 0$

\item \quad $\k > 0$

\item \quad $e^{-4\e} u^{q' /2+\z\e} \leq  \k \leq  u^{q' /2-\z\e}$.
\end{itemize}}

\medskip
This lemma is also proved momentarily. Accept it for now to complete the
analysis of (7.13).

In particular, take an open set $U'\subset  U$ with compact
 closure in
$U$, and take some very small $\e > 0$ in Lemma 7.6 to obtain $s_\e$. Then,
modify
$m$ from (7.13) with a cut off function so that the result (still to
be called
$m$) is unchanged when $(t, u) \in U'  \x (0, s_\e/2)$ but vanishes
when $(t, u)$ is sufficiently close to the frontier of $U \x  (0, s_\e)$.
With this
$m$ understood, use Lemma 7.7 to produce $\k$.  Then, note that $\th\,'
= \k\dt\th$ obeys the equation $\p\th\,'  = 0$ on $U'\x
(0, s_\e/2)$ and extends as 0 on $U'\x\{0\}$. The latter
implies that $\th\,' \equiv 0$ as this is the only anti-holomorphic
function which vanishes on a line segment. Thus, $\th\equiv  0$ as
 claimed.

The proof of the Proposition 7.2 is thus completed with the proofs of Lemmas
7.6 and 7.7. These will be taken in reverse order.

\medskip
{\bf Proof of Lemma 7.7}\qua Write $\k = u^{q' /2}\dt e^{\tau}$ and
then $\tau$ must obey the equation
\begin{equation}
\p\tau +i\dt m/4u=0.
\end{equation}
Here is one solution to (7.16):
\begin{equation}
\tau(w)=-(4\pi)^{-1} \int_{im(w')>0} \! m(w')
(( \bar w-\bar w')^{-1}\!-(\bar w - w')^{-1})/im(w)d\bar w'\!\wedge dw'
\end{equation}
To see that this integral is well defined, note that the integral above can
be rewritten as
\begin{equation}
\tau(w)=-i(2\pi)^{-1} \int_{y\geq 0} m(w')/
(( \bar w-\bar w')(\bar w-w'))d\bar w'\wedge dw' .
\end{equation}
To estimate the norm of $\tau$, consider breaking the domain of integration
into two parts. The first part has $|w - w' | \leq {\mbox{im}}(w) = 2u$;
and the second has $|w - w' | \geq  2u$. Because $|w - \bar{w}'|\geq u$
when $w$ and $w'$  lie in the upper half
plane, the contribution from the first region has absolute value no larger
than
\begin{equation}
(2\pi)^{-1}u^{-1} \int_{|w-w'|\leq 2u} |m|
|w-w'|^{-1} id\bar w'\wedge dw' .
\end{equation}
The claim here is that (7.19) is no greater than $\z \e$, where $\z$
depends only
on the constant $c$, but not on $\e, \ R$ or on other properties of $m$.
 Here is a
proof of this claim: Break the integration region into concentric annuli
indexed by $n \in \{0, 1,\dots\}$, where the $n$'th annulus is
characterized by the requirement that $2^{-n}u \leq |w - w'|
\leq 2^{-n+1}u$. The contribution to (7.19) from the $n$'th annulus
is no
greater than $\pi^{-1}2^n u^{-2}$  times the $L^1$ norm of $|m|$
over the $n$'th
annulus.  This number is no greater than $\sqrt{2} u^{-1}$ times the
$L^2$ norm of
$|m|$ over this same annulus. And, according to Lemma 7.6, the $L^2$ norm in
question is no greater than $\e (2^{-n})^{c/2} u$. Thus, the n'th
annulus contributes
no more that $\sqrt{2} \e 2^{-nc/2}$ to the integral in (7.19). With this
last point
understood, it follows that (7.19) is no greater than $\sqrt{2} \e
\sum_{n\geq 0}
 2^{-nc/2}$. As this last sum converges, the asserted uniform bound of
(7.19) by $\z \e$ is obtained.

Meanwhile, the absolute value of the contribution to (7.18) from the second
region is no larger than
\begin{equation}
(2\pi)^{-1}\dt \int_{|w-w'|\leq 2u} |m(w')|/|w'-w|^2 i\dt d\bar w'
\wedge dw' .
\end{equation}
The latter integral can be bounded by breaking the integration region
into $u\geq 0$ half annuli, where  these half annuli are indexed
by $n \in\{1, \dots , N\}$, here $N \leq 1 + 4 \ln(R/u)$. In this
case, the $n$'th half annulus is characterized by $2^nu \leq |w - w' |
\leq 2^{n+1}u$. The contribution to (7.20) from the $n$'th half
annulus is no greater than $ (2\pi)^{-1}2^{-2n}u^{-2}$ times the $L^1$
norm of $|m|$ over said half annulus.  Meanwhile, this $L^1$ is no
greater than $(2\pi )^{-1} 2^{-2n}u^{-2}$ times the $L^2$ norm of
$|m|$ over the half annulus in question. Apply the second point in
Lemma 7.6 to bound this last $L^2$ norm by $\e 2^{n+1} u$, and thus
bound the contribution from the n'th half annulus by a constant $4
\e$. With the preceding understood, it follows that (7.20) is no
greater than $4 \e N \leq 16 e |\ln u|$.  Thus,
\begin{equation}
|\tau|\leq\zeta\e\dt(1+|\ln u|) .
\end{equation}
where $\z$ depends only on the constant $c$ in the statement of
the lemma.

This section ends with the proof of Lemma 7.6.

\medskip
{\bf Proof of Lemma 7.6}\qua The proof of Lemma 7.6 is broken into six steps.

{\bf Step 1}\qua It proves useful for a subsequent application in
Section~9, to establish a version of Lemma 7.6 under slightly weaker
assumptions
then those given. Thus, the proof starts with a digression to review the
background for these assumptions.

To start the digression, suppose that $s_0 \in (0, 1/4)$ has been fixed,
together with a constant $\a \geq  0$. Let $W$ denote the domain in the
$(t, u)$
plane where
\begin{itemize}
\item\quad $t^2+u^2 < s_0$
\item\quad $u > \a|t| $.\itemnum{7.22}
\end{itemize} 
Let $W' \subset S^1 \x  (B^3 -\{0\})$ denote the domain where $z > 0$,
$\rho/z < 1/4$ and where (7.22) holds using $u = (-f)^{1/2}$.
(The function $u$ is
well defined where $0 \leq  \rho/z < 1/4$.)

Now, suppose that $C'\subset W'$  is a finite energy,
pseudo-holomorphic subvariety for which the following regularity assumptions
hold: Given $\e > 0$, there exists $s_\e > 0$ such that when $s <
 s_\e$, then
\begin{itemize}
\item\quad $\rho/u < \e$
\item\quad $\int_{C'\cap\{u < s\}} d\phi \wedge dh\leq \e s^3$.\itemnum{7.23}
\end{itemize} 
Note that Lemma 7.4 and (7.7) guarantee (7.23) for the subvariety in Lemma
7.6.

With $C'$ as described above, the association of $w = t + i u$ to
its points defines a continuous map, $\pi$ , from $C'$ to the $u
> 0$ half of the $(t, u)$ plane. As with the example from Lemma 7.6, the
function $w$ (via the map $\pi$) defines a complex coordinate near
all but a
set, $\L$, of at most countably many points of $C'$ . This set $\L$
contains the singular points of $C'$ and the critical points of
$u$ (or $t$).

The push forward by this map $\pi$  of the function $\eta\equiv (x + i y)$
defines an unordered set $\{\eta_j\}$ of some number, $q'$, of
locally defined functions on the domain $W$ in (7.22). (Note that
$q'$  is a constant, as it measures the intersection number of $C'$ with
certain
members of the ($t = c$, \ $f = c'$) case of Example 1.6.) These
$\{\eta_j\}$ can be
ordered as smooth functions near any point which is not in
$\pi(\L)$. Note that each
$\eta_j$ satisfies (7.9) at points not in $\pi(\L)$. The coorresponding
 $m_j$ is still
given by (7.15).

With the preceding understood, set $m = \sum_j m_j$. Then, the assertions of
Lemma 7.6 make sense in the present context if $A$ is suitably interpreted. In
particular, consider the following lemma:

\medskip
{\bf Lemma 7.8}\qua{\sl{Let $W, W'$  and $C'$ be as
described above. Let $m = \sum_j m_j$. Then $|m|$ is locally square
integrable. Moreover, there exists $c \subset  (0, 1)$ and, given $\e > 0$,
there exists $s_\e > 0$ such that the following is true: Suppose that
$s, r
\in (0, s_\e)$. Let $A$ denote the intersection of $W$ with an annulus
whose
center has $|t| < \min(c, s/(2\a))$ and $u = s$, and whose inner and outer
radii are $r/2$ and $r$, respectively. Then:}

\begin{itemize}
\item {If} $r \leq  s/4$,
{then} $\int_A |m|^2dt\, du \leq \e^2 (r/s)^c s^2$.

\item {If $r \geq  s/4$, then} $\int_A |m|^2 dt \, du
 \leq \e^2 r^2$.
\end{itemize}}

\medskip
This lemma is reduced in the next step to a related lemma.

\medskip
{\bf Step 2}\qua It follows from (7.15) that
\begin{equation}
|m|^2|_A\leq\zeta\dt[\sup_{w\in A} (\sum_j |\eta_j|/u)]^2
\sum_j |d\eta_j|^2 .
\end{equation}
Then, (7.22) implies that given $\d > 0$, there exists $s'_{\delta}$ such
that when $r, s \in (0, s'_{\d})$ then the left-hand side of
(7.24) is no greater than $\z \d^2 \sum_j |d\eta_j|^2$. Thus, it is enough
 to
bound the integral of $\sum_j |d\eta_j|^2$ over the region $A$.

With the preceding understood, introduce the function $\l \equiv
 \sqrt{z}\eta$  and
the corresponding set $\{\l_j = \sqrt{z}\eta_j\}$. Since $\rho/u <<
1$ on $C'$ in the
region of interest, one has
\begin{equation}
|d\eta_j-z^{-\frac 12} d\l_j|^2 \leq\zeta |d\eta_j|/u .
\end{equation}
Hence, it is more than sufficient for the proof of Lemma 7.8 to establish
the following lemma:

\medskip
{\bf Lemma 7.9}\qua{\sl{Let $W, W'$  and $C\subset W'$  be as
described above. The function $\sum_j |d\l_j|^2$ on $W$ is locally
integrable. Moreover, there exists $c > 0$, and, given $\e$, there exists
$s_\e >0$ which are such that when $r, s \in (0, s_\e)$, then:}
\begin{itemize}

\item {Suppose that $r \leq  s/4$. Let $D\subset  W$ denote the
 disk with
center where $|t| \leq  \min(c, u/(2\a))$ and $u = s$, and with
 radius $r$.
Then,} $\int_D \sum_j |d\l_j|^2 \leq  \e^2 (r/s)^c s^3$.

\item {In general,} $\int_{\a|t|\leq u\leq r\atop|t|\leq
c}
\sum_j |d\l_j|^2 \leq  \e^2 r^3$.
\end{itemize}}

\medskip
{\bf Proof of Lemma 7.8 given Lemma 7.9}\qua The local integrability of
$\sum_j |d\eta_j|^2$  follows from
that of $\sum_j |d\l_j|^2$, and the first point follows directly
from the
corresponding point of Lemma 7.9. The second point of Lemma 7.8 follows
from that of
Lemma 7.9 by decomposing the region $A$ of Lemma 7.8 into
subregions indexed by $n
\in\{1, 2, \dots\}$, where the $n$'th subregion is defined by the
constraint $2^{-nr}
\leq  u \leq  2^{-n+1}r$. Then, the integral of $|z|^{-1} \sum_j
|d\l_j|^2$ over
the $n$'th subregion is, according to the second point in Lemma 7.9, no
greater than
the following: $\e^2(2^{-n}r)^{-1}(2^{-n+1}r)^3=8\e^2 2^{-2n}r^2$.
Summing the latter
contributions over $n$ bounds the desired integral by $8/3 \e^2 r^2$.
Thus, the
second point of Lemma 7.8 follows if $s_\e$ in the statement of said
lemma is
taken to equal $s_{\e'}$  from Lemma 7.9 where $\e'  = \e
(3/8)^{1/2}$.

\medskip
{\bf Step 3}\qua This step and the subsequent three steps
contain the proof of Lemma 7.9.

\medskip
{\bf Proof of Lemma 7.9}\qua Consider first the assertion that $\sum_j |d\l_j|^2$
is locally integrable. For this purpose, note that the behavior of $\l$
near points in $\L$ can be discerned from the known local structure of a
pseudo-holomorphic map as described, for example, in \cite{PW} or \cite{Ye}. (See also
\cite{MS}.) In particular, the claim here is that if $w_0 \in\L$, then there
exists a constant $\z \geq  1$ such that
\begin{equation}
\sum_j |d\l_j|^2 \leq \zeta|w-w_0|^{-2+2/\z}
\end{equation}
at points $w$ near, but not equal to $w_0$.

To prove (7.26), suppose that $o \in\L$. Then, the  analysis in \cite{Ye} or \cite{PW}
shows that $C'- o$ near $o$ is a disjoint union of smooth
submanifolds. Moreover, the closure of any such submanifold is the image of
the standard disk in ${\Bbb C}$ via an almost everywhere 1--1,
pseudo-holomorphic map $f$ into $S^1  \x (B^3 - \{0\}$). This map
$f$ can
be assumed to send the origin in ${\Bbb C}$ to $o$.

In addition, the coordinate $v$ on ${\Bbb C}$ can be chosen so that the
pull-back via $f$ of the function w obeys an equation of the form
\begin{equation}
\bar\p w+\g_1\p w +\g_2\bar\p\bar w=0  ,
\end{equation}
where $\g_{1,2}$ are smooth functions which obey $\g_1(0) = 0$ and
$|\g_2|\leq  \z
|\eta|^2/u^2$. Here, $\p$   is the complex s derivative on ${\Bbb C}$.
(This last
equation follows readily from (7.6).)  Now, $w$ is a smooth function of $v$
 and
so has a Taylor's series.  By virtue of (7.27), the latter must be such that
\begin{equation}
w+\g_2(0)\bar w=a_0+a_1 v^n+{\cal O}(|v|^{n+1})  ,
\end{equation}
where $n$ is a positive integer and where $a_0$ and $a_1$ are
constants, with
neither non-zero when $|\g_2| < 1$.

Now, to prove (7.26), first use (7.28) to express $v$ as a function of $w$.
Then, as $\l$ pulls back via $f$ as a smooth function, the desired
result
follows using the chain rule.

\medskip
{\bf Step 4}\qua To consider the two points of the
lemma, remark that by
virtue of (7.9) and (7.15), each $\l \in\{\l_j\}$ obeys an equation of the
form
\begin{equation}
\p\l+{\cal P}=0  ,
\end{equation}
at any point not in $\pi(\L)$. Here, $|{\cal P}|\leq \z (|\l|/u+|d\l|)|\eta|/u$
and now $\p$   is the complex $w$ derivative. Then, using
(7.23), one can assume the following: Given $\d$, there exists $s'_{\d}
> 0$ such that when $u < s'_{\d}$, then
\begin{equation}
|{\cal P}|\leq\d(u^{\frac 12}+|d\l|) .
\end{equation}
This last equation implies that where $u \leq s'_{\d}$, then
\begin{equation}
(\zeta^{-1} |d\l|^2-\zeta\d^2 u)dt \wedge du\leq
-id\l\wedge d\bar\l  ,
\end{equation}
where $\z \geq 1$ is a fixed constant.

Equation (7.31) is the basic inequality for the subsequent arguments.

\medskip
{\bf Step 5}\qua This step considers the second integral
inequality in
Lemma 7.9. For this purpose, note that $\{\l_j\}$ can also be written as the
pushfoward via $\pi$  from $C\,'$ of $|h|^{1/2} e^{i\phi}$. Thus, $i
 d\l\wedge
d \bar{\lambda } $  is equal to $dh\wedge  d\phi$. With this understood, (7.29)
implies that
\begin{equation}
\int_{C'\cap\{u\leq r\}}
\sum_j |d\l_j|^2\leq\zeta (\d^2 r^3+
\int_{C'\cap\{u\leq r\}} d\phi\wedge dh).
\end{equation}
It then follows from (7.23) and (7.32) that given $\e > 0$, there exists
 $s_\e > 0$
for which the final point of Lemma 7.9 holds when $r < s_\e$.

\medskip
{\bf Step 6}\qua This step considers the first integral
inequality in
Lemma 7.9. The domain $D$ here is the full disk of radius $r$.  With this
understood, introduce $f(r) = \int_D\sum_j |d\l_j|^2$. Let $C' (r) \subset
C'$ denote
the $\pi$--inverse image of $D$. Then, (7.29) implies that
\begin{equation}
f(r)\leq \zeta\d^2 r^2s-i\int_{C'(r)} d\l\wedge d\bar \l .
\end{equation}
When $r$ is such that $\p C' (r)$ is disjoint from $\L$, then the
latter is the union of a set, $\{\G\}$, of smoothly embedded circles. With
this last observation, use Stoke's theorem to rewrite the integral on the
right-hand side above as
\begin{equation}
-i\int_{C'(r)} d\l\wedge d\bar \l =-i
\sum_{\G\subset\p C'(r)} \int_{\G} (\l-\l_{\G})d\bar\l .
\end{equation}
Here, the sum in question is over the set $\{\G\}$ of components of $\p C'
(r)$, and
$\l_\G$ can be any constant, although a convenient choice is made
below. Note that
there are no anomalous boundary terms in (7.34) arising from the
 points in $\L$. (The
absence of such terms can be verified by pulling the integral on the
 right-hand side
of (7.34) back to the smooth model for $C'$ via its pseudo-holomorphic map
into $S^1
\x (B^3 - \{0\}$). Up on the smooth model, the function $\l$ is smooth while
the inverse image of the points in $L \cap  D$ is also a finite set.)

With regard to (7.34), note that the number of components of $\p C'(r)$ is no
greater than the number, $q'$, of elements in the set $\{\eta_j\}$.
 In addition, the
projection $\pi$  maps each component $\G$ of $\p C'(r)$ as a covering
map onto $\p
D$ with some covering degree $q_\G \leq q'$.

With the preceding understood, here is how to define $\l_\G$: Let $t \in
[0, 2\pi)$
denote a standard angle parameter on $\p  D(r)$. Then, set
\begin{equation}
\l_{\G}\equiv (2\pi r' q_{\G})^{-1} \int_{\G}\l\pi^* d\tau .
\end{equation}
It then follows from the fundamental theorem of calculus that the right-hand
side of (7.34) has absolute value no greater than
\begin{equation}
\zeta\dt\sum_{\G\subset\p C'(r)} r\int_{\G} |\l_\tau|^2\pi^* d\tau
\leq \zeta r\sum_j\int_{\p C'(r)} |d\l_j|^2 d\tau .
\end{equation}
Here, $z \geq 1$ is independent of $r$ and of the center
 which defines $D$.
(The subscript ``$\tau$'' in (7.36) indicates the $\tau$--partial
derivative.)

Together, (7.33) and (7.36) imply the differential inequality
\begin{equation}
f\leq\zeta (\d^2 r^2 s+rf') .
\end{equation}
This last equation can be integrated from any value of $r < s/4$ to $r =
 s/4$
and the result is the inequality
\begin{equation}
f(r)\leq 4^{1/\zeta} (r/s)^{1/\zeta} (f(s/4)+\zeta\d^2 s^3).
\end{equation}
To end the story, note that it has already been established that there
exists, given $\d > 0$, a number $s'_{\d} > 0$ such that $f(s/4)
\leq \d^2 s^3$ when $s < s'_{d}$. Thus, the right side of (7.38)
is no greater than $\z (r/s)^{1/\z} d^2 s^3$ when $s < s'_{\d}$. With this
last fact given, take $c$ in Lemma 7.9 to equal $1/\z$ and take
$s_\e$ to equal
$s'_{\d}$ for $\d = \e/\sqrt{\z}$.

\section{A proof that $C_{\rho}$ is smooth}
\setcounter{equation}{0}

The purpose of this section is to prove Assertions 3(b) and 3(c) of Proposition
7.1 under the assumption of Assertion 3(a). Assertion 3(a) is also proved here
in the case where $p = 1$, but the proof of Assertion 3(a) for $p > 1$ is
relegated to the Appendix. (The latter proof is quite lengthy.)

The remainder of this section is divided into two subsections. Subsection
8.a establishes some general continuity properties of $C_\rho$ near $\rho
= 0$ when
parametrized as a multi-valued map of the $(t, \rho)$ coordinates.
Subsection 8.b
uses the results from subsection~8.a to prove the $p = 1$ version of
Assertion 3(a),
and also Assertions 3(b) and 3(c) of Proposition 7.1 under the assumption that
Assertion 3(a) holds.

Note that because of the second assertion of Proposition 7.2, no generality
is lost by assuming that $C = C_{\rho}$ near the given point $t_0 \in
Z_R$. This
assumption will be made implicitly below.

\subsection*{8.a\qua Holder continuity}

This subsection establishes some general continuity results for $C$ near a
regular point $t_0 \in Z$. These are used subsequently here, in
Section~9 and
in the Appendix. The discussion in this subsection is broken into five steps.

\medskip
{\bf Step 1}\qua First, agree to simplify the notation by taking
linear coordinates on $Z$ near the regular point $t_0$ so that $t_0
= 0$.

Observe that the first assertion of Proposition 7.2 gives an open interval
$(-\d, \d)$ of 0 in $Z_R$, and, with $\e > 0$ chosen, a number
$s_\e > 0$ so that
$|z|/\rho < \e$ where $C$ intersects the set where $t \in (-\d, \d)$ and
 where $\rho^2
+ z^2 < s_\e$

With the preceding understood, fix $\e > 0$ but with $\e << 1/100$ and let
$C'$ be a
component of the intersection of $C$ with the set where $t \in
 (-\d, \d)$ and $\rho^2
+ z^2 < s_\e$. On $C'$ , the pair of functions $(t, u \equiv (2
f)^{1/2} = (\rho^2 -
2z^2)^{1/2})$ restrict as coordinates to some disk neighborhood
of all but a countable
set, $\L$, of points of $C'$ . This set $\L$ consists of the singular
points of $C'$
and the critical points of the restriction of $t$. (The critical points
 of $t$ and $u$
are identical.) Note that $\L$ has no accumulation points on $C'$ . The
 fact that $\L$ is
countable with no accumulation points is proved by an argument which
copies, almost
verbatim, the proof of Lemma 7.5. In fact, the only change in the
 argument is that
$\k$ in (7.6) should be taken to be $\k = -(\rho^2 + z^2)^{1/2}/u = 1 +
O(|z|^2/\rho^2)$.

Note that the association of $(t, u)$ to the points in $C'$ defines a
continuous map,
$\pi$ , from $C'$ to the $u > 0$ half of the $(t, u)$ plane. One
then has:

\medskip
{\bf Lemma 8.1}\qua {\sl There exists $s_0 > 0$ and $\d > 0$ and a positive
integer
$p' \leq p$ with the following significance: Suppose $(t, u)
\in (-\d, \d) \x (0, s_0)$ and $(t, u) \not\in\pi (\L)$. Then
$\pi^{-1}(t, u)$
consists of precisely $p'$  points.}

\medskip
{\bf Proof of Lemma 8.1}\qua The number of points in $\pi^{-1}(t_0,
 u_0)$ counts
the intersection number of $C'$ with the $(t, f) = (t_0, u_0^2)$
version of Example 1.6. This intersection number is no greater than $p$ by the
definition of $p$ in Section~4. The fact that this intersection number is
constant on $p'$  follows from the invariance of intersection
number under perturbations.

The functions $(t, u)$ are preferred coordinates on $C'$ for the
following reason: The fact that $C'$ is pseudo-holomorphic
implies (and is implied by) a reasonably simple differential equation for
$\phi$
and $h = \rho^2z$:
\begin{itemize}
\item\quad $h_u=\rho^2\dt u\dt\phi_t$
\item\quad $h_t=-g^2\dt\rho^2\dt u^{-1}\dt\phi_u$.\itemnum{8.1} 
\end{itemize}
Here, the subscripts ``$u$'' and ``$t$'' signify the partial derivative by the
indicated coordinate.

Introduce the function $\nu = h/u^3$ on the subset of $S^1  \x (B^3 -
\{0\}$)
where $|z|/\rho < 1/4$. Where defined, $\nu$ is an analytic functional
of $z/u$
which obeys $\nu = z/u + {\cal O}(|z|^3/u^3)$. Also:
\begin{itemize}
\item\quad $|\nu|\leq\zeta\dt\e \ \ {\mbox{on}} \ \ C'$
\item\quad $\nu \ {\mbox{tends uniformly to zero as $u$ tends to
zero}}$.\itemnum{8.2}
\end{itemize} 
With $\nu$ understood, (8.1) is equivalent to the system:
\begin{itemize}
\item\quad $\nu_u+3\nu/u-(1+\kappa_1)\dt\phi_t=0$
\item\quad $\nu_t+(1+\kappa_2)\dt\phi_u=0$\itemnum{8.3}
\end{itemize} 
Here, $\k_1 = \rho^2/u^2 - 1 = \nu^2 + {\cal O}(\nu^4)$ and $\k^2 = g^2\dt
\rho^2/u^4
- 1 = 8\nu^2 + {\cal O}(\nu^4)$ are analytic functionals of $\nu^2$
where the latter
is less than 1/16.

\medskip
{\bf Step 2}\qua If $C'$ is not smooth and $\pi$  is not
1--1, then $\phi$ and $\nu$ cannot be considered as globally defined
functions of $(t,
u)$ where $t \in (-\d, \d)$ and $u \in (0, s_\e)$. Rather, they must be
thought
of as multivalued functions whose values over any point $(t, u)$ specify the
$p'$ different sheets of $C'$ over the $(t, u)$ plane.

The preceding can be said differently as follows: As remarked, the pair
 $(t, u)$
define local coordinates $C'$ except at the points of $\L$. This
set $\L$ can be assumed to be infinite or empty, for it $\L$ is
finite, one can
decrease $s_0$ to obtain a $C'$ which admits a proper fibration by $\pi$
over. In any event, when $(t, u)\not\in\pi (\L)$, then the
point $(t, u)$ has a
neighborhood that parametrizes $p'$  sheets of $C'$ via a set, $\{(\phi_,
\nu_1),\dots (\phi_{p'},\nu_{p'})\}$ of pairs of functions. Here, no two
pairs in this
set agree at any point of the given neighborhood; and each pair obeys
(8.2) and (8.3)
on the given neighborhood. However, as $(t, u)$ approaches a point in
$\pi(\L)$,
two or more pairs from this set approach the same value.  What is more,
as $(t, u)$ is
circled around any point in $\pi(\L)$, monodromy may permute some
non-trivial subset
of $\{(\phi_,\nu_1),\dots (\phi_{p'},\nu_{p'})\}$.

\medskip
{\bf Step 3}\qua The purpose of this step is to prove that
\begin{equation}
\sup_j |\nu_j|\leq c\dt u
\end{equation}
for some constant $c$. The argument for (8.4) starts with (8.3).
Differentiate the latter to eliminate $\phi$, and so obtain the
following second
order equation for $\nu = \nu_j$ on the complement of $\pi(\L)$:
\begin{equation}
((1+\kappa_1)^{-1} (\nu_u+3\nu/u))_u+ ((1+\kappa_2)^{-1} \nu_t)_t=0
\end{equation}
This is to say that:
\begin{eqnarray}
&& (1+\k_1)^{-1}\dt (\nu_{uu}+3\nu_u/u -3\nu/u^2-2\dt\k'_1\dt
(1+\k_1)^{-1}(\nu\nu_u^2-3\nu^2\nu_u))\nonumber \\
 &&\quad + (1+\kappa_2)^{-1}\dt( \nu_{tt}-\k'_2\dt(1+\k_2)^{-1}
\nu\nu_t^2)=0
\end{eqnarray} 
Now consider (8.6) where $\nu \geq  0$ and where $u$ is such that $\sup_j
|\nu_j|\leq
1/32$. Here, (8.6) implies the differential inequality
\begin{equation}
(1+\k_1)^{-1}\dt (\nu_{uu}+3\nu_u/u -3\nu/u^2+\z\dt\nu)
+(1+\k_2)^{-1}\nu_{tt}\geq 0,
\end{equation}
where $\z$ is a constant which is independent of $u$ and $\nu$.

Note that this last equation holds when $\nu = \nu_j$ for each $j$ as long
as $\nu_j
\geq  0$ and $u$ is small. With this last point understood, introduce the
function
$\nu = \sup \{0, \nu_1, \dots , \nu_j\}$. This function is
Lipschitz and it also
obeys (8.7), albeit in the weak sense.  (In this regard, interpret $\k_1 =
\k_1(\nu)$ and $\k_2 = \k_2(\nu)$.) With the preceding understood, the
point now is to
consider applying the maximum principle to obtain an upper bound for $\nu$.
 For
this purpose the constant $\d$ can be chosen positive and such that:
\begin{itemize}
\item\quad $[-\d,\d]\subset Z_R$,
\item\quad ${\mbox{Equation (8.7) is valid on a neighborhood
 of}}$ 
\item[]\quad$A\equiv \{(t,u):t\in [-\d,\d] {\mbox{ and }} u\in (0,\d)\}$,
\item\quad $0\leq \nu\leq 1/32 \ {\mbox{on a neighborhood of}} \ A$,
\item\quad $|\k_1|+|\k_2|\leq 1/1000 \ {\mbox{on}} \ A$,
\item\quad $\d\leq 1/(1+4\dt\z), \ {\mbox{where $\z$
appears in (8.7)}}$.\itemnum{8.8}
\end{itemize}
With the preceding understood, choose $t_1 \in (-\d/2, \d/2)$ and introduce the
function
\begin{equation}
s\equiv 256\dt\d^{-1}\dt (u-u^2+(t-t_1)^2).
\end{equation}
Note that (8.8) implies that
\begin{itemize}
\item $(1\!+\!\k_1(\nu))^{-1}\!\dt\!
(s_{uu}\!+3s_u/u\!-\!3s/u^2+\z\!\dt\! s)\!+\!(1+\k_2(\nu))^{-1}\!\dt
s_{tt}\leq 0 \ {\mbox{on}} \ A$,
\item $s\geq\nu \ {\mbox{on the boundary of}} \ A$.\itemnum{8.10}
\end{itemize} 
It then follows from the comparison principle that $s \geq  \nu$ at all points
of $A$. And, since $t_1$ in (8.9) can be chosen at random in $[-\d/2,
\d/2]$, it
follows that
\begin{equation}
\sup_j \nu_j\leq 256\dt\d^{-1}\dt u
\end{equation}
at points $(t, u)$ where $t \in [-\d/2, \d/2]$ and where $u \leq  \d$.

A similar argument proves that $\inf_j \nu_j \geq  -256\dt \d^{-1}\dt u$ on
this
same domain.

\medskip
{\bf Step 4}\qua This step serves up the following
technical lemma:

\medskip
{\bf Lemma 8.2}\qua {\sl {Let $\d$ and $A$ be as in} (8.8). {Then the
function
$\sum_j (|d\phi_j|^2 + |d\nu_j|^2)$ is integrable over the subset of
points $(t, u)
\in A$ where $t \in [\d/2, \d/2]$. Moreover, there exists $\z \geq  1$ with the
following significance: Let $s \in (0, \d/8)$ and let $D$ be a disk of
radius $s$
in the $(t, u)$--plane with center at a point $(t_1, u_1)$ with
$\d/8 \leq  t_1 \leq
\d/4$ and $0 \leq  u_1 \leq \d/4$. Let $D' \equiv D \cap  A$. Then}
\begin{equation}
\int_{D'}\sum_j (|d\phi_j|^2+|d\nu_j|^2)\leq\z\dt s^{1/\z} .
\end{equation}}

\med
{\bf Proof of Lemma 8.2}\qua To see that $\sum_j (|d\phi_j|^2 + |d\nu_j|^2)$
is locally integrable near points of $\pi(\L)$ argue as follows: Let $w
 = t
+ iu$ and let $w_0 \in\pi(\L)$. Then, the proof of (7.26) can be repeated
almost verbatim here to find $\z > 0$ such that $\sum_j (|d\phi_j|^2 +
|d\nu_j|^2)
\leq  \z |w - w_0|^{-2+2/\z}$ on some neighborhood of $w_0$.

In order to prove the global integrability of $\sum_j (|d\phi_j|^2 +
|d\nu_j|^2)$,
first reintroduce the cut-off function $c$ as described prior to (2.1).
Then, fix a pair of numbers $s_0, s_1$  which satisfy $0 < s_0 < s_1
\leq
\d/4$. Now, let $\b$ denote the function $\chi(2|t|/\d)\dt \chi(u/s_1)\dt (1 -
\chi(u/s_0))$. Thus, this function vanishes where either $|t| \geq  \d$ or
$u\leq  s_0$ or $u \geq  2s_1$. On the other hand, this function equals 1 where
both $|t| \leq  \d/2$ and $2s_0 \leq  u \leq  s_1$.

Square both sides of (8.3), multiply the result by $\b^2$ and then sum the
result over $j$. Next, integrate this sum over $A$. The result,
after a
judicious application of the triangle inequality and appeals to (8.4) and
the fourth point of (8.9), yields
\begin{equation}
\int_A\b^2\sum_j (|d\phi_j|^2+|d\nu_j|^2)\leq\z
(s_1\d +|\int_A \b^2\sum_j
(\nu_{ju}\phi_{jt}-\nu_{jt}\phi_{ju})|) .
\end{equation}
Here, $\z$ is a constant which is independent of  $s_0$ and 
$s_1$. Now,
integrate by parts on the right-hand side of (8.13) to lift the derivatives
off of $\nu_j$. The resulting right-hand integal reads:
\begin{equation}
-2 \int_A\sum_j (\b_u \nu_j\b\phi_{jt}-b_t\nu_j\b\phi_{ju}) .
\end{equation}
The fact that $\{(\phi_j, \nu_j)\}$ defines a multivalued function on the
complement of $\pi(\L)$ has no effect on the application of integration by
parts. Furthermore, there are no anomalous boundary terms in (8.14) from
the points in $\pi(\L)$.

These last points can be verified with the following argument: To begin,
the data $\{(\phi_j, \nu_j)\}$ can be thought of as the push-forward via
$\pi$
of the restriction to $C$ of the functions $(\phi,\nu)$. Furthermore,
the support of
$\pi^*\b^2$ on $C$ is the image of a smooth manifold $\Sig$ with
boundary via a map,
$f$, which is pseudo-holomorphic on the interior of $\Sig$. This map
identifies
$C - \pi^{-1}(\L)$ with the complement of a finite number of points
in $\Sig$.
Now, note that the integral on the right side of (8.13) is equal to
$\int_{\Sig} f^*
(\pi^*\b^2 d\nu\wedge d\phi)$. Moreover, the functions and forms in this
integral are
smooth on $\Sig$, and have support in $\Sig$'s interior. Thus, Stokes'
theorem applies
and equates this integral with $-2\dt \int_{\Sig} f^*(\nu\pi^*
d\b\wedge\pi^*\b d\phi)$,
which is (8.14).

With (8.14) understood, remark that $|d\b|\dt |\nu_j|$ is uniformly
bounded courtesy of (8.4). Thus, Holder's inequality implies that (8.14) is
no larger in absolute value than
\begin{equation}
\z\dt (s_1\dt\d)^{\frac 12}\dt (\int_A\b^2\sum_j |d\phi_j|^2).
\end{equation}
where $\z$ is independent of both $s_0$ and $s_1$. And, this last
bound on (8.14) can
be substituted into (8.13) to deduce that
\begin{equation}
\int_A\b^2\sum_j (|d\phi_j|^2 +|d\nu_j|^2)\leq\z\dt s_1\dt\d  ,
\end{equation}
where $\z$ is, once again, independent of both $s_0$ and $s_1$.

Since $\z$ is independent of both $s_0$ and $s_1$  in (8.16), one can
take $s_0$ to 0
in (8.16) to conclude that $\sum_j (|d\phi_j|^2 + |d\nu_j|^2)$ is
integrable as
asserted in Lemma 8.2. Moreover, (8.16) implies that when $s'\leq \d/4$, then
\begin{equation}
\int_{|t|\leq\d/4, \ u\leq s'}
\sum_j( |d\phi_j|^2 +|d\nu_j|^2)\leq\z\dt s'\dt\d .
\end{equation}
Now consider the bound in (8.12). In this regard, it proves convenient to
consider separate arguments for the two cases where the center point $(t_1,
u_1)$ satisfies either $u_1 \leq 4s$ or $u_1 \geq  2s$. In
 the former, case,
the integral of $\sum_j (|d\phi_j|^2 + |d\nu_j|^2)$ over $D'$  is
no greater than its integral over the subset where $\d/2 \leq  t \leq
\d/2$ and $u \leq  5s$. The latter integral is bounded in (8.17) using
$s'  = 5s$.

With the preceding understood, consider the case where $u_1 > 2s$. In this
case, $u > 0$ on $D'$  and so $C$ is the image of a smooth,
pseudo-holomorphic map over $D'$. To begin, suppose that $s_1\in
(0, u_1/2]$ and let $D \equiv  D(s_1)$ denote the disk in the $(t,
u)$ plane
which is concentric with $D'$  but has radius $s_1$. Note that the
integral of $\sum_j (|d\phi_j|^2 + |d\nu_j|^2)$ over the radius $s_1 = u_1/2$
disk has been bounded apriori by $\z\dt u_1\dt \d$. Now, consider squaring
 both
sides of (8.3) and integrating over the radius $s_1  < u_1/2$ disk
 $D$. The
result, after an application of Holder's inequality reads
\begin{equation}
\int_D\sum_j ( |d\phi_j|^2 +|d\nu_j|^2)\leq\z\dt
|\int_D \sum_j d\nu_j\wedge d\phi_j|+\z\dt s^2_1 .
\end{equation}
The next step in the proof of (8.12) analyzes the right-hand side of (8.18)
via an integration by parts to obtain the inequality
\begin{equation}
\int_D\sum_j ( |d\phi_j|^2 +|d\nu_j|^2)\leq\z\dt
|\int_{\p D} \sum_{j\in\O'} \nu_j\dt d\phi_j|+\z\dt s^2_1 .
\end{equation}
(The arguments which justified (8.14) also justify (8.19).)

Further analysis of the right-hand side of (8.19) requires the following
claim: There is a constant $\z$ which is independent of $s_1$  and of the
chosen
center point of $D$ such that the right side of (8.19) is no larger than
\begin{equation}
\z\dt(s_1\dt\int_{\p D}\sum_j ( |d\phi_j|^2 +|d\nu_j|^2)+s^2_1).
\end{equation}
when $\p  D \cap  \L = \emptyset$. This last claim will be justified
momentarily, so
accept it for the time being. Then, let $f(s_1)$ denote the integral on the
right side of (8.19). Coupled with the preceding claim, (8.19) implies that
$f$ satisfies, for all but finitely many values of $s_1$, the differential
inequality
\begin{equation}
f\leq\z\dt s_1\dt f'+\z\dt s^2_1 .
\end{equation}
Here, $f'$  denotes the derivative of $f$. This last
equation can be integrated directly to find that $f(s_1) \leq
\z\dt s_1^{1/\z}$ where $\z$ is independent of $s_1$  and of the center
point of $D$.
(The proof of this last point requires the apriori bound for $f(u_1)$.)

With the preceding understood, the justification of the claim that (8.20)
bounds the right-hand side of (8.19) completes the proof of Lemma 8.2. To
start this justification, note that as $\p  D$ is assumed to be disjoint from
$\L$, the map $\pi$  from $C$ to the $(t, u)$ plane is a $p$--fold
covering map near
$\pi^{-1}(\p  D)$ and so $\pi^{-1}(\p  D)$ is a smooth, embedded 1--manifold
in $C$.
Let $\{\G\}$ denote the set of components of $\pi^{-1}(\p  D)$. Then, the
integral on the right side of (8.19) is equal to
\begin{equation}
\sum_{\G\in\{\G\}}\int_{\G}\nu\wedge d\phi .
\end{equation}
Now, each boundary component $\G$ is homologically trivial in the
$\rho > 0$
subspace of $S^1\x {\R}^3$ as each deforms to a multiple cover of
$\p D$.
 This means that
\begin{equation}
\int_\G\nu\wedge d\phi=\int_\G (\nu-\a)\wedge d\phi
\end{equation}
where $\a$ is any constant. In particular, take $\a = (2\pi \dt
p_\G\dt s_1)^{-1}
\int_\G \nu\pi^*(d\s)$, where $p_\G$  is the degree of the covering map
$\pi\co\G
\to  \p  D$ and where $d\s$ is the induced Euclidean line
element on $\p D$.
(Thus, $\int_{\p  D} d\s = 2\pi s_1$.) Note that $\pi^*d\s$ trivializes
$T\G$ and so
defines a Riemannian metric for $\G$.  With respect to this metric, there is a
standard eigenvalue estimate which finds that $\int_\G (\nu-\a)^2 \pi^*d\s
\leq
s_1^2 p^2\dt \int_\G |d\nu|^2 \pi^*d\s$. Then, with the proceding understood,
(8.23) implies that
\begin{equation}
\sum_{\G} |\int_{\G}\nu\wedge d\phi| \leq \z\dt s_1\dt\sum_\G
(\int_{\G} |d\nu|^2\pi^*d\s)^{\frac 12} (\int_\G |d\phi|^2\pi^*d\s)^{\frac
12} .
\end{equation}
Meanwhile, the right-hand side of (8.24) is equal to 
$$\z\dt s_1 \dt \sum_j
(\int_{\p  D} (|d\nu_j|2)^{1/2} (\int_{\p  D} |d\phi_j|^2)^{1/2}.$$
Thus, the claim prior to (8.20) follows from (8.24).

\medskip
{\bf Step 5}\qua Let $\d$ be as in (8.8). The set
$\{(\phi_j, \nu_j)\}$
defines a continuous map, $\Phi$, from the set $A'  = \{(t,
u): t \in
[-\d/8, \d/8]$ and $u \in (0, \d/8]\}$ into the $p$'th symmetric
product of
$\R^2$. This concerns the extension of $\Phi$ to the $u = 0$ line in the
closure of $A'$:

\medskip
{\bf Lemma 8.3}\qua{\sl {The map $\Phi$ extends as a continuous map to
the $u = 0$ line
in the closure of $A'$. Furthermore, the extension is Holder
continuous in the sense that there exist constants $\a > 0$ and $\z$ such that
if $w, w'  \in [-\d/8, \d/8]\x [0, \d/8]$ then the distance from
each $(\nu,\phi) \in \Phi(w)$ to $\Phi(w' )$ and from each $(\phi,\nu)\in
\Phi(w' )$ to $\Phi(w)$ is no greater than} $\z\dt |w - w'|^{\a}$}.

\medskip
{\bf Proof of Lemma 8.3}\qua First, introduce the function $k = \sup_j |\phi_j|$.
Then it follows from Lemma 8.2 that $|dk|^2$ is integrable over $[-\d/4,
 \d/4]$
and that there is a constant $\z \geq  1$ with the following significance:
Let $s, \ D$ and $D'  = D \cap  A$ be as in the statement of Lemma
8.2. Then
\begin{equation}
\int_{D'} |dk|^2\leq z\dt s^{1/\z} .
\end{equation}
With (8.25), a minor modification of the proof of Theorem 3.5.2 in \cite{Mo}
proves that the function $k$ extends to the closure of $A'$  as a
Holder continuous function in the sense that $|k(w) - k(w' )|
\leq \z' \dt |w - w' |^{1/z'}$
when $w, w'$  lie in the closure of $A'$. Here,
$\z'$  is a fixed constant. (Theorem 3.5.2 in \cite{Mo} does not
consider the case where $k$ is defined on a half space. However, the
proof is
easily modified to handle this case too. Alternately, one can define $k$ on
$[-\d/4, \d/4]\x [-\d/4, -\d/4]$ by declaring that when $u < 0$, then $k(t,
u) =
k(t, -u)$. This extension of $k$ will obey (8.25) where $D'$  is now
any disk of radius s which lies entirely in $[-\d/4, \d/4]\x [-\d/4, \d/4]$.
Then, Theorem 3.5.2 in \cite{Mo} applies without modification.)

Now consider the map, $\Psi$, from Sym${}^n({\Bbb C})$ to ${\Bbb C}^n$
which sends
the unordered set of complex numbers $\{\l_1, \dots , \l_n\}$ to the
coefficients $(a_1, \dots , a_n)$ of the monic, $n$'th order
polynomial whose roots
are $\{\l_1, \dots  , \l_n\}$. (Thus, $a_j$ is a homogenous, $j$'th order
polynomial
in the variables $\{\l_j\}$. For example, $a_1 =\sum_j \l_j$ and $a_n =
\prod_j
\l_j$.) This map is a homeomorphism with Holder continuous inverse.
This is to say that there exists $\z \geq 1$ (which depends on $n$)
and is such that if $a, a' \in {\Bbb C}^n$, then the distance from each
$\l \in\Psi^{-1}(a)$ to $\Psi^{-1}(a' )$ and from each $\l
\in\Psi^{-1}(a' )$ to $\Psi(a)$ is no greater than $\z\dt |a - a'
|^{1/\z}$.

With the preceding understood, the association to each point in
$A'$  of the element $\Psi(\{\phi_j\}) \in {\Bbb C}^n$ defines a
map, $\Psi_A$, from $A'$  to ${\Bbb C}^n$. Because $\sup_j \{|\phi_j| +
|\nu_j|\}$ is bounded and because $\Psi$ is polynomial in the $\{\phi_j\}$, it
follows from Lemma 8.2 that $\Psi_A$ is a Sobolev class $L^2_1$ map.
Furthermore,
Lemma 8.2 implies that each component of $\Psi_A$ obeys (8.25) for an
 appropriate
choice of $\z$. With this understood, the proof of Theorem 3.5 in \cite{Mo}
can be
readily modified to prove that $\Psi$ extends to the closure of $A'$
as a Holder continuous function. Finally, since $\Psi^{-1}$ is Holder
continuous,
Lemma 8.3's claim follows.

\subsection*{8.b\qua Assertions 3(b), 3(c) and the $p = 1$ case of Assertion 3(a) of
 Proposition 7.1}

The purpose of this subsection is to establish the $p = 1$ version of
Assertion 3(a) of Proposition 7.1 as well as Assertions 3(b) and 3(c) of
Proposition 7.1 under the assumption that $C$ intersects a neighborhood in
$S^1\x (B^3
- 0)$ of the regular point $t_0$ as the disjoint union of finite energy, $p
= 1$
subvarieties.

Here is the argument: Under the given assumption above, $C$ intersects a
neighborhood of the point $t_0$ in $S^1\x (B^3 - 0)$ as $p$ disjoint
components,
and so it is sufficient to focus attention on any one component. The latter
 is
described as in Step 1 of subsection 8.a by a pair of functions $(\phi,\nu)$
of the
variables $(t, u)$ subject to (8.2) and (8.3) for some suitably small
choice
of the constant $\e$. Furthermore, it follows from (8.4) that $|\nu| \leq
c\dt u$ for some constant $c$; and it follows from Lemma 8.3 that
 $(\phi,\nu)$
extend to the $u = 0$ line as Holder continuous functions for
some positive
Holder exponent.

Note that the function $\nu$ obeys (8.5), and by virtue of (8.3), the function
$\phi$
obeys the equation
\begin{equation}
u^{-3}((1+\k_2)\dt u^3\phi_u)_u + ((1+\k_1)\dt \phi_t)_t=0.
\end{equation}
Surprisingly enough, (8.5) and (8.26) can be thought of as equations for
functions on a domain in $\R^5 = \R\x \R^4$ as follows:
Let $(v_0, v_1, v_2, v_3, v_4)$ be Euclidean coordinates on ${\R}\x {\R}^5$
where $v_0$ corresponds to the first factor of $\R$. Then, map
$\R\x \R^4$ to the $(t, u)$ half plane by writing $t = v_0$ and $u
= (v_1^2 +
\dots v_4^2)^{1/2}$. The pull back of $\phi$ and $\nu$ by this map gives
functions $\underline{\phi}$ and $\underline{\nu}$, respectively, on
${\R}\x{\R}^4$. It also proves convenient to introduce the angular
 coordinate
$\th\in  [0, \pi]$ on the $\R^4$ factor of $\R^5$ by declaring that $\cos\th
= v_1/u$. With $\th$ understood, set $\underline{\l}
\equiv\underline{\nu}\cos\th$.

A domain, $D'$, for $\underline{\phi}, \underline{\nu}$ and
$\underline{\l}$ can be
taken to be the product of an open line segment centered at zero in
the $\R$
factor with the complement of the origin of some small ball in the $\R^4$
factor. Let $\e_1$ denote the radii of this segment and this ball.
Let $D\subset
D'$  now denote the product of this same interval with the whole
ball in the $\R^4$ factor.

Now consider:

\medskip
{\bf Lemma 8.4}\qua {\sl {Let $\{\nabla_\a\}_{\a=0,\dots 4}$ denote the
derivatives with respect to the coordinates $\{v_\a\}_{\a=0,\dots
4}$ on
${\R}\x
\R^4$. Then $\underline{\phi}, \ \underline{w}$ and $\underline{\l}$ all
extend from
$D'$  to $D$ as Holder continuous (with some positive exponent) and Sobolev
class
$L^2_1$ functions. Furthermore, $\underline{\phi}$ and
$\underline{\l}$ obey the
equations}
\begin{itemize}
\item\quad $\displaystyle\D\underline{\phi}+\sum_{\a,\b}\nabla_\a
 (a^{\a\b}\nabla_\b\underline{\phi})=0$
\item$\displaystyle\quad \D\underline{\l}=\nabla_1\underline{\k}_1
\nabla_0\underline{\phi}-\nabla_0\underline{\k}_2\nabla_1
\underline{\phi}$.\itemnum{8.27}
\end{itemize} 
{Here, $\D = \sum_\a \nabla_\a\nabla_\a$ is the standard
Laplacian on $\R^5$ and
$a^{\a\b}$ is a diagonal matrix with $a^{00} = \k_1$ and with $a^{\a\a} =
\k_2$ for}
$\a > 0$.}

\medskip
{\bf Proof of Lemma 8.4}\qua The fact that these functions are Sobolev class
$L^2_1$ follows directly from Lemma 8.2. Holder continuity follows from
Lemma
8.3. The equation for $\underline{\phi}$ is simply (8.26), while that for
$\underline{\l}$ follows from (8.5) with two extra observations:
 First, $\D(\cos\th) =
-3\dt \cos\th/u^2$. Second, if $f$ is a function on $\R\x \R^4$
of only $t$ and $u$, then $\nabla_1f = f_u\dt \cos\th$.

The equations in (8.27) for $\underline{\phi}$ and $\underline{\l}$
 and their initial
$L^2_1$ and Holder bounds can now be used with fairly standard (though tedious)
elliptic regularity arguments to prove that the latter are analytic
functions on the
domain $D$ in ${\R}\x \R^5$. (The elliptic regularity arguments are
essentially the same as those of Morrey which yield Theorems 6.2.5 and 6.7.6
in \cite{Mo}.)

The analyticity of $\underline{\phi}$ and $\underline{\l}$ imply directly
Assertion 3(b)
of Proposition 7.1 for the original functions $(\phi, z)$. Then, the
analyticity of
$(\phi, z)$ and the equations in (8.3) give directly Assertion 3(c) of
 Proposition
7.1. This analyticity also implies Assertion 3(a) in the $p = 1$ case. Indeed,
for $\d > 0$ but small, the map $(-\d, \d)\x (0,\d)$ to $S^1\x (B^3
- \{0\})$
which assigns the coordinates $(t, x = r \cos\phi, y = \rho
\sin\phi , z)$ where $\phi
= \phi(t, u)$, $z = z(t, u)$ and $\rho = 2^{1/2} u\dt (1 +
 z^2/u^2)$ defines an
embedding which is 1--1 and onto the intersection of $C$ with
 neighborhood of $t_0$.

\section{The number of singular points}
\setcounter{equation}{0}

When studying the behavior of $C$ near the singular set $Z_S\subset
Z$, it is
sufficient, after Propositions 5.1 and 7.1, to restrict attention to the
case where $C$ is a pseudo-holomorphic subvariety with $n_{\pm}  = 0$
at all
points of $Z$ and with $q_{\pm}  = 0$ at all regular points of
$Z$. Furthermore,
one can assume that C has no irreducible components which lie in the $z = 0$
half space as all of these are described by Example 1.3.

The results in this section are summarized by:

\medskip
{\bf Proposition 9.1}\qua{\sl The set of singular points is finite.}

\medskip
The proof of this proposition occupies the remainder of this section.

\medskip
{\bf Proof of Proposition 9.1}\qua The proposition is an immediate corollary of
Proposition 5.2 and Proposition 9.2, below. First, Proposition 5.2 produces
a uniform bound by $\z\dt s^3$ for the function $\mu(s)$ in (5.1). Then,
Proposition 9.2 implies the following claim:
\begin{itemize}
\item Let $t_0\in Z_s$ and let $\nu >0$. If $s$ is small, then the subset
of $C$ where $t\in [t_0-\nu, \ t_0+\nu]$ contributes at least
$\z^{-1}\dt s^3$ to the function $\mu$. Here, $\z\geq 1$
is independent of $t_0$, $\nu$, and $s$. \hspace{\fill} (9.1)
\end{itemize} 
\addtocounter{equation}{1}

\noindent The bound on $\mu$ from Proposition 5.2 plus (9.1) imply that
$Z_S$ is
finite.

As just remarked, Proposition 9.2 leads to (9.1); however a brief digression
is required in order to state this proposition. To start the digression,
agree to simplify the notation by taking $t_0$ to be the origin, 0, in a
Euclidean coordinate system for a neighborhood of $t_0$ in $Z$. Having
done so,
fix $\nu > 0$ and small so that $Z$ has neither $t =
 -\nu$ nor $t = \nu$ as singular
points. When $r_1$ is positive, but small (much less that $\d$ in
(5.1)); and when
$s \in (0, r_1)$, let $D(s)$ denote the subset of $S^1\x \R^3$ where the
following constraints hold:
\begin{itemize}
\item $ \quad t\in [-\nu,\nu]$
\item $ \quad \rho\leq r_1$
\item $ \quad |z|\leq r_1$
\item $ \quad 0\leq h\leq s$\itemnum{9.2}
\end{itemize}
Remember that when $r_1$ is generic and $s$ is small relative to $r_1$,
 then $C$ has
empty intersection with the boundary of $D(s)$ where $|z| = r_1$. Thus, the
intersection of $C$ with the boundary of $D(s)$ consists of the regions
where:
\begin{itemize}
\item $ \quad h=s^3 \ {\mbox{while}} \ \rho\leq r_1, \
  |z|\leq r_1 \ {\mbox{and}} \ t\in [-\nu,\nu] $
\item $ \quad h=0 \ {\mbox{while}} \ \rho\leq r_1, \
  |z|\leq r_1 \ {\mbox{and}} \ t\in [-\nu,\nu] $
\item $ \quad  \rho= r_1 \ {\mbox{while}} \
  |z|\leq s^3/r^2_1 \ {\mbox{and}} \ t\in [-\nu,\nu] $
\item $ \quad t=\pm\nu \ {\mbox{while}} \ \rho\leq r_1, \
  |z|\leq r_1 \ {\mbox{and}} \ 0\leq h\leq s^3$\itemnum{9.3}
\end{itemize}
Use $\p  C$ to denote the intersection of $C$ with the boundary of $D(s)$.
Also,
use $\p' C$ to denote the intersection of $C$ with the part of the
boundary of $D(s)$ where $t=\pm\nu$ and use $\p''C$ to denote the
 intersection of
$C$ with the part of the boundary of $D(s)$ where $\rho = r_1$.

Now consider the function
\begin{equation}
\mu_0(s)=\int_{C\cap D(s)} d\phi\wedge dh.
\end{equation}
Integration by parts can be used to equate
\begin{equation}
\mu_0(s)=-s^3\int_{\p C\cap \{h=s\}} d\phi-
\int_{\p' C}h\dt d\phi .
\end{equation}
End the digression and consider the following proposition:

\medskip
{\bf Proposition 9.2}\qua {\sl {Suppose that 0 is in $Z_S$ and that $q_+(0) >
0$. When
$r_1 > 0$ is small and when $s > 0$ is very small, the following is true:}

\begin{itemize}
\item{The first term in} (9.5) {is greater than} $s^3$.

\item{The second term in} (9.5) {is bounded in
absolute value by
$\z' \dt r_1\dt s^3$. Here, $\z'$  is independent of $s$
and} $r_1$.

\item{The third term in} (9.5) {is bounded in absolute
value by $\z\dt
s^6$. Here, $\z$ is independent of $s$ but may depend on} $r_1$.

\item{Thus, $\mu_0(s) > s^3/2$ for small $r_1$ and} $s$.
\end{itemize}}

\medskip
Note that the fourth assertion of Proposition 9.2 implies (9.1) because when
$r_1$ is small and $s$ is small, then the set $\O(s)$ in (5.1) will
contain $D(s)$.

\medskip
{\bf Proof of Proposition 9.2}\qua The proof of this proposition is quite
lengthy and so is broken down into steps, fourteen in all.

\medskip
{\bf Step 1}\qua This step proves the second assertion
of Proposition
9.2. Here is the argument: Since $\pm\nu$ are regular points, Proposition 7.1
finds $r_0 > 0$ and a constant $\z$ such that $|d\phi| \leq
\z$ on $\p C'$ when $r_1
\leq  r_0$. Thus, since $h \leq  s^3$, the second assertion of
Proposition 9.2 follows
by integration.

\medskip
{\bf Step 2}\qua This step states and proves Lemma 9.3,
below, which
implies the third assertion of Proposition 9.2. Note that Lemma 9.3 also
plays a pivotal role in the proof of Proposition 9.2's first assertion.

\medskip
{\bf Lemma 9.3}\qua {\sl {There exists $r_0 > 0$ such that the following is true:
Given $r_1\in (0, r_0]$, let $I(r_1; s) = \int_{\p''C}|d\phi|$.
 Then, for all but finitely many values for} $r_1\in [r_0/2, r_0]$:

\begin{itemize}
\item{The function} $s^{-3}\dt I(r; s)$ {is bounded as
s tends to zero.}

\item{The quantity} $s^{-3}\dt$ length$(\p\,'' C)$ {is bounded
as} $s \to  0$.

\item{The number of components of $\p\,''C$ becomes constant
as $s \to  0$.
Furthermore, each component of $\p\,''C$ is an interval on which $dh \neq
0$ when $s$
is small.}
\end{itemize}}

\medskip
{\bf Proof of Lemma 9.3}\qua Since $|d\phi| \leq  r_1^{-1}$ on $\p\,'' C$, the
size
of
$I(s)$ will be ${\cal O}(s^3)$ when length$(\p\,'' C)$ is also ${\cal
O}(s^3)$. Thus,
the first assertion follows from the second.  A uniform bound for
 $s^{-3}\dt$
length$(\p\,'' C)$ will occur (for small $s$) if $C$ intersects
the constraint surface where $r = r_1$ and $h = 0$ in smooth
 points, and if $dh \neq
0$ at all points of this intersection.

The fact that $r_1$ can be chosen to make the latter true follows from:

\medskip
{\bf Lemma 9.4}\qua {\sl{There are at most finitely many singular points
of $C$ in any compact subset of $C$. In addition, there are at most
finitely many critical points of the restriction of $h$ to the smooth
part of any compact subset of} $C$.}

\medskip
Given Lemma 9.4, one can choose $r_1$ from the complement of a finite
set to insure that $dh \neq 0$ where $h = 0$ and $r = r_1$. This makes
both the second and third points of Lemma 9.3 true. (The final point
follows from the implicit function theorem when $dh \neq 0$ at the
points in $C$ where $h = 0$ and $r = r_1$.)

\medskip
{\bf Proof of Lemma 9.4}\qua The first assertion about the singularities of
$C$ follows from general regularity results about pseudo-holomorphic
submanifolds. Indeed, one can parametrize $C$ using the smooth model curve
$C_0$
and the pseudo-holomorphic map $f\co C_0\to C$ as in Definition
1.1. Then, the result follows from general regularity results about $f$ as
in \cite{PW}, \cite{Ye} or \cite{MS}.

As for the statement about $h$, consider $f$ and $C_0$ as above and consider
the pull-back $f^*h$. This satisfies an elliptic equation on $C_0$ with the
schematic form $d^*d(f^*h) + a\dt d(f^*h) = 0$, where $a$ is an
appropriate vector field on $C_0$. Indeed, because $f$ is
pseudo-holomorphic,
this follows from the fact that the restriction of $h$ to the smooth
part of $C$ obeys
the equation $d((g\rho^2)-1\dt Jdh = 0$. It then follows using the unique
continuation
theorem of Aronzajn \cite{Ar} as applied to $f^*h$ that the latter can have
only finitely
many critical points in any compact subset of $C_0$. (Argue here as in the
proof of
Lemma 7.5.)

\medskip
{\bf Step 3}\qua This step begins the proof of the first
assertion of
Proposition 9.2. In this regard, the simplest case occurs when the
intersection of $\p C$ with the surface $h = s^3$ contains a
 closed loop. Then,
the first term in (9.5) is non-negative and has the form $2\pi s^3\dt m$,
where $m$ is a positive integer because $(2\pi)^{-1}\dt d\phi$
restricts to $\p C$ to
define an integral valued cohomology class. (The possibility that $m =
0$ is ruled out
with the observation that in this case, either $q_+ = 0$ or else $C$ is
identical to a
surface from either Example 1.3 or 1.4.)

Note that the argument for the first assertion of Proposition 9.2 in the
general case also involves a cohomological argument using $d\phi$. However, in
the general case, these cohomological arguments are coupled with certain
analytic results. The first of these analytic results concerns some special
properties of the level sets of $h$ on $C$ and is discussed here.

To begin, choose $r_1 > 0$ and small so that the conclusions of
Lemma 9.3
hold. Fix $s'  \in (0, s]$ for the moment and consider the
subset of $C \cap  D(s)$ where $h = {s'}^3$. When $s'$  is
chosen from a certain dense and open set, then $h^{-1}({s'}^3) \cap
(C \cap  D(s))$ will lie in the smooth portion of $C$ and it will miss all
critical points of $h$. In this case, the set
\begin{equation}
\g(s')=h^{-1}({s'}^{
3}\cap (C\cap D(s))
\end{equation}
is a smoothly embedded curve.
Furthermore, if $dh$ is used to orient the normal bundle to $\g(s' )$, then the
restriction of $d\phi$ to $\g(s' )$ is nowhere vanishing and negative
definite. Thus,
one has the following crucial observation:
\begin{itemize}
\item The integral of $-d\phi$ over $\g(s')$ is bounded from
below by its integral over any subset of $\g(s')$.\hspace{\fill} (9.7)
\end{itemize} 

\medskip

{\bf Step 4}\qua The following lemma is a key corollary
to Lemma 9.3:

\medskip
{\bf Lemma 9.5}\qua{\sl {Suppose that $r_1 > 0$ is chosen as in Lemma 9.3, but
small;
and suppose that $s > 0$ is very small. Let $s'  \in (0, s]$ be
such that $-\g(s' ) d\phi > 5/9$. Then $-\int_{\g(s')} d\phi > 1/2$ for
 all $s''
\in (0, s]$. In particular, $\mu_0(s'')$ from (9.5)
satisfies $\mu_0(s'') > {s''}^{3}/2$ for all such} $s''$.}

\medskip
{\bf Proof of Lemma 9.5}\qua The proof starts by connecting the curve
$\g(s' )$ to $\g(s'' )$ to form a closed
curve. However, this task requires a digression to consider the endpoints of
$\g(s' )$. To start the digression, note that $\g(s')$
will have endpoints in $\p\,'C$ and in $\p\,'' C$. To consider the endpoints of
$\g(s')$ in $\p\,' C$, remark first that the $t = \nu$ portion of $\p\,'C$
consist of
some number $p \geq  0$ of arcs, each with one endpoint in $Z$ and
the other endpoint
where $h = s^3$. Furthermore, if $r_1$ is small and s is very small, then
the $h =
s^3$ endpoint of an arc component of $\p\,'C$ does not sit in $\p\,''C$ as
$C$ is
assumed to have no $z = 0$ components. With $\p\,'C$ understood, then the
boundary of
$\g(s')$ in the $t = \nu$ portion of $\p\,'C$ will consist of $p$
points, each on one of
the corresponding arc components of $\p\,' C$. Thus, the endpoints of
$\g(s')$ and
$\g(s'' )$ in the $t = \nu$ portion of $\p\,'C$ are paired by the
condition that they
lie in the same component arc of $\p\,'C$. Then, the intervening segment of this
component arc can be used to connect these paired points.

Now consider the endpoints of $\g(s')$ in $\p\,''C$. As the latter is a
collection of
arcs where $dh \neq   0$, there is a point of $\g(s') \cap  \p\,''C$ in
each such
arc, and thus a natural pairing between the points of $\g(s') \cap  \p\,''C$ and
$\g(s'' ) \cap  \p\,''C$ with an arc in $\p\,''C$ connecting the paired
points.

With the preceding understood, connect the endpoints of $\g(s')$
in $\p\,'C$ with those of $\g(s'')$ by
arcs in $\p\,'C$ and likewise connect the endpoints of
$\g(s')$ with those of $\g(s'')$ by arcs
in $\p\,''C$. Think of the result as a 1--cycle in
$C$. (This cycle traverses $\g(s')$ in the positive direction and
$\g(s'')$ in the negative direction.) This 1--cycle
in $C$ is null homologous as one can consider the homology which moves
$s'$  to $s''$. Thus, $-d\phi$ has integral
zero over this cycle. Then, the contribution to the integral of $-d\phi$
over the
portion of the 1--cycle in $\p\,'C$ is uniformly small (${\cal O}
(r_1)$, courtesy of Proposition 7.1), and that over the portion in
$\p\,''C$ is also small $({\cal O}(s)$, courtesy of
Lemma 9.3). Thus, the contributions to the integral from $\g(s')$
and $\g(s'')$ must be almost the same.

\medskip
{\bf Step 5}\qua The first assertion of Proposition 9.2
follows
immediately using Lemma 9.5 with the following:

\medskip
{\bf Lemma 9.6}\qua  {\sl {Choose $r_1 > 0$, but small so that the
 conclusions of Lemma
9.3 hold. Suppose that $s$ is sufficiently small. Then, there exists
$s'  \in (0, s]$ and an arc, $\g$, which lies in $C  \cap  D(s)
 \cap
\{h = {s'}^{3}\}$ and which obeys $-\int_{\g} d\phi > \pi/2$. Then,
because $-d\phi > 0$ on the set $C \cap \{h = {s'}^{ 3}\}$, it
follows that $-\int_{\g(s')} d\phi > \pi/2$ as well.}}

\medskip
{\bf Proof of Lemma 9.6}\qua The proof of this lemma occupies the remaining
steps in the proof of Proposition 9.2.

\medskip
{\bf Step 6}\qua This step constitutes a digression of sorts
 to state an
auxilliary result which will be used later:

\medskip
{\bf Lemma 9.7}\qua {\sl {Either $h =$ constant on $C$, or else there
 exists $N < \i$
with the following significance: Let $X$ denote the subset of $S^1\x
(B^3-\{0\})$
where $(x^2 + y^2 + z^2)^{1/2} < {1/2}$. Then $N$ bounds the
intersection number
of $C$ with any pseudo-holomorphic curve $S$ on which $h$ is constant.}}

\medskip
{\bf Proof of Lemma 9.7}\qua The proof starts with a standard observation:
\begin{itemize}
\item Let $K$ be a compact subset of $S^1\x (B^3-\{0\})$.
Then $C$ has only finitely many intersections with any
$h=$ constant pseudo-holomorphic curve $\Sigma$ in $K$.
 \hspace{\fill} (9.8)
\end{itemize}\addtocounter{equation}{2}
One can use (9.8) to draw the following conclusion: Given $\e > 0$, there
exists $N_\e$ such that $C$ has intersection number $N_\e$ or less in $X$ with
any $h =$ constant pseudo-holomorphic curve where $|h| > \e$. So,
at issue are
the $h =$ constant curves with small $|h|$.

To consider the case of an $h =$ constant curve with small $|h|$,
first use
(9.8) with the fact that $h$ is continuous on $C$ to draw the following
conclusion: Unless $h \equiv  0$ on $C$, there exists $\e > 0$ and an open
interval $I \subset ({1/2}, 1)$ such that $|h| > \e$ at all points of $C$
where
$r\equiv
(x^2 + y^2 + z^2)^{1/2} \in I$. This last remark is equivalent to the
assertion that $C$ has empty intersection with all curves $\Sig$ with
$h =$
constant at points where $r \in I$ if the constant in question has
norm less
than $\e$.

With these last points understood, use (9.8) with the preceding remark to
conclude that unless $h =$ constant on $C$, then $C$ has constant
intersection
number with all pseudo-holomorphic curves $\Sig$ where $h =$ constant $\in (0,
\e)$. Conclude likewise that unless $h =$ constant on $C$, then $C$ has
constant
intersection number with all pseudo-holomorphic curves where $h =$ constant
$\in (0, -\e)$. Thus, the only $h =$ constant curve where $C$ might have
unbounded intersection number is with some $h = 0$ pseudo-holomorphic curve.
However, if $C$ has more than $N$ intersection points with such a curve,
then $C$
will have more than $N$ such with any $h =$ constant
pseudo-holomorphic curve if
$|h|$ is sufficiently small.

\medskip
{\bf Step 7}\qua Now return to the proof of Lemma 9.6.
This step
summarizes some consequences of the results from Sections 3 and 4 together
with:

\medskip
{\bf Lemma 9.8}\qua {\sl{Given $\e > 0$ and $\a > 0$, there exists $s_0 > 0$ such
that when $s \in (0, s_0)$, then the following is true: The intersection of $C$
with the set of points where $t^2 + (\rho^2 + z^2) \leq  s^2$ and $\a|t|
\leq
(\rho^2 + z^2)^{1/2}$ is the disjoint union $C_{\rho}\cup C_z$, where:}
\begin{itemize}

\item $|z|/\rho < \e$ \ {on} \ $C_{\rho}$ {and}
$\rho/|z| < \e$ on
$C_z$.

\item {If $C'$ denotes either $C_{\rho}$ or $C_z$, then} $\int_{C'}
d\phi\wedge  dh\leq  \e s^3$ .

\item {The functions $t$ and $u = (-f)^{1/2}$
 restrict to a neighborhood
of all but a set of at most countably many points in $C_z$ as local
coordinates. This countable set, $\L$, consists of the singular points
of $C_z$
and the critical point of $dt$ (or equivalently $du$).  This set has no
accumulation points in} $C_z$.

\item {The map $\pi$ to the $(t, u)$ plane given by
restriction of
$t$ and $u$ to $C_z$ has the following additional property: Given
$\e$, $\a$ and
$s$, there exists $\nu_0$ such that when $\nu \in (0, \nu_0)$, then $\pi$  maps
onto the set $A$ where $\a|t| < u < \nu$. In addition, when $(t, u)\in A - A
\cap\pi (\L)$, then $\pi^{-1}(t, u)$ has precisely $q_+(0)$ elements.}
\end{itemize}}

\medskip
{\bf Proof of Lemma 9.8}\qua The first two assertions follow directly from
Proposition 3.3. The third assertion is proved by copying verbatim the proof
of Lemma 7.5. The existence of $\tau_0$ which makes the fourth assertion true
follows from Propositions 3.2 and 4.1.

\medskip
{\bf Step 8}\qua Let $C_z$ be as in Lemma 9.8 with $\e$,
$\a$ and
$s$ all small. (Some precise bounds will be given subsequently.) Then, by
definition, the functions $t$ and $u = (z^2 - \rho^2/2)^{1/2}$
restrict as
coordinate functions near any point of $C_z - \L$. This is to say that each
such
smooth point has a neighborhood which is specified by writing
\begin{itemize}
\item $ \quad x=x(t,u), \ y=y(t,u)$
\item $ \quad z^2=u^2+2^{-1}\dt (x^2+y^2).$\itemnum{9.9} 
\end{itemize}

As in the proof of the second assertion of Proposition 7.2, the ${\Bbb C}$--valued
function $\eta = x + i\dt y$ is constrained to obey (7.9) with $a$ and $b$
obeying the bounds $|a| + |b| \leq  \z\dt |\eta|^2/u^2$. In particular,
according to Lemma 9.8, one has $|a| + |b| \leq  \e$ when $\e$ is small.

As with the discussion in Section~7, the parametrization in (9.9) is only
local in the following sense: The association of the functions $(t, u)$ to
points in $C_z$ defines a map $\pi$  to the $(t, u)$ plane which is
$q'$  to 1 where the number $q' = q_+(0)$. Thus, where
(9.11) holds, $C$ is at best described by a multivalued map from the
$(t, u)$
plane to ${\Bbb C}$. Put differently, the portion of $C$ where (9.11) holds is
described globally by a map from the indicated domain in the $(t, u)$ plane
into the Sym${}^{q'} {\Bbb C}$. The coordinates of the latter map
consist of a $(t, u)$ dependent, unordered set of complex numbers
 $\{\eta_j\}$;
and each $\eta\in\{\eta_j\}$ is a smooth function which obeys (7.9)
where it
is unambiguously defined. The set $\{\eta_j\}$ describe the sheets of $C_z$
over
the image of $\pi$ . Alternately, one can think of $\{\eta_j\}$ as the
$\pi$--push-forward of the function on $C$ which assigns to the value of $x +
i\dt y$ to a point.

With the preceding understood, let $\th\equiv\prod_j \eta_j$. The domain of
$\th$ can be taken to be any small $\a, \ \e$ and $\nu$ version of the
subset $A$ of Lemma 9.8. Here, $\th$ is a continuous function which
is smooth
away from $\pi(\L)$.

As in Section~7, when $\a, \ \e$ and then $\nu$ are chosen and small, the
function $\th$ obeys (7.12) and  (7.13) on $A$. Moreover, the proof of
Lemma
7.7 can be repeated more or less verbatim (using Lemma 7.8 to control $m$
instead of Lemma 7.6) to write $\th$ as
\begin{equation}
\th=\k^{-1}\th\,'
\end{equation}
where $\th\,'$  obeys $\p  \th\,'  = 0$ on its domain of
definition; and where $\k$ is described in the three points of Lemma
7.7. (As
in Section~7, $\p   = 2^{-1}\dt (\p_t - i\cdot \p_u)$.)

To be more precise, Lemma 7.7 writes $\k = u^{q' /2}\dt e^{\tau}$, where $\tau$
obeys
\begin{equation}
\p\tau +i\dt\chi_*\dt m/4u=0  ,
\end{equation}
where $\chi_*$ is a suitable cut off function with values in [0, 1] which
equals one on the domain in $A$ in Lemma 9.8 and which vanishes
where $|t| \geq
2 u/\a$ and $u \geq  2\nu$. Otherwise, the precise details of $\chi_*$
are not
relevant to the subsequent discussion.

\medskip
{\bf Step 9}\qua For almost all choices of $v \in (0,
\nu)$, the
arc $a(v)$ in the $(t, u)$ plane where $(t^2 + u^2)^{1/2} = v$ and where
$\a |t|
\leq  u$ will miss $\pi(\L)$. For such $v$, the $\pi$--inverse image of the arc
$a(v)$ is a smooth submanifold in $C$.

Suppose that $v$ is such that $a(v)$ misses $\pi(\L)$ and consider the
integral
of the imaginary part of $- d\th\,' /\th\,'$  over $a(v)$. Let
$L(v)$ denote this integral. To estimate this integral, first write
$\th\,'  = (t + i\dt u)^{q'} \dt \th\,''$. Then
\begin{equation}
L(v)\geq q'\dt\pi-\z\dt\a-\int_{a(v)}
{\mbox{im}}(d\th\,''/\th\,'').
\end{equation}
Here, the constant $\z$ is independent of $v$ and $\a$ as long as both are
small.

To estimate the right most term in (9.12), introduce polar coordinates $v =
(t^2 + u^2)^{1/2}$ and $\s = \tan^{-1}(u/t)$ in the $(t, u)$ plane.
Also, agree to
write $\th\,'' = e^{\D}$, where $\D$ is also annihilated by $\p$.
Then, write $\D = \D_1 + i\dt \D_2$ where $\D_1$ and $\D_2$ are real.
 (Note that
$\th\,''$ is not zero on $A$ when $\d$ is small -- this is
guarranteed by Lemma 9.7.) With the preceding understood, the right most
term in (9.12) is given by
\begin{equation}
L_1(v)\equiv -\int_{a(v)} \p_\s\D_2 \ d\s .
\end{equation}
Since $\p  \D = 0$, one has
\begin{equation}
v^{-1}L_1(v)=\int_{a(v)} \p_v\D_1 \ d\s .
\end{equation}
Given (9.14), take a pair of numbers, $0 < v' < v'' <\nu$
and integrate $v^{-1}L_1(v)$ over the
interval $[v' , v'' ]$. The result,
according to (9.14), is
\begin{equation}
\int_{[v',v'']} v^{-1}\dt L_1(v)\dt dv
=\int_{a(v'')} \D_1 \ d\s-\int_{a(v')}\D_1 d\s .
\end{equation}
The claim now is that if $v'$  is sufficiently small, then the
right-hand side of (9.15) is positive. To see that such is the case, use the
condition $|\eta|/z \leq  \e$ to see that $|\eta|/u < \e/(1 - \e^2/2)$ on
$C|_A$. Then, use the points in Lemma 7.7 to conclude that $|\th\,''|
 \leq  u^{1/4}$ on its domain of definition. Since $\th\,''$
does not vanish
when $\nu$ is small (due to Lemma 9.7), it follows that, given $v'' >
0$, there
exists $v' > 0$ such that the supremum of $\D_1$ where $v \in (0,v')$
 is strictly less than the infimum of $\D_1$ where $v =v''$.

Given that there exists $v' > 0$ where (9.15) is positive,
there exists $v > 0$ and as small as desired, where $L_1(v) > 0$. For such $v$,
the number $L(v)$ in (9.11) obeys the bound $L(v) > (q'  -1/100)$.

\medskip
{\bf Step 10}\qua Let $a' (v) = \pi^{ -1}(a(v))\subset C$.
Then, $L(v)$ is not equal to $-\int_{a' (v)} d\phi$, but almost so. In fact, it
follows from (7.22) that $-\int_{a' (v)} d\phi$ differs from $L(v)$ by no
more than a uniform multiple of $\e\dt |\ln(\a)|$. The argument for this is
as follows: First of all, the integral of $-d\phi$ over $a'
(v)\equiv\pi^{-1}(a(v))$ is equal to that of $\pi_*(-d\phi)$ over $a(v)$. Then,
the difference between the $L(v)$ and the integeral of $\pi_*(-d\phi)$ over
$a(v)$ is proportional to the integral over $a(v)$ of $-d(im(\tau))$. To
estimate this last integral, remark first that $\tau(\dt )$ can be written as
in (7.18) except that $\chi_*\dt m$ should replace $m$ inside
the integral. With
the preceding understood, let $\chi$ be a standard, non-increasing cut-off
 function
on $[0, \i)$ which takes value 1 on [0, 1] and 0 on $[2,\i)$. Then,
decompose the $\chi_*\dt m$ version of (7.22) into a sum of two
integrals; these being:
\begin{equation}
\begin{array}{l}
\displaystyle{\tau'_1=-i(2\pi)^{-1}\int_{u'\geq 0}\chi
 (|w-w'|/2u)\chi_*m(w')/((\bar w-\bar w')\!\dt\!(\bar w-w'))d\bar w'
  \wedge dw' }\vspace{1\jot}\\  
\displaystyle{\tau_2=-i(2\pi)^{-1}\int_{u'\geq 0}\!\!(1\!-\!\chi
 (|w-w'|/2u))\chi_*m(w')/((\bar w-\bar w')\!\dt\!(\bar w-w'))d\bar w'
  \wedge dw'} \end{array}
\end{equation}

Now, note that $|\tau_1| \leq \z \e$ for some constant $\z$
 which is independent
of $\e, \ \a$ and $\nu$. This follows using Lemma 7.8 and the following:
\begin{itemize}
\item $ \quad |w-\bar w'|\geq u \ {\mbox{when $w$ and $w'$
 lie in the upper half plane}}.$
\item $ \quad |w- w'|\geq 4u \ {\mbox{is required for a non-zero
contribution}}.$\itemnum{9.18}
\end{itemize}
The details here are identical to those in Section~7 which establish the $\z\e$
bound for (7.19).

As for $\tau_2$, the claim here is that
\begin{equation}
| d\tau_2|\leq \z\dt\e/u  ,
\end{equation}
where $\z$ is independent of $\e, \ \a$ and $\d$. Indeed, this
follows readily
by differentiating the second line of (9.16) and again invoking Lemma
7.8. In
this case, argue by breaking the integration region into annuli as described
in the derivation of the $\z \ln(\e)$ bound for (7.20) in Section~7.

Given the preceding, it follows from (9.18) and the estimate $|\tau_1| <
\z\dt \e$ that
\begin{equation}
|\int_{a(v)}d\tau| \leq \z\dt (\sup_{a(v)} |\tau_1|+\e \dt
\int_{a(v)} u^{-1}\leq z\dt (\e +\e\dt |\ln(\a)|).
\end{equation}

\medskip
{\bf Step 11}\qua This step digresses momentarily to
 justify a
simplifying assumption. The digression begins with the observation that the
subsequent arguments for Lemma 9.6 need only consider the case where, for
small $s'$, the corresponding $\g(s')$ (as defined in
(9.6)) has no circle component. Indeed, as already remarked in Step 3, if
there is a circle component, the integeral of $-d\phi$ over this component is a
positive, integral multiple of $2\pi$ since $-d\phi$ restricts to
$\g(s')$ as a volume form. Thus, it is sufficient to make the
following assumption: Let $s'  \in (0, s]$ and suppose that the
corresponding $\g(s')$ misses the critical points of $h$ on $C$ and
the singular points of $C$. Then $\g(s')$ is a finite union of
disjoint, embedded arcs with endpoints on $\p\,' C$.

\medskip
{\bf Step 12}\qua Choose $v > 0$ and small where
$-\int_{ a' (v)}
d\phi > (q'  - 1/100)$. A slight (and arbitrarily small) change of $v$
will insure that $a' (v)$ misses all critical points of $h$ on $C$,
and this condition will henceforth be assumed. Additional constraints on $v$
are as follows: First, by choosing $v$ small, the points in $C$ where
$(t^2 +
u^2)^{1/2} = v$ and $t = - u/\a$ will assuredly lie in $C \cap  D(s)$ on level
curve
$\{h = s_{1j}\}_{1\leq j\leq q'}$  with each $s_{1j} < s$.
Furthermore, after possibly a slight change of $v$, each of the corresponding
level sets $\{h^{-1}(s_{1j})\}$ misses all singular points of $C$, and each
$\{s_{1j}\}$
is a regular value of $h$. This implies that each level set $h =
s_{1j}$ is a
smoothly embedded curve, $\g(s_{1j})$, in the smooth part of $C$.
Likewise, no
generality is lost by assuming that the points in $C$ where both $(t^2 +
u^2)^{1/2} = v$ and $t = u/\a$ lies in $C \cap  D(s)$ on
 level curves $\{h =
s_{2j}\}$ with $s_{2j} < s$, and where each such level curve is a
 smoothly
embedded curve, $\g(s_{2j})$, in the smooth part of $C$.

There is some direction on each $\g(s_{1j})$ where $-\phi$ is
increasing, and
likewise on each $\g(s_{2j})$. Thus, $a(v)$ can be completed by
attaching the
relevant portions of $\g(s_{1j})$ and $\g(s_{2j})$ to a give a piecewise
smooth
arc in
$C \cap  D(s)$, where each endpoint lies either in $\p\,'C$ or $\p ''C$,
 and over which the integral of $-d\phi$ is
at least $(q'  - 1/100)\dt \pi$. Note that this arc has at
most $q'$  components, so there is at least one component, say
$\g'$, with $\int_{\g'} (-d\phi) > 15\dt \pi/16$.

\medskip
{\bf Step 13}\qua The construction of $\g'$  does not complete
the proof of Lemma 9.6 because $\g'$ does not necessarily lie in
a level set of $h$. To achieve the latter requirement, it is necessary to
modify $\g'$. This modification is carried out via a Morse
theoretic argument on $C$ which uses the function $h$ as the Morse
function and
uses the level sets of $\phi$ as the flow lines of a pseudo-gradient vector
field. However, before starting, remark that the argument that follows is
complicated by three technical issues:
\begin{itemize}
\item $C$ may have some singular points.
\item $h$ may have degenerate critical points.
\item The pseudo-gradient flow lines may not be in general
position.\hspace{\fill} (9.20)
\end{itemize} 
The first issue is dealt with by pulling back the story to the universal
model $C_0$ via the pseudo-holomorphic map $f\co C_0\to  S^1\x (B^3 -\{0\})$
as described in Definition 1.1. The second and third issues in
(9.20) will be dealt with by considering small perturbations of both $f^*h$
and $f^*\phi$.

In any event, to begin, choose $s_0 > 0$ but less than the infimum
of $h$ on
$\g'$  and so that $\g(s_0)$ avoids all critical points of $h$ and all
singularities of $C$. Let $f\co C_0\to  S^1\x \R^3$ denote the
smooth model for $C$ from Definition 1.1.  Let $D_0 =
f^{-1}(D(s) \cap \{h
\geq s_0\})$. Thus, the boundary of $D_0$ consists of a finite union of arcs,
where each arc lies either in a component of $f^{-1}(\p\,'C\cap\{
h \geq  s_0\})$ or a component of $f^{-1}(\p\,''C\cap\{
h \geq  s_0\})$ or else
in a level set of $f^*h$ with level either $s_0$ or $s$. (The
 latter coincide
with $f^{-1}(\g(s_0))$ and $f^{-1}(\g(s))$, respectively.)

On the complement of the critical set of $f^*h$, the tangents to the level
sets of $f^*\phi$ define a vector field which is pseudo-gradient like
for the
function $f^*h$. Note that $f^{-1}(\p\,' C)$ may not be a union
of level sets of $f^*\phi$, but Proposition 7.1 implies that no generality is
lost by requiring $f^*h$ to restrict to  $f^{-1}(\p\,' C)$
without critical points. Meanwhile, the requirements for the chosen s insure
that $f^*h$ also restricts to $f^{-1}(\p\,' C)$ without critical points. Thus,
one can choose a pseudo-gradient for $f^*h$ which is tangent to the level
sets of $f^*\phi$ except in some arbitrarily small, but apriori fixed
neighborhood of $f^{-1}(\p\,' C)\cup f^{-1}(\p\,''C)$.

The simplest case to consider assumes that $f^*h$ has only non-degenerate
critical points and that the tangents to the level sets of $f^*\phi$ give
 a
suitably generic pseudo-gradient for $f^*h$. This genericity assumption
will be enforced in the remainder of this step; the general case is
considered in Step 14.

Under the preceding genericity assumption, it follows by standard Morse
theoretic considerations that $D_0$ has a 1--skeleton which is
constructed as
follows: Start with the arcs of $(f^*h)^{-1}(s_0)$ and add a finite
set of arcs
in level sets of $f^*\phi$ whose endpoints lie on the arcs of
$(f^*h)^{-1}(s_0)$.

This last picture of the one skeleton of $D_0$ implies that there is a
continuous, embedded arc, $\l_1$, in $D_0$ with the following
properties:
\begin{itemize}
\item $\l_1$ is a union of components, each of which is either
an arc in $f^{-1}(\g(s_0))=(f^*h)^{-1}(s_0)$ or else an arc in a level
set of $f^*\phi$.
\item The boundary points of $\l_1$ lie in
$f^{-1}(\p\,' C\cap\p\,''C\cap\{h=s_0\})$ and each boundary point of
$\g'$ is connected by an arc in $f^{-1}(\p\,' C\cap\p\,''C\cap\{h\geq s_0\})$
to a boundary point of $\l_1$. \hspace{\fill} (9.21)
\end{itemize} 
Join the endpoints of $\g'$  and $\l_1$ in $f^{-1}(\p\,' C\cap \{h \geq  s_0\})$
and call the resulting curve $\l$. Orient $\l$ so that the part which
coincides
with $\g'$  is correctly oriented. Then the integral of $-d\phi$ over $\l$ has
the form $m\dt 2\pi$  for some integer $m$.

Note that if $m \leq  0$, then it follows (using Lemma 9.3) that the
integral of $-d\phi$ over a disjoint union of arcs in $\g(s_0)$ is at least
$3\pi /4$ (assuming that $s$ is small) and this implies Lemma 9.6. Thus,
 one can
assume that $m$ is positive.

Now, remark that using standard Morse theory, a basis for $H_1(D_0;{\Bbb Z})$
can be constructed so that each basis element is a continuously embedded
loop in $D_0$ which is a union of embedded arcs. Here, each arc lies in
$f^{-1}(\g(s_0)) = (f^*h)^{-1}(s_0)$ or in a level set of $f^*\phi$. Since the
integer $m$ is positive, there is at least one such generator on which
$-f^*d\phi$ integrates to a positive multiple of $2\pi$. This also
implies
Lemma 9.6.

\medskip
{\bf Step 14}\qua Now consider the case where $f^*h$ has some
degenerate critical points and/or where the level sets of $f^*\phi$
are not in
general position. In either case, one can perturb $f^*h$ and the
pseudo-gradient
vector field in small disks about the critical points of $f^*h$ so that the
resulting function, $h'$, has non-degenerate critical points;
and so that the resulting pseudo-gradient flow lines are in general position.

In this regard, note that these perturbations can be assumed to be supported
in a union, $V \subset D_0$, of disks of arbitrarily small (but positive)
 radius
about the critical points of $h$. This means, in particular, that the new
function,
$h'$, can be assumed to agree with h except in $V$. And, the new
gradient flow lines are level sets of $f^*\phi$ except in $V$ and except
(possibly) near $f^{-1}(\p\,' C)$. Furthermore, since $f^*h$ is
a harmonic function on $C_0$ (with respect to an appropriate metric), the
function $h'$  can be assumed to have only index 1 critical
points also. Thus, the number of such critical points can be taken to be
independent of the diameter of the disks which comprise $V$.

(The point here is that there is some complex coordinate for $C_0$
 near each
critical point of $f^*h$ for which the latter is the real part of a
holomorphic function on a small disk in $\Bbb C$. One can then consider
standard perturbations of such functions.)

With the preceding understood, standard Morse theory arguments can be used
to construct the 1--skeleton of $D_0$ as follows: Start with the arcs of
$f^{-1}(\g(s_0)) = (f^*h)^{-1}(s_0)$ and add additional arcs with end points on
those of $(f^*h)^{-1}(s_0)$. The latter arcs lie in level sets of $f^*\phi$
except where they intersect $V$. Furthermore, the number of such arcs
can be
assumed independent of the diameter of $V$, and the length in $V$ of each such
arc can be assumed bounded by a uniform multiple of the diameter of $V$.
(Remember that $f^*h$ is the real part of a holomorphic function with
respect to some complex coordinate on a neighborhood of each critical point.)

Given this 1--skeleton, repeat the arguments of the previous step.  The only
difference here comes from the fact that $-f^*d\phi$ may have non-zero
integral over an arc segment in $V$. However, the contribution from the
segments in $V$ can be made negligibly small in absolute value because of
the following facts:
\begin{itemize}
\item Each arc segment in $V$ has length bounded by a multiple
of the radius of a disk in $V$.
\item $f^*dh=0$ at the center of each disk component of $V$,
so $f^*d\phi=0$ there as well. Thus, $|f^*d\phi|$ is small in $V$ when the
disks of $V$ have small radii.
\end{itemize} \hspace{\fill} (9.22)

\med {\bf Acknowledgement}\qua The author is supported in part by the
National Science Foundation.

\newpage

\appendix\small\parskip 4pt plus 2pt minus 2pt

\makeatletter
\let\@@itemize@\itemize
\def\itemize{\@@itemize@\parskip 0pt\relax}
\let\@@enumerate@\enumerate
\def\enumerate{\@@enumerate@\parskip 0pt\relax}
\def\@listi{\leftmargin25.5pt\parsep 0pt\topsep 0pt 
 \itemsep2pt plus2pt minus1pt}
\let\@listI\@listi
\@listi
\makeatother

\section*{Appendix\\The proof of Assertion 3(a) of Proposition 7.1}

\setcounter{section}{1}
\renewcommand{\thesection}{\Alph{section}}
\renewcommand{\theequation}{\thesection.\arabic{equation}}
\setcounter{equation}{0}

The proof is broken into three parts, where each part is
further decomposed into separate steps.  Part I of this
Appendix discusses in some detail the behavior of
$C$ in a small neighborhood of one of its $u > 0$ singular
points. Part II uses the analysis from the first part
of the Appendix to prove Assertion 3(a) of Proposition 7.1
under a restrictive assumption. The final part  of this
Appendix proves that the restrictive assumption in part II
is always valid.

Before starting, a digression is in order to set the
stage. In particular, attention is focused on a
neighborhood in $S^1\x B^3$ of a point
$t_0\in Z=S^1\x\{0\}$. By assumption $C$ is a finite
energy, pseudo-holomorphic curve and $t_0$ is a regular
point for $C$ where $q_{\pm}=0$. Thus, $|z|/\rho << 1$
on $C$ near $t_0$.

\subsection*{Part I\qua The local model where $u > 0$}

This part of the Appendix constitutes a digression of
sorts to present a local model for $C$ in the neighborhood
of $u > 0$ singular
point. This part of the proof is broken into six steps.

\medskip{\bf Step 1}\qua To begin, consider a
regular point $t_0\in Z$. As in subsection 8.a, introduce the
$(t,u)$ coordinate system on a neighborhood of $t_0$ as
above so that $t_0=(0,0)$. For $\d >0$ and small, let
$A=\{(t,u): t\in [-\d,\d]$ and $u\in [0,\d]\}$. According
to the discussion in subsection 8.a, when $\d$ is small, $C$
(in a neighborhood of $t_0$) can be described by a
multivalued graph $\Phi=\{(\phi_j,\nu_j)\}_{1\leq j\leq
p}$ over $A$, where each $(\phi_j,\nu_j)$ is a locally
defined function on the complement of a set $\pi(\L)$, in
the interior of $A$. Here $\L\subset C$ is a certain set
with at most countably many elements and  no accumulation
points, while $\pi$ denotes the projection induced map from
the range of $\Phi$ in $C$ to $A$. (The set $\L$ consists
of the singular points of $C$ and the critical points of
the function $t$.)  To conserve notation, henceforth use
$C$ to denote only the image of $\Phi$.

According to Lemma 8.3, the map $\Phi$, as a map into the
$p^{\mbox\footnotesize{th}}$ symmetric product of ${\R}^2$,
is Holder continuous with some exponent $\a >0$. In
addition, given $\e > 0$, one can choose $\d$ so that each
$(\phi_j,\nu_j)$ obeys (8.2) and so that (8.3) holds on
the complement of $\L$ in the interior of $A$. One can
also assume (8.4).

Note that $\{(\phi_j,\nu_j)\}$ are obtained by restricting
the functions $(\phi,\nu)$ on $S^1\x B^3$ to $C$. Further
progress requires a digression to introduce a replacement,
$b(\nu)$, for the function $\nu$ on the subset of $S^1\x
B^3$ where $z/\rho < 1/32$. Here is $b$:
\begin{equation}
b(\nu) = \int^\nu_0 (1+6s^2)^{-\frac 12}(1+2s^2)^{-1}ds
=\nu+{\cal O}(\nu^3) .
\end{equation}
Note that $b$ is an analytic function of $\nu$ where the
latter is less than 1/32.

The advantage of replacing $\nu$ with $b$ arises when
(8.3) is written in terms of the complex variables
$w=t+i\dt u$ and $\eta=\phi+i\dt b$. The resulting
equation has the schematic form
\begin{equation}
\bar\p\eta +\a_1\p\eta -\textstyle\frac{3}{2} u^{-1}
(1+\a_2)b=0  ,
\end{equation}
where $b={\mbox{im}}(\eta)$ while $\a_1$ and $\a_2$ are
real valued, real analytic functions of $b^2$ where
$|b|\leq 1$. Furthermore, both are ${\cal O}(b^2)$ near
$b=0$. (The analogous equation for $\eta'=\phi+i\dt \nu$
has a term with derivatives on $\bar\eta'$. The presence
of such a term complicates the already complicated
discussion to come.)

To derive (A.2), consider $\nu$ as a function of $b$
instead of vice-versa. Then the $b$--derivative of $\nu$
obeys $\nu'=(1+\kappa_1)^{1/2}(1+\kappa_2)^{1/2}$, where
$\kappa_1$ and $\kappa_2$ appear in (8.3). With this
understood, one finds (A.2) with $\a_1=(u-g)/(u+g)$
and $\a_2=u^3\rho^{-2}(u+g)^{-1}-1$.

Note that in terms of $w$ and $\eta$, the sheets of $C$
over $A-\pi(\L)$ are described by a set $\{\eta_j\}_{1\leq
j\leq p}$ of locally defined solutions of (A.2).
Alternately, $\{\eta_j\}$ can be thought of as a Holder
continuous map from $A$ into Sym${}^p(\Bbb C)$.

\medskip
{\bf Step 2}\qua
Introduce now the smooth model $C_0$ as described in 
Definition 1.1. Restricting to a subset of $C_0$ gives a
smooth, complex curve $C_A$ with a proper,
pseudo-holmorphic map $f$ from $C_A$ into $S^1\x
(B^3-\{0\})$ whose image is $\pi^{-1}(A-\{u=0\})$. By
assumption, $f$ is 1--1 on the complement of a countable
subset. Note that the composition $\pi\dt f$ is generally
not a holomorphic map $f$ from $C_A$ into $A$. However,
this map is close to being holomorphic where $u$ is small.
Indeed, let $v$ be a local, complex coordinate for $C_A$.
Then, $f_*\p_v= t_v\p_t+u_v\p_u+h_v\p_h+\phi_v\p_\phi$
and
\begin{equation}
J\dt f_*\p_v=-t_v gu^{-1}\p_u +u_v ug^{-1}\p_t+h_v(g\rho^2)^{-1}
\p_{\phi}-\phi_v g\rho^2\p_h .
\end{equation}
The condition that $f$ is $J$--pseudo-holmorphic requires
$J\dt f_*\p_v= -i\dt\p_v$; and the latter requires (in
part) that $t_v-i\dt u_v(u/g)=0$. The complex conjugate of
this last equation can be rewritten as
\begin{equation}
\bar\p_v w-\a_1\bar\p_v\bar w=0  ,
\end{equation}
where $\a_1=(u-g)/(u+g)$. Since $\a_1$ is ${\cal
O}(b^2)={\cal O}(u^2)$, this last equation implies that
$\pi\dt f$ is nearly holomorphic where $u$ is small.

By the way, this last equation with (A.2) implies that that
$\eta$ as a function of $v$ obeys
\begin{equation}
\bar\p\eta -\textstyle\frac{3}{2}(\bar\p_v\bar w)u^{-1}
(1+\a_2)b=0.
\end{equation}
Now suppose that $v$ is a complex parameter for some disk
$D\subset C_A$ and that $f$ maps the origin (where $v=0$)
to a point with coordinates $w=w_0$ and $\eta=\eta_0$.
Of course (A.4) implies that $w$ is a real analytic
function of $v$ (when $|v|$ is small). In particular, this
last fact with (A.4) implies that the parameter $v$ can be
chosen so that
\begin{equation}
w(v,\bar v)=w_0+v^k +\a_1(b_0)\bar v^k+{\cal R}(v).
\end{equation}
Here $k\in\{1,\dots p\}$ and $\cal R$ is a real analytic
function of $v$ which obeys $|{\cal R}|\leq \z \dt |v|^{k+1}$.
Note that the upper bound of $p$ here stems from the
following two facts: First, $f$ is almost everywhere 1--1
onto $C$. Second, the projection $\pi$ from $C$ to $A$ has
at most $p$ sheets.

In the case where $w_0\not\in\pi(\L)$, then
$f^*\pi^*dw\neq 0$ at $v=0$, and this implies that $k=1$
is the only possibility. However, when $w_0\in\pi(\L)$,
then $k > 1$ is allowed. Even so, when $k=1$, then (A.6)
implies that $f$ embeds a neighborhood of $v=0$ into
$S^1\x (B^3-\{0\})$. In this case, there is a disk $U$
surrounding $w_0$ with the property that $C|_U$ has at
least two irreducible components, and at least one of
these is an embedded disk. (An irreducible component of
$C$ over an open set $U\subset A$ is the closure in $C|_U$
of  a connected component of $C|_U-\L|_U$.)

\medskip
{\bf Step 3}\qua Once again, let $\pi\co C\to A$
denote the projection induced map and use $\pi$ to define
the pull-back space $\pi^*C$. The latter consists of the
space of triples $(w,\eta,\eta')$ where $w\in A$ and where
both $\eta$ and $\eta'$ are points in $\pi^{-1}(w)$. This
space $\pi^*C$ is a manifold at points $(w,\eta,\eta')$
where neither pair $(w,\eta)$  nor $(w,\eta')$ is a
singular point of $C$. In particular $\pi^*C$ is smooth
near $(w,\eta,\eta')$ when $w\not\in \pi(\L)$. The space
$\pi^*C$ has two canonical projections to $C$ which are
denoted by $\pi_{\pm}$ below. Here $\pi_-$ sends
$(w,\eta,\eta')$ to $(w,\eta)$ while the image of this
point under $\pi_+$ is $(w,\eta')$. Moreover, note that $C$
embeds canonically in $\pi^*C$ as $\{(w,\eta,\eta)\}$,
that is, where $\pi_+=\pi_-$.

The remainder of this step describes the 
irreducible components of a neighborhood of the point
$(w_0,\eta_0,\eta'_0)$ in $\pi^*C$. (Remember: an 
irreducible component of a neighborhood ${\cal
U}\subset\pi^*C$ of $(w_0,\eta_0,\eta'_0)$ is defined to be
the closure in ${\cal U}$ of a component of 
${\cal U}(w_0,\eta_0,\eta'_0)$.)

The simplest case to consider is when $w_0\not\in\pi(\L)$.
Here there is a disk $U\subset A$ centered about $w_0$
which has empty intersection with $\pi(\L)$. Then $C|_U$
is the disjoint union of $p$ copies of $U$, and then
$\pi^*C$ over $C|_U$ is the disjoint union of $p^2$ copies
of $U$, where $p^2-p$ copies lie in $\pi^*C-C$. In
particular, $\pi^*C|_U$ near $(w_0,\eta_0,\eta'_0)$ is
diffeomorphic to $U$ via the obvious projection.

The structure of $\pi^*C$ near $(w_0,\eta_0,\eta'_0)$ is
more complicated than $w_0\in\pi(\L)$. In this case it
proves useful to consider separately the cases where
$\eta_0\neq\eta'_0$ and $\eta_0=\eta'_0$. In the first
case, $(w_0,\eta_0)$ and  $(w_0,\eta'_0)$ lie in distinct
components of $C|_U$, and then a neighborhood of 
$(w_0,\eta_0,\eta'_0)$ over $U$ is
homeomorphic to the fiber product of these two components.
In particular, if the first component intersects
$\pi^{-1}(U-w_0)$ as a $k$--sheeted cover of $U-w_0$, and
if the second component likewise defines a $k'$--sheeted
cover, then there is a  neighborhood of 
$(w_0,\eta_0,\eta'_0)$ in $\pi^*C$ which intersects
$\pi^{-1}(U-w_0)$ as a $k\dt k'$ sheeted cover of $U-w_0$.
Furthermore, the number of irreducible components of
$\pi^*C$ near $(w_0,\eta_0,\eta'_0)$ is equal to the
greatest common divisor of $k$ and $k'$.

Note that the preceding assertions follow from (A.6)
in as much as the latter implies that when $|w-w_0|$ is
small but not zero, there are precisely $k$ distinct
points in the $v$--disk which map to $w$. Furthermore,
each such point is nearly a $k^{\mbox\footnotesize{th}}$
root of $w-w_0$. In fact, according to (A.6), each $v$
which maps to the given $w$ differs from a  $k^{\mbox\footnotesize{th}}$
root of $w-w_0$ (and vice versa) by ${\cal O}(\e)$ where
$\e=|w-w_0|^{1/k}(|u|^2|+|w-w_0|^{1/k})$.

The second case to consider has $\eta_0=\eta'_0$. In this
case $\pi^*C-C$ intersects a  neighborhood of 
$(w_0,\eta_0,\eta_0)$ in $\pi^*C|_U$ in components which
can be described as follows: The components are partially
labeled by ordered pairs of irreducible components of
$C|_U$. This label for a point $(w,\eta,\eta')$
corresponds to the ordered pair of irreducible components
which  respectively contain $(w,\eta)$ and $(w,\eta')$.
With the preceding understood, let $B$ be an 
irreducible component of
$C|_U$ whose intersection with $\pi^{-1}(U-w_0)$ defines a
$k$--sheeted cover of $U-w_0$. Let $B'$ be a distinct
irreducible component of
$C|_U$ and define $k'$ from $B'$ in an analogous manner.
Then the pair $(B,B')$ label a certain number of
irreducible components of some neighborhood in $\pi^*C$ of
the given point $(w_0,\eta_0,\eta_0)$. And this number is
(as in the first case) equal to the greatest common
divisor of $k$ and $k'$.

On the other hand, if $B=B'$, then the pair $(B,B)$ labels
$k-1$ distinct components of the intersection of $\pi^*C$,
and thus $k$ irreducible components of some neighborhood 
$(w_0,\eta_0,\eta_0)$ in $\pi^*C$. (Precisely one of these
components lies in $C$.)

Note that the preceding assertions also follow from (A.6).

\medskip
{\bf Step 4}\qua The preceding description
of $\pi^*C$ can be made somewhat more explicit with the
help of a certain parametrization of $\pi^*C$ near points
in $\L$. Consider first the parametrization for $\pi^*C$
near a point $(w_0,\eta_0,\eta'_0)$ where $w_0\in\pi(\L)$.
For this purpose, fix respective, irreducible components
$B$ and $B'$ for $C$ near $(w_0,\eta_0)$ and for 
 $(w_0,\eta'_0)$ as in Step 3. Let $a$ denote the
greatest common divisor of $k$ and $k'$, and write
$k=\l\dt a$ and $k'=\l'\dt a$, where $\l$ and $\l'$ are
relatively prime.

As remarked in Step 3, a component of $\pi^*C$ near
 $(w_0,\eta_0,\eta'_0)$  is labeled by an integer
$s\in\{0,\dots ,a-1\}$ in addition to the label by the
pair $(B,B')$. A coordinate for this component is given
by a complex parameter $\tau$ for a small radius disk in
$\Bbb C$ centered at 0. The following lemma describes the
parametrization of points in $\pi^*C$ by $\tau$.

\medskip
\noindent{\bf Lemma A.1}\qua{\sl{There is a constant
$\z\geq 1$ with the following significance: Suppose that
$w_0\in\pi(\L)$ is a point with $u_0\leq\z^{-1}$.
Let $(w_0,\eta_0,\eta'_0)\in\pi^*C$ be given, and let
$(B,B')$ and $s$ parametrize an irreducible component of
a neighborhood of $(w_0,\eta_0,\eta'_0)$ in $\pi^*C$
as described above. Write a neighborhood of $(w_0,\eta_0)$
in $B$ as the image via the pseudo-holomorphic map $f$ of
a disk $D\subset C_0$ with coordinate $v$ as in} (A.6),
{and likewise parametrize 
a neighborhood of $(w_0,\eta'_0)$ in $B'$ by the $f$ image
of a disk $D'$ with complex coordinate $v'$. Then a 
neighborhood of the point $(w_0,\eta_0,\eta'_0)$ in the
given component of $\pi^*C$ is parametrized by a complex
coordinate $\tau$ on a disk about the origin in $\Bbb C$
by writing}
\begin{eqnarray}
v &=& \tau^{\lambda'} \nonumber \\
v'&=& \exp(2\pi i s/k')\tau^{\lambda}(1+m(\tau)),
\end{eqnarray}
{where $m$ is a certain function of $\tau$ which is
smooth where $\tau\neq 0$, Lipschitz near $\tau=0$ and
satisfies the bound} $|m|\leq\z (|\tau|^\l
+|\tau|^{\l'})$.}

\medskip
{\bf Proof of Lemma (A.1)}\qua The function $m(\tau)$ is
chosen so that $w(v(\tau))=w(v'(\tau))$. In particular, if
(A.6) is used for both $v$ and $v'$ then this last
condition can be written as follows
\begin{eqnarray}
&& v^k[(1+m)^{k'}-1]+\a_1(b_0)\bar v^k[(1+\bar m)^{k'}-1]= \nonumber \\
&& {\cal R}(\tau^{\lambda'})-{\cal R}(\exp(2\pi i s/k')\tau^{\lambda} (1+m(\tau))).
\end{eqnarray}
This algebraic equation for $m$ can be written in fixed
point form as
\begin{equation}
m=v^{-k}k^{'-1}({\cal R}(\tau^{\lambda'})-{\cal R}(e^{2\pi i
s/k'}\tau^{\lambda}))-\a_1(b_0)\bar v^k/v^k\bar m+{\cal K}  ,
\end{equation}
where ${\cal K}$ is an analytic function which satisfies
$|{\cal K}|\leq\z(u|m|^2+|\tau|^\l|m|)$.

The appropriate solution to (A.9) will be found with the
help of the contraction mapping theorem. For this purpose,
it proves useful to first replace $\tau$ in (A.11) by an
indeterminant $x_1$ and to likewise replace $v/|v|$ by an 
indeterminant $x_2$.  Then (A.11) has the schematic form
\begin{equation}
m={\cal F}(m,x_1,x_2)  ,
\end{equation}
where ${\cal F}$ is a real analytic function of its entries
when $u_0$ is small, when $m$ and $x_1$ take values in a
small radius disk in $\Bbb C$, and when $x_2$ take values
in a disk of radius 2. Furthermore,
\begin{eqnarray}
|{\cal F}| & \leq & \zeta (|x_1|^{\min(\l,\l')}+ u_0|m|) \nonumber \\
|\p_m{\cal F}| & \leq & \zeta (|x_1|^{\l}+ u_0),
\end{eqnarray}
where $u_0$, $|x_1|$ and $|m|$ are all smaller than some fixed constant
$\z^{-1}$.

Given (A.11), the contraction mapping theorem finds $\z\geq 1$ such that
when $u_0$ and $|x_1|$ are both less than $\z^{-1}$, and
$|x_2|\leq 2$, then (A.10) has a unique, small norm solution
$m=m(x_1,x_2)$. In addition, the contraction mapping theorem guarantees
that this solution obeys $|m|\leq\z\dt |x_1|^{\min(\l,\l')}$ and varies
with $x_1$ and $x_2$ as a real analytic function. With this last point
understood, then $m$ in Lemma A.1 is obtained from this fixed point
$m(x_1,x_2)$ of (A.10) by setting $x_1=\tau$ and $x_2=v/|v|$.

\medskip
{\bf Step 5}\qua
The set $\{\eta_j\}$ is the push-forward via the map $\pi$ of  a
well-defined function on $C$, namely the complex function
$\eta=\phi +i\dt b$. Likewise, the set of differences
$\{\eta_j-\eta_k\}$ is the push-forward from $\pi^*C$ of the complex
function $\Delta\equiv\pi_-^*\eta-\pi^*_+\eta$ via the map $\pi$ from
$\pi^*C$ to $A$. In particular, this function $\Delta$ maps
$\pi^*C-C$ to ${\Bbb C}-\{0\}$. Of interest in this step is the behavior of
$\Delta$ on the neighborhoods in $\pi^*C$ which are described in Lemma
A.1. For this purpose it proves convenient to use (A.5) to obtain an
equation for $\pi^*C$ for $\Delta$. Consider:

\medskip\noindent{\bf Lemma A.2}\qua{\sl{There exists $\z\geq 1$ with the
following significance: Suppose that $w_0\in\pi(\L)$ is a point with
$u_0 <\z^{-1}$. Let $(w_0,\eta_0,\eta'_0)\in\pi^*C$ be given, and let
$(B,B')$ and $s$ parametrize an irreducible component of
a neighborhood of $(w_0,\eta_0,\eta'_0)$ in $\pi^*C$
as described above. As in Lemma} (A.1), {let $\tau$
parametrize 
a neighborhood of $(w_0,\eta_0,\eta'_0)$ in the given
irreducible component. Then  the function $\Delta$ pulls back as a
function of $\tau$ to obey an equation of the form}
\begin{equation}
\bar\p_\tau\Delta +{\cal J}\dt{\mbox{im}}(\Delta)=0,
\end{equation}
{where ${\cal J}$ is a function of $\tau$ which is smooth except at
$\tau=0$, but is, in any event, bounded near $\tau=0$. In fact,}
\begin{equation}
|{\cal J}| \leq \zeta'(u_0)|\tau|^{\lambda'-1}  ,
\end{equation}
{where $\z'$ depends, as indicated, on} $u_0$.}

\medskip{\bf Proof of Lemma A.2}\qua To obtain (A.12), first subtract the
$\eta=\eta_k$ version of (A.2) from the $\eta=\eta_j$  version, pull the
result up to $\pi^*C$, and then pull back to the $\tau$ disk to obtain an
equation for $\Delta$. The latter has the following schematic form:
\begin{eqnarray}
&& \bar\tau_{\bar w}\bar\p_\tau\Delta +\tau_{\bar w}\p_\tau\Delta
  +\a_{1-}(\bar\tau_w\bar\p_\tau\Delta+\tau_w\p_\tau\Delta) 
 +(\a_{1+}-\a_{1-})
(\bar\tau_w \bar\p_\tau\eta_+ +\tau_w\p_\tau\eta_+)\nonumber\\
&& \qquad -\textstyle\frac{3}{2} u^{-1}(b_- -b_+ +(\a_2 b_-) -(\a_2 b_+))=0 .
\end{eqnarray}
Here $b_{\pm}={\mbox{im}}(\eta_{\pm})$, $\a_{1\pm}\equiv\a_1(b_{\pm})$
and $\a_{2\pm}\equiv\a_2(b_{\pm})$. Meanwhile,
$\tau_w\equiv\p_w\tau$, \ $\bar\tau_w=\p_w\bar\tau$, etc. To simplify
this last equation,  remark first that (A.4) and the first point of (A.7)
imply that $\tau_{\bar w}=-\a_{1-}\tau_w$. Thus the terms with
$\p_\tau\Delta$ in (A.14) cancel. Then using Taylor's theorem with
remainder, (A.14) can be rewritten as
\begin{equation}
\bar\tau_{\bar w}(1-|\a_{1-}|^2)\bar\p_\tau\Delta +{\cal J}'
{\mbox{im}}(\Delta)=0
\end{equation}
where $|{\cal J}'|\leq\z(u_0(1+\sup_j |\p_w\eta_j|)+u_0^{-1})$.
In fact, because $\eta$ pulls back to $C_A$ as a smooth function of $v$,
while $w$ obeys (A.6), it follows that 
$\sup_j |\p_w\eta_j|\leq\z|w-w_0|^{-1+1/k}$ near $w_0$. This last bound
implies that $|{\cal J}'|\leq\z|\tau|^{-\l\l'a+\l'}$ since $u_0$ is bounded
away from zero near $w_0$. Meanwhile, it follows from (A.4) and the first
point of (A.7) that $\bar\tau_{\bar w}(1-|\a_{1-}|^2)=1/w_\tau$, so (A.15)
implies (A.12) with ${\cal J}=w_\tau{\cal J}'$, and so (A.6) and (A.7)
give $|{\cal J}|\leq\z|\tau|^{\l'-1}$ as claimed.

\medskip
{\bf Step 6}\qua Equations (A.12) and (A.13) can now be used
to analyze the behavior of $\Delta$ as a function of the parameter $\tau$
from Lemmas A.1 and A.2.  The results are summarized by:

\noindent 
{\bf Lemma A.3}\qua{\sl{Continue with the assumptions of Lemma A.2.
The constant $\z$ in said lemma can be chosen so that the following
additional conclusion holds: Suppose that $\tau$ does not parametrize a
subset of $C$ in $\pi^*C$. Where the parameter $\tau$ obeys
$|\tau|\leq\z^{-1}$, the function $\Delta$ can be written as}
\begin{equation}
\Delta ={\mbox{re}}(\eta_0-\eta'_0)+e^\kappa\Delta_0  ,
\end{equation}
{where $\Delta_0$ is a non-zero, holomorphic function of $\tau$ with}
re$(\Delta_0(0))=0$. {Also, the function $\kappa$ is continuous, and
smooth except possibly at $\tau =0$; in any event,}
$|\nabla\kappa|\leq\z|\ln(\tau)|$. {In addition, if}
${\mbox{im}}(\eta_0-\eta'_0)=0$, {then}
\begin{equation}
\Delta_0=c_n\tau^n +{\cal O}(|\tau|^{n+1}) ,
\end{equation}
{where $n\geq 1$ is an integer, and where $c_n$ is a nonzero
constant. Furthermore, in this last case, $|\nabla^j\kappa|$ is bounded
for all $j\in\{0,\dots ,\l'-1\}$ and vanishes at the origin.}}

\medskip {\bf Proof of Lemma A.3}\qua Let
$\underline{\Delta}\equiv\Delta-{\mbox{Re}}(\eta_0-\eta'_0)$. Note that
im$\underline{\Delta}={\mbox{im}}(\Delta)$. With this understood,
choose $\kappa$ to obey the equation
$\p_\tau\kappa+ {\cal
J} \dt{\mbox{im}}\underline{\Delta}/\underline{\Delta}=0$ near the origin
of the $\tau$ plane. This  and (A.12) imply that $\Delta_0$ is
holomorphic and Re$(\Delta_0(0))=0$. In this regard, note that the
condition $\Delta_0\equiv 0$ implies that $\pi_-\eta-\pi_+\eta$ is
constant near $\tau=0$. Since $\Delta$ is real analytic on $\pi^*C-C$
away from $\pi^{-1}(\L)$, then $\pi_-\eta-\pi_+\eta$ is
constant on the whole of a component of $\pi^*C$, and this is forbidden
by the given assumptions. As $\Delta_0$ is not identically zero, then
$\Delta_0$ cannot be constant when im$(\eta_0-\eta'_0)$ vanishes.
This implies the assertion about (A.17).

Now consider the assertions about the derivatives of $\kappa$.
To begin remark that the equation in the preceding paragraph for $\kappa$
has infinitely many solutions since any holomorphic function can be added
to any one solution. Nonetheless, there is at least one solution with the
required behavior, as the following argument shows: Extend the
differential equation for $\kappa$ to the whole of the $\tau$ plane by
replacing ${\cal J}$ with $\chi_r {\cal J}$, where $\chi_r(\tau)\equiv
\chi(|\tau|/r)$ with $\chi$ the standard bump function and $r > 0$ a
small, positive number.  The solution for $\kappa$ is 
$\kappa(\tau)=\pi^{-1}\int (\tau'-\tau)^{-1}\chi_r{\cal J}
{\mbox{im}}\underline{\Delta}/\underline{\Delta} d^2\tau'+\kappa_0$,
where $\kappa_0$ is a polynomial function of the complex variable $\tau$.
Here $\kappa_0$ is chosen so that $\kappa$ and its derivates to order
$\l'-1$ vanish at the origin. In this regard (A.16) is invoked to estimate
the derivatives to order $\l'-1$ of the function
$\pi^{-1}\int (\tau'-\tau)^{-1}\chi_r{\cal J}
{\mbox{im}}\underline{\Delta}/\underline{\Delta} d^2\tau'$ at $\tau=0$.

Note that the procedure just outlined fails to produce the asserted bound
for $|\nabla\kappa|$ only in the case when $\l'=1$. The argument for this
case rests on the following claim: The function ${\cal J}
{\mbox{im}}\Delta/\Delta$ can be written as $\nu(\tau,\tau/|\tau|)$,
where $\nu(\dt \ , \ \dt)$ is a smooth function on ${\Bbb C}^2$.
Granted this claim, differentiate the integral expression for $\kappa$ 
given in the previous paragraph to find that as long as $\tau\neq 0$, one
has $|\nabla\kappa|\leq\z\int_{|\tau'-\tau|\leq 2r}
|\tau'-\tau|^{-1} |\tau'|^{-1} d^2\tau'$.
The integral here is bounded by $\z|\ln(\tau)|$ which is the required
bound.

Argue as follows for the claim in the preceding paragraph:
First, say that a function of $\tau$ is {two-variable smooth} if it
can be written as $\nu(\tau,\tau/|\tau|)$,
where $\nu$ is a smooth function on ${\Bbb C}^2$. A basic fact to be used
below is the following: if a function $\sigma$ is 2--variable smooth,  and 
$|\tau|^{-1}\sigma$ is bounded near $\tau =0$, then $|\tau|^{-1}\sigma$
is also two variable smooth. (Use Taylor's theorem to write 
$\nu=\tau\nu_++\bar\tau\nu_-$, where $\nu_{\pm}$ are two-variable smooth.)

With the preceding understood, note that $\Delta$ will fail as a 
smooth function of $\tau$ only when $v'$ in (A.7) does not have smooth
dependence on $\tau$; and the latter can happen due to the presence of
the function $m(\dt)$. However, it follows from (A.10) that $m$ is
two-variable smooth. However, because $\p\Delta$ is bounded, so it must
also be two-variable smooth. This implies the claim directly in the case
where $\Delta_0(0)\neq 0$. When $\Delta_0(0)=0$, then according to (A.17),
one has $\Delta_0=\tau^n\iota(\tau)$, where $\iota$ is complex analytic
and non-vanishing at 0. Thus $\underline{\Delta}$ has the form
$\tau^n\sigma$ where $\sigma$ is 2--variable smooth and non-vanishing at
$\tau =0$. In this case,
$\p\underline{\Delta}=|\tau|^{-1}\tau^n\sigma_1$,
where $\sigma_1$ is 2--variable smooth. This last observation implies that
${\cal J}{\mbox{im}}(\underline{\Delta})/\underline{\Delta}$
(which equals $\underline{\Delta}^{-1}\bar\p\underline{\Delta}$)
has the form $|\tau|^{-1}\sigma_2$ where $\sigma_2$ is 2--variable smooth.
However, ${\cal J}{\mbox{im}}(\underline{\Delta})/\underline{\Delta}$
is bounded, so it must also be 2--variable smooth.

\subsection*{Part II\qua A finiteness proposition for $\L$
under a density assumption}

The main result here, Proposition A.7, proves Assertion 3(a) of
Proposition 7.1 under an auxiliary assumption whose justification is
deferred to Part III of this Appendix. This portion of the proof of
Assertion 3(a) of Proposition 7.1 is broken into ten steps.

\medskip
{\bf Step 1}\qua To begin, consider a regular point $t_0\in
Z$ and reintroduce  $\d >0$, the set $A$ and the coordinates 
$w=t+i\dt u$ and $\eta=\phi+i\dt b$ as in Step 1 of Part I.
According to Lemma 8.3, the assignment of $w$ to
$\{\eta_j(w)\}=\pi^{-1}(w)\subset C$ defines a continuous map from $A$ to
Sym${}^p\Bbb C$, and the continuity of this map facilitates a
stratification of a neighborhood of the point $t_0$ as
$\bigcup_{1\leq n\leq p}A_n$, where a point $w\in A_n$ if $n$ or more
elements from $\{\eta_j(w)\}$ coincide. Note that each $A_n$ is closed;
and that $A_n-A_{n+1}$ is relatively open in $A_n$. Thus each $w\in
A-A_{n+1}$ has an open neighborhood which is also in $A-A_{n+1}$.
(This all follows from Lemma 8.3.) Furthermore, without loss of
generality, one can assume that $A_2$ consists of the union of $\pi(\L)$
with some subset of $Z$.

With $A_n$ understood, let $Z_n\equiv A_n\cap Z$.

This stratification of $Z$ will now be used to formulate an induction
proof that $\L$ is finite. Here is how: Suppose that $n\geq 2$ and that
for all $m < n$, it has been established that all limit points of $A_m$
are in $Z_n$. The induction step will prove that all limit points of
$\bigcup_{m\leq n}A_m$ lie in  $Z_{n+1}$. For this purpose it is
sufficient to assume that the given point $t_0$ lies in $Z_n$ and then
prove that $t_0$ is an isolated point of $A_n$. With this understood,
remark that $\d$ is chosen small in Step 1 of Part I, then
$A\cap A_{n+1}=\emptyset$ and this implies that it is sufficient to
consider only the case $n=p$ in the induction step.

Here is one last remark for this step: By taking $\d$ smaller if
necessary, one can arrange, with no generality lost, that each $\phi_j$
takes values only in (--1/4,1/4) on $A$.

\medskip
{\bf Step 2}\qua This step introduces the following basic
assumption:
\begin{equation}
{\mbox{The set $Z-Z_p$ is dense in $Z$ near $t_0$.}}
\end{equation}
The justification of this assumption is deferred to Part III of this
Appendix.

The assumption in (A.18) is used in the subsequent arguments solely to
deduce the existence of points $t_+\in (0,\d/100]$ and $t_-\in
[-\d/100,0)$ which lie in $Z$ where all $\{\phi_j\}$ are  distinct.
Since the set of such points is open and since $\L$ has at most countably
many elements, there is no obstruction to choosing $t_{\pm}$ and then
$s_0 >0$ so that the line segments where $t=t_{\pm}$, \ $u\in [0,s_0]$
miss all points of $\pi(\L)$ and such that the set of functions $\{b_j\}$
are distinct where $t=t_{\pm}$ and $u\in (0,s_0]$.

\medskip
{\bf Step 3}\qua Now introduce the function $\D$ on
$\pi^*C|_A$ as in Step 5 of Part I, and denote the real part of $\Delta$
by $\underline{\phi}$ and the imaginary part by $\underline{b}$.
This step considers the behavior of $d\underline{\phi}$ where
$\underline{b}=0$. The next lemma summarizes the results.

The statement of the lemma introduces the notion of an embedded graph
$\G$ in a manifold. For the present purposes this term is defined as
follows: Call the manifold $X$.  Then $\G$ is a locally compact subset of
$X$ with certain additional properties. First, $\G$ has a distinguished
subset $\G_v$ which is a locally finite set of points. (Each point $\G_v$
is called a vertex.)  Second, $\G-\G_v$ is a locally finite collection of
properly embedded, open intervals in $X-\G_v$. (The closure of each
component of $\G-\G_v$ in $X$ is called an edge.)  Third, each vertex of
$\G$ has a neighborhood in $X$ whose intersection with $\G$ consists of a
finite union of properly embedded images of [0,1) by an embedding which
sends 0 to the vertex in question. Moreover, these embedded intervals
intersect pairwise only at the given vertex.

\noindent {\bf Lemma A.4}\qua {\sl {Let}
$\G\subset\pi^*C-(C\cup\pi_+^{-1}(\L)\cup\pi_-^{-1}(\L))$
{denote the zero set of the function $\underline{b}$.
Then $\G$ is an embedded graph with the following properties.}
\begin{itemize}
\item {The vertices of $\G$ are the $\underline{b}=0$ critical
points of $\underline{b}$ and there are at most finitely many of these in
any subset of $\pi^*C-(C\cup \pi_+^{-1}(\L)\cup \pi_-^{-1}(\L))$
which has compact closure.}

\item {The one form $d\underline{\phi}$ restricts without
zeros to each edge of} \ $\G$.

\item {At each point of $\G-\G_v$, the two form
$d\underline{\phi}\wedge d\underline{b}$ orients $\pi^*C$ so that $\pi$
is an orientation preserving map to $A$.}

\item {There are an even number and more than two edges of
$\G$ ending in each critical point of $\underline{b}$. These edges are
oriented so that $d\underline{\phi}$ increases moving towards the vertex
on precisely half of these edges.}
\end{itemize}}

\medskip{\bf Proof of Lemma A.4}\qua The fact that $\G$ is a graph whose vertices
are the critical points of $d\underline{b}$ follows from the assertion
that $d\underline{b}$ defines a real analytic function on the complement
in $\pi^*C-C$ of $\pi_+^{-1}(\L)\cup\pi_-^{-1}(\L)$. To prove this claim,
consider a small radius disk $D\subset A$ where $u > 0$ which avoids all
points of $\pi(\L)$. Then $C|_D$ consists of $p$ disjoint copies of $D$,
and each $\eta\in\{\eta_j\}$ is a bonafide function on $D$ which obeys
(A.2). It then follows from this last point with standard elliptic
regularity results (as in Chapters 5 and 6 of \cite{Mo}) that each such $\eta$
defines a real analytic function on $D$.

Meanwhile, $\pi^*C|_D$ consists of $p^2$ copies of $D$ with $p^2-p$
copies being disjoint from $C$. With the preceding understood, let
$D'\subset\pi^*C-C$ denote  one of the latter copies of $D$.
Then $w$ pulls back via $\pi$ as a bonafide coordinate of $D'$ and then
the analyticity of each $\eta_j$ as a function of $w$ implies that
$\Delta$ is also an analytic function of $w$.

To see that $\underline{b}$ has at most finitely many critical points on
$D'$, observe that there will be an infinite number of such critical
points only if $\underline{b}=$ constant on $D'$. Thus
$\underline{b}=$ constant on the component of $\pi^*C$ which contains
$D'$. And this last possibility is ruled out by the assumptions in Step 1.

The analysis of the behavior of $d\underline{\phi}$ on $\G$ requires the
observation that $\Delta$ on $D'$ obeys an equation having the schematic
form
\begin{equation}
\bar\p\Delta +\a_{1-}\dt\p\Delta +{\cal R}\dt\underline{b}=0  ,
\end{equation}
where $\a_{1-}\equiv\a_1(\pi_-^*b)$ and where ${\cal R}$ is a real analytic
function on $D'$. In particular, at points of $\G$, the complex function
$\Delta$ obeys the equation $\bar\p\Delta +\a_1\dt\p\Delta=0$.
Since $|\a_1|\leq\z\dt u^2$, it follows that where $u$ is small, 
$|\bar\p\Delta|^2 << |\p\Delta|^2$, which implies both the second and
third points of Lemma A.4.

The final point of (A.19) can be deduced using (A.19) and the analyticity
of $\Delta$ to write down a local model for $\Delta$ near a critical
point where $\underline{b}=0$. Indeed, it follows from (A.19) that near
such a point, $\Delta=\underline{\phi}_0+c\dt(w^m-\a_1(b_0)\dt\bar w^m)
+{\cal O}(|w|^{m+1})$ for some $m > 1$. (Here $\underline{\phi}_0$
is the value of $\underline{\phi}$ at the critical point, and $b_0$ is
the value of $\pi_-(b)$ at the critical point.) It then follows from this
last expression that the condition im$(\Delta)=0$ defines exactly $2m$
embedded line segments emanating from the critical point, these being
approximately the $2m$ rays where the $w^m$ is real. The question of the
orientations is left to the reader to verify from this local expression
for $\Delta$.

\medskip
{\bf Step 4}\qua This step considers the behavior pf $\G$
near the points of the form $(w_0,\eta_0,\eta'_0)$ where
$w_0\in\pi(\L)$ and $\eta_0\neq\eta'_0$. The first observation in this
regard is that $\G$ will miss this singular point of $\pi^*C$ unless the
imaginary parts of $\eta_0$ and $\eta'_0$ agree, so suppose that this is
the case. Introduce from Lemma A.1 the complex parameter $\tau$ on an
irreducible component of a neighborhood of $(w_0,\eta_0,\eta'_0)$ in
$\pi^*C$. It follows from (A.17) and the bound in Lemma A.3 on
$|d\kappa|$ by $\z|\ln(\tau)|$ that there are no $\underline{b}=0$
critical points of $\underline{b}$ in some small neighborhood of 
$(w_0,\eta_0,\eta'_0)$. These same points from Lemma A.3 also imply that
there is some non-zero, and even number of edges of $\G$ whose closure
contains $(w_0,\eta_0,\eta'_0)$, and that half of these point away from
$(w_0,\eta_0,\eta'_0)$ with the $d\underline{\phi}>0$ orientation, while
half point towards
$(w_0,\eta_0,\eta'_0)$ with this orientation.

\medskip
{\bf Step 5}\qua This step considers the behavior of $\G$
near points $(w_0,\eta_0,\eta'_0)$ with $w_0\in\pi(\L)$. For this
purpose, let $\tau$ be the local parameter from Lemma A.1 on an 
irreducible component of a neighborhood of $(w_0,\eta_0,\eta'_0)$ in
$\pi^*C$. Suppose, in addition, that this irreducible component does not
lie in $C$. Once again, let $r >0$ be small and let $\g$ denote the
circle where $|\tau| =r$. And, for almost all $r$, the circle $\g$ will
intersect $\G$ transversely. (In fact, it follows from (A.17) and that
bound on $|d\kappa|$ in Lemma A.3 that $\g$ and $\G$ are transverse
whenever $r$ is sufficiently small.)

Now let $\kappa,\kappa'$, \ $\l,\l'$ and $a$ be as in Lemma A.3. Then it
follows from Lemma A.3 and (A.17) that there are as many segments of $\G$
which intersect $\g$ and point outward with the $d\underline{\phi}>0$
positive orientation as point inward. Furthermore, there is at least one
such segment pointing in each direction with the $d\underline{\phi}>0$
orientation.

\medskip
{\bf Step 6}\qua  This step extends the graph $\G$ to allow
for vertices and edges on the $u=0$ line. For this purpose it is
necessary to first digress to consider the behavior of $\G$ near points
in the $u=0$ line which do not lie in $Z_p$. Let $t_1$ be such a point.
It follows from the induction hypothesis that $C$ is smooth near any such
point and thus over some semi-disk neighborhood $D'=\{(t,u):
u\geq 0$ and $u^2+(t-t_1)^2 <r^2\}$ of such a point, $C$ consists of some
$p$ sheets, with each being diffeomorphic to $D'$ via the projection
$\pi$. In particular, the coordinates $(t,u)$ give coordinates on each
such sheet. More to the point, over $D'$, the functions $\{\eta_1,\dots
,\eta_p\}$ are distinguishable and each is analytic.

Furthermore, according to Assertion 3(c) of Proposition 7.1, each such
$\eta_j$ has the form $\eta_j=\phi_{0j}(t)+i\dt 4^{-1}\phi'_{0j}(t)\dt
u+{\cal O}(|u|^2)$ where $\phi_{0j}(t)$ is an analytic function of $t$ near
$t_1$. This last observation implies that $\Delta$ restricts to a
component of $\pi^*C|_{D'}$ to have the form
$\Delta=\underline{\phi}_0(t)-8^{-1}\underline{\phi}''_0\dt u^2+i\dt
4^{-1}\underline{\phi}'_0(t)\dt +{\cal O}(u^3)$, where $\underline{\phi}_0$
is an analytic function of $t$ near $t_1$. Note that
$d\underline{\phi}\wedge d\underline{b}=
4^{-1}|\underline{\phi}'_0|^2dt\wedge du+{\cal O}(u)$ and thus even where
$u=0$, this form gives the same orientation for $\pi^*C$ as does
$dt\wedge du$.

With the preceding understood, agree to extend the graph $\G$ to the
portion of the closure $\pi^*C|_A$ which lies over the complement of
$Z_p$ in the $u=0$ line. The extra vertices are the points over the complement of
$Z_p$ in the $u=0$ line where $\underline{\phi}'_0=0$. The interiors of
the extra edges are the component intervals of the set of points where
$u=0$ both $t\neq 0$ and also $\underline{\phi}'_0\neq 0$.

\medskip
{\bf Step 7}\qua This step digresses to consider the
behavior of $\G$ near one of these new vertices. In the simplest case, 
$\underline{\phi}_0$ has a non-degenerate critical point over $t_1\not\in
Z_p$ on the $u=0$ line. Near such a critical point, Assertion 3(c) of
Proposition 7.1 has
\begin{equation}
\Delta =c\dt[(t-t_1)^2-4^{-1}\dt u^2+i\dt 4^{-1}(t-t_1)\dt u]+
{\cal O}(u^3) ,
\end{equation}
with $c\neq 0$. Thus $\G$ has three edges which end at the vertex
$(t_1,0)$, these being the segments where
\begin{itemize}
\item $u=0 \ {\mbox{and}} \ t < t_1$
\item $u=0 \ {\mbox{and}} \ t > t_1$
\item $u>0 \ {\mbox{and}} \ t = t_1\  
{\mbox{(to leading order in $u$)}}$.\itemnum{A.21}
\end{itemize}

The $d\underline{\phi}$ orientations for the latter are equivalent to that
given, respectively, by --sign$(c)dt$, sign$(c)dt$ and --sign$(c)du$.
Thus the positive $d\underline{\phi}$ orientation for the first two
segments make both point either towards or away from $t_1$ together, and
then the third segment points opposite, away or towards $t_1$,
respectively.

The following lemma states the situation in general.

\medskip{\bf Lemma A.5}\qua{\sl{Fix a vertex of $\G$ over the $u=0$ line. This
vertex lies over a point $t_1\not\in Z_p$ and is a critical point of 
$\underline{\phi}_0$. If the order of vanishing of $d\underline{\phi}_0$
at $t_1$ is even, then there are at least four and an even number of
edges of $\G$ which are incident to this vertex. If the order of
vanishing is odd, there are at least three and an odd number of edges of
$\G$ which are incident to this vertex. In either case, two of these
project to the $u=0$ line, one to where $t<t_1$ and the other to where $t
>t_1$. The remaining edges lie, except for their vertices, where $u >0$.
Furthermore, for $r > 0$ but small, each such edge intersects the half
circle where $(t-t_1)^2+u^2=r^2$ and $u\geq 0$ transversely in exactly one
point. By following this half circle from its start at $u=0$ and
$t=t_1-r$ to its end at $u=0$ and $t=t_1+r$, the edges are met
consecutively, and consecutive edges point in opposite directions (either
towards or away from $t_1$) when oriented to have} $d\underline{\phi}_0
>0$. {(If the order of vanishing of
$d\underline{\phi}_0$ at $t_1$ is $N$, then it can be proved that there
are exactly $N+2$ edges of $\G$ which are incident on the vertex $t_1$.)}}

\medskip{\bf Proof of Lemma A.5}\qua Introduce polar coordinates 
$(r,\theta)\in [0,1)\x [0,\pi]$ near $t_1$ by writing $t-t_1=r\cos\theta$
and $u=r\sin\theta$. Then, as $\Delta$ is analytic near $t_1$, it can be
written as
\begin{equation}
\Delta =\underline{\phi}_0(t_1)+ r^{N+1}\dt (f(\theta)+i\dt g(\theta))
+{\cal O}(r^{N+2}).
\end{equation}
Here, $f$ and $g$ are real analytic functions on $[0,\pi]$ which, because
of (A.2), are constrained to obey
\begin{itemize}
\item $\sin^{-3}\theta (\sin^3\theta f_{\theta})_{\theta}
 +(N+1)(N+4)f=0$
\item $f_{\theta}(0)=f_{\theta}(\pi)=0 $
\item $g=f_{\theta}/(N+4)$.\itemnum{A.23}
\end{itemize}

Because $f$ is constrained here to be an analytic function on $[0,\pi]$,
the first line in (A.23) implies the following:
\begin{itemize}
\item The coefficients in a power series solution for $f$ about any point
are determined by the value of $f$ and its derivative at that point.
\hspace\fill (A.24)\end{itemize} 
In particular, this last point implies that there is a unique (up to
multiplication by real numbers) solution to the first two constraints
(A.23). The fact that there is a non-trivial solution is not obvious, but
such is indeed the case. It is left to the reader to verify that  the
first two points in (A.23) are satisfied by a polynomial of order $(N+1)$
in $\cos\theta$ which is either even or odd (depending on the parity of
$N+1$) under the involution $\theta\to\pi-\theta$. 

Note that (A.24) also precludes both $f_{\theta\theta}$ and $f_{\theta}$ 
vanishing simultaneously.  This implies that the zeros of $g$ are
non-degenerate. Then (A.22) implies that for small, but positive $r$,
the graph $\G$ has no vertices save that at $t_1$, and it implies that
each edge of $\G$ which is incident on this vertex has exactly one
intersection point with the half circle of constant $r$, and that the
intersection there is transverse.

Here is an argument for the equality between the parity of the number of
zeros of $f_0$ and the parity of the order of vanishing of
$d\underline{\phi}_0$ at $t_1$: Consider the case where the parity of the
latter is odd (the even case is argued along similar lines). Odd order of
vanishing of $d\underline{\phi}_0$ implies that $f(0)=f(\pi)$. On the
other hand, the analyticity of $f$ and the first two points of (A.23)
imply that $f_{\theta\theta}$ and $f$ have  opposite signs at both 0 and
$\pi$. Thus, $f_{\theta\theta}(0)=f_{\theta\theta}(\pi)$ also.
In as much as $f_{\theta}$ has transverse zeros, this last condition
implies that $f_{\theta}$ must have an odd number of zeros in the
interval $(0,\pi)$.

The following is an argument for the existence of at least three zeros of
$f_{\theta}$ (counting those at $\theta=0$ and $\theta=\pi$): When
$f$ and $f'$ are solutions to the first two lines of (A.23) which
correspond to different values of $N$, then integration by parts shows
that $f_{\theta}$ and $f'_{\theta}$ are orthogonal with respect to the
measure $\sin^3\theta \ d\theta$ on $[0,\pi]$. As the $N=0$ solution to
the first two lines of (A.23) is $\cos\theta$, this orthogonality implies
that each $N>0$ solution $f$ obeys $\int f_{\theta}\sin^4\theta \
d\theta=0$. Hence $f_{\theta}$ must change sign on the interval $[0,\pi]$
and so it has at least three zeros.
 
\medskip
{\bf Step 8}\qua This step serves as a digression to
introduce the family of maps $\psi$ from either [0,1] or [0,1) or (0,1]
into $A\times {\Bbb C}\times{\Bbb C}$ with the following properties:
\begin{itemize}
\item Im$(\psi)\subset{\mbox{closure}}(\G)$.
\item $\psi$ is 1--1 onto its image.
\item Let $e\subset\G$ be the interior of an edge.
 Then $\psi^{-1}(e)$ is either empty or connected, and in the latter case, 
$\psi$ restricts to $\psi^{-1}(e)$ as a smooth, orientation preserving
embedding.  Here, the orientation on $e$ is defined so that
$d\underline{\phi} >0$. \hspace\fill (A.25)
\end{itemize}

\noindent A map $\psi$ which obeys the constraints in (A.25) will be
called an oriented, graphical path.
\addtocounter{equation}{1}

\medskip
{\bf Step 9}\qua To begin the discussion here, introduce
$t_{\pm}$ and $s_0$ from Step 2. If $s_0$ is sufficiently generic, then
$\pi(\L)$ has empty intersection with the line segment with
$t\in [t_-,t_+]$ and where $u=s_0$. Furthermore, if $s_0$ is  generic,
then $\G$ has no vertices which project to the $u=s_0$ line segment, nor
does it have non-transversal intersections with the inverse image of this
segment. Thus, with $s_0$ generic, the subset of $\G$ where $u\leq s_0$
intersects the level set $u=s_0$ in some finite number, say $N$ points,
Also by Step 2's assumption, $t\in (t_-,t_+)$ on $\G$.

The first claim is:

\medskip{\bf Lemma A.6}\qua{\sl{Assume that $t_0$ is an isolated point of $Z_p$.
If $N>0$, then there exists a set  $\{\psi_j\}$ of no more than $N$ and
at least $N/2$ oriented graphical paths in $A\times\Bbb C\times\Bbb C$
with the following properties:}

\begin{itemize}
\item {The map $\psi_j$ sends any endpoints of its domain to the
$u=s_0$ line.}

\item {If the domain of $\psi_j$ is [0,1) or (0,1], then the
closure of the image of $\psi_j$ consists of the  union of the image of
$\psi_j$ with a non-empty set of points which map by $\pi$ to $Z_p$.}

\item {If the image of $\psi_j$ intersects a vertex where
$t=t_1$ on the $u=0$ line, then precisely one of these possibilities
occurs:}

\leftskip 25pt

\item[\bf a\rm)] $\psi_j$ {both enters and exits this vertex along edges where} \ $u
> 0$.

\item[\bf b\rm)] $\psi_j$ {both enters and exits this vertex along edges where} \
$u\equiv 0$.

\item[\bf c\rm)] $\psi_j$  {enters  along the $u=0$ edge where $|t|>|t_1|$ and it
exits an edge where} \ $u > 0$.

\item[\bf d\rm)] $\psi_j$ {enters along an edge where $u > 0$ and
exits along the $u=0$ edge where} \ $|t|<|t_1|$.

\leftskip 0pt

\item {The image of $\psi_j$ intersects at most one $\Delta =0$
point of $\pi^{-1}(\L)$. And, if such a point is in the image of
$\psi_j$, then the domain of $\psi_j$ is [0,1].}

\item {Every point of $\G$ on the $u=s_0$ line is in the image
of some} \ $\psi_j$.

\item {When $j\neq i$, then the image of $\psi_j$ intersects
that of $\psi_k$, if at all, only at vertices or at points of}
$\pi^{-1}(\L)$.
\end{itemize}}

\medskip
The remainder of this step is occupied with the proof of Lemma A.6.

\medskip
{\bf Proof of Lemma A.6}\qua The maps $\{\psi_j\}$ will be constructed
consecutively starting with $\psi_1$. To construct $\psi_1$, introduce
the subgraph $\G_1\subset\G$ which is obtained from the latter by
deleting the following subset of $u\equiv 0$ edges:
\begin{itemize}
\item An edge where $t > 0$ is deleted if the
$d\underline{\phi} >0$ orientation {agrees} with the $dt >0$
orientation.
\item An edge where $t < 0$ is deleted if the
$d\underline{\phi} >0$ orientation {disagrees} with the $dt >0$
orientation.\hspace\fill (A.26)
\end{itemize} 
Note that $\G$ and $\G_1$ agree where $u >0$.

The graph $\G_1$ is an example of a subgraph $\G'\subset\G$ with the
following properties:
\begin{itemize}
\item Every vertex of $\G'$ intersects an even number of edges,
half of which are incoming and half of which are outgoing when oriented
so that $d\underline{\phi} >0$.  If a vertex is on the $u=0$ line, say
with $t$ value equal to $t_1$, then an incident $u\equiv 0$ edge is
outgoing if and only if $|t|\leq |t_1|$ on this edge.
\item Every point in $\pi_+^{-1}(\Lambda)\cup\pi^{-1}(\Lambda)$
intersects the closure in $A\times\Bbb C\times\Bbb C$ of an even number
of edges of $\G'$. Furthermore, half of these are incoming and half are
outgoing when they are oriented so that $d\underline{\phi} >0$.
\item The closure of $\G'$ in $A\times\Bbb C\times\Bbb C$ 
is the union of the subgraph, a set of points which map by $\pi$
to $Z_p$, and a subset of $\pi^{-1}(\Lambda)$.
\item The  subgraph $\G'$ has non-empty intersection
with the $u=s_0$ line.\hspace\fill (A.27)
\end{itemize} 
(The first point here follows from Lemmas A.4 and A.5, and the second point
from the discussions in Steps 4 and 5.)

Let $\G'\subset\G$ be a subgraph which obeys the points of (A.27).
The subsequent construction in this step produces an oriented, graphical
path  $\psi$ which satisfies the first four points in Lemma A.6 and which
lies in the closure of $\G'$. To begin, select a point where $\G'$
intersects that $u=s_0$ line segment and follow the corresponding edge of
$\G'$ to where $u <s_0$. This begins to trace out the image of graphical
path $\psi\co [0,\d')\to\G'$ for some $\d' >0$ which is either orientation
preserving or reversing, depending on whether the edge is oriented inward
or outward along the $u=s_0$ segment. For the sake of argument, suppose
that this $\psi$ is orientation preserving. (If not, reparametrize the
interval [0,1] by inverting about 1/2.) The map $\psi$ can be extended
indefinitely by walking down the given edge unless one of the following
occur:
\begin{itemize}
\item\quad The $u=s_0$ segment is intersected.
\item\quad   A $u >0$ vertex of $\G'$ is reached.
\item\quad  The image of $\pi\dt\psi$ limits to a point in
$\pi_+^{-1}(\Lambda)\cup\pi^{-1}_-(\Lambda)$.
\item\quad  A $u=0$ vertex of $\G'$ is reached.\hspace\fill (A.28)
\end{itemize} 
If the first case in (A.28) occurs, then stop because up to
reparametrization of the domain of definition, a graphical path $\psi$
obeying the first four points of Lemma A.6 has been constructed.

Assume that the second point of (A.28) occurs. In this case, the image of
$\psi$ enters the  vertex on an inward oriented edge, and therefore the
map $\psi$ can be extended to a larger domain so that the extension
leaves the given vertex on one of the outward oriented edges. (Use the
first point of (A.27) here.) This extends $\psi$ as an oriented, graphical
path.

Likewise, if the third point of (A.28) occurs, then $\psi$ can also be
extended as an oriented, graphical
path. This follows from the second point of (A.27).

If the fourth point of (A.28) occurs, let $t_1$ denote the $t$ coordinate
of this vertex. Note that $\psi$ entered this vertex from an edge which
has $u >0$ on its interior. If there is an outgoing edge with $u >0$ on
its interior, then $\psi$ can be extended by continuing along this end.
If not, then there is an outgoing edge with $u\equiv 0$ and that one has
$|t|<|t_1|$. (Use the second point in (A.27).) Extend $\psi$ by continuing
along this last outgoing edge.

Now one can continue to extend $\psi$ by walking along the interior of
the edge just chosen until one of the possibilities in (A.28) reoccurs. If
the first point occurs, then stop and the resulting $\psi$ is an
oriented, graphical path in the closure of $\G'$ which obeys the first
four points of Lemma A.6. If the second or third points in (A.28) occur,
extend $\psi$ as before. If the last point in (A.28) occurs, then $\psi$
can also be extended as described previously except in the case where the
incoming edge has $u=0$, as this case was not previously considered. Note
that in this case, the entering edge has $|t|>|t_1|$ due to the chosen
direction for leaving a $u=0$ vertex. In any case, if there is an
outgoing edge from this vertex which has $u >0$ on its interior, then
extend $\psi$ by preceding along this vertex. If not, then the other
$u=0$ edge (which has $|t|<|t_1|$) is outgoing, and one can extend
$\psi$ by proceeding along the latter.

These last steps can be iterated repeatedly to extend $\psi$ as an 
oriented, graphical
path  indefinitely unless point one in (A.28) occurs. In this case,
$\psi$ is an 
oriented, graphical path in the closure of $\G'$ which obeys the first
four points of Lemma A.6. (Note that the image of $\psi$ can intersect at
most one $\Delta =0$ point in $\pi^{-1}_+(\L)\cup\pi^{-1}_-(\L)$ because
$\underline{\phi}={\mbox{Re}}(\Delta)$ is zero at such points but is
increasing along the image of $\psi$. Also, as $\underline{\phi}$ is 
increasing along the path, its image can intersect any given vertex at
most once.) If the first point in (A.28) never occurs, then the
construction just described  proceeds without end to construct a map,
$\psi$, from [0,1) into $A\times\Bbb C\times \Bbb C$. However, in this
last case, the closure of the image will be the union of the image of
$\psi$ with a non-empty subset of $Z_p$. Thus, in this case, the first
four points of Lemma A.6 are also obeyed. (In this regard, note that if
the domain of $\psi$ is equal to [0,1), then no $\Delta =0$ points are in
the $\psi$'s image as $\Delta=0$ at the points which project by $\pi$ to
$Z_p$ and $\underline{\phi}={\mbox{re}}(\Delta)$ is increasing along the
image of $\underline{\phi}$.)

With the preceding construction understood, apply it to the case where
$\G'=\G_1$ and let $\psi_1$ denote the corresponding $\psi$. Now consider
the construction of $\psi_2$ and, sequentially, $\psi_j$ for $j\geq 3$.
The construction here is inductive. Given that $\psi_{j-1}$ for $j\geq
2$, as has been constructed, set $\G_j$ to be the subgraph of $\G_{j-1}$
which is obtained by deleting all interiors of edges which intersect the
image of $\psi_{j-1}$. The construction just described of $\psi_{j-1}$
insures that $\G_j$ obeys all the points of (A.27) except possibly the
last one. If the last point is also obeyed, then repeat the recipe of the
preceding paragraphs with $\G'=\G_j$ and take $\psi_j$ to be the
corresponding $\psi$. Since each $\psi_j$ maps into $\G_j$, the final
requirement for the set $\{\psi_j\}$ is automatically satisfied.

\medskip
{\bf Step 10}\qua Consider the following proposition:

\medskip\noindent{\bf Proposition A.7}\qua {\sl{Suppose that $Z-Z_p$ is dense
in $Z$ near the point $t_0$. Then $\L$ is finite.}}

\medskip The remainder of this step proves the preceding proposition.

\medskip{\bf Proof of Proposition A.7}\qua Assume, to the contrary, that $\L$ has
an infinite number of points. This assumption will be seen to imply a
contradiction to the fact that $\underline{\phi}$ is increasing along any 
oriented, graphical path.
To obtain such a contradiction, it is necessary to digress for a moment
to introduce the subgraph $\G'\subset\G$ which is obtained from the
latter by deleting various edges and vertices. A vertex $v\in\G$ is
missing from $\G'$ when all edges which intersect $v$ either satisfy one
of the conditions of (A.26) or else intersect the image of some $\psi_j$
in their interior.  Meanwhile, an edge $e\in\G$ is missing from $\G'$ if
either of the conditions of (A.26) hold, or if the interior of said edge
intersects the image of some $\psi_j$. Note that $\G'\neq\emptyset$ if
$\L$ is infinite because each $\psi_j$ intersects no more than one point
in the set $\pi^{-1}_+(\L)\cup\pi^{-1}_-(\L)$ where $\Delta =0$, and
there are infinitely many such points. Furthermore, the local picture of
$\G$ near a $\Delta =0$ point  as given in Step 5 indicates that each
such point is in the closure of a non zero and even number of edges of
$\G$. In particular, $\G'$ here satisfies all but the last condition in
(A.27).

With the preceding  understood, choose a $\Delta =0$ point which is in
the closure of an edge from $\G'$. It follows from the discussion in Step
5 that there is at least one outgoing edge and one incoming edge (both
with the $d\underline{\phi}>0$ orientation) of $\G'$ which have the given
point in their closure. Select such an outgoing edge and proceed along it
starting from the given $\Delta =0$ point. Then the construction in the
previous step produces a map $\psi$ from [0,1) into $A\times\Bbb C\times
\Bbb C$ which is a graphical embedding whose image is in the closure of
$\G'$. By construction, $\psi$ sends 0 to the chosen $\Delta =0$ point.
Note that the first point in (A.28) cannot occur becuase $\G'$ is disjoint
from the $u=s_0$ line. It then follows from the construction that the
closure of the image of $\psi$ is the union of the image of $\psi$ with
some non-empty subset of points which project by $\pi$ to $Z_p$. However,
this is impossible as $\Delta=0$ at any such point, while
$\underline{\phi}$ is already zero at $\psi(0)$ and then increases along
the image of $\psi$.

\subsection*{Part III\qua The density of $Z-Z_p$ in $Z-Z_{p+1}$}

The purpose of this last section of the appendix is to prove the
following.

\noindent {\bf Proposition A.8}\qua{\sl {Let $t_0$ be a regular point of $Z$
which lies in $Z-Z_{p+1}$. Then $Z-Z_p$ has non-empty intersection with
any given neighborhood of} $t_0$.}

\medskip Note that this proposition justifies the assumption in (A.18).

The proof of this proposition is broken into fourteen steps.

\addtocounter{equation}{4}
\medskip
{\bf Step 1}\qua This step outlines the proof.
To begin, reintroduce the domain $A$ and then the coordinates
$w=t+i\dt u$ and $\eta=\phi +i\dt b$ as in Step 1 of part I.
Thus $t_0$ corresponds to the point $w=0$. Next, introduce the data
$\{\eta_j\}$ as in Step 1 of part I and let $\underline{\eta}\equiv
p^{-1}\sum_j\eta_j$. Note that  $\underline{\eta}$ is a bonafide, complex
valued function on $A$. Set $\{\underline{\eta}_j\equiv\eta_j-
\underline{\eta}\}$ and observe that all elements in
$\{\underline{\eta}_j\}$ vanish at each point of $A_p$.

Fix $s_0 >0$ but much less than the number $\d$ which defines the size of
the half-square $A$ and let $D\subset A$ denote the $u\geq 0$ radius
$s_0$ half disk with center at $w=0$. For positive integer $n$, define
data $\{\underline{\eta}_{j,n}\}$ on $D$ by dilating and rescaling
$\{\underline{\eta}_j\}$ according to the rule
\begin{equation}
\underline{\eta}_{j,n}(w) =\underline{\eta}_j(2^{-n}w)/\lambda_n .
\end{equation}
Here $\l_n$ is a suitably chosen real number. Up to this factor
$\{\underline{\eta}_{j,n}\}$ is obtained  from
$\{\underline{\eta}_j\}$ by pulling the latter back from the ball of
radius $2^{-n}s_0$. In any event, the data 
$\{\underline{\eta}_{j,n}\}$ defines a continuous map from the radius
$s_0$ half-disk about $w=0$ to Sym${}^p\Bbb C$.

This construction in (A.29) should be done for each positive integer $n$
to obtain a sequence $\{\{\underline{\eta}_{j,n}\}\}_{n=1,2,\dots}$.
For the purposes of the next lemma, this sequence should be thought of as
a sequence of maps into the space Sym${}^p\Bbb C$. And, for the purposes
of the next lemma, Sym${}^p\Bbb C$ will be identified with ${\Bbb C}^p$
via the homeomorphism from the latter to Sym${}^p\Bbb C$ which associates
the $p$ roots of a monic, $p^{\mbox{\footnotesize{th}}}$ order polynomial
on $\Bbb C$ with the coefficients of its $p$ lower order terms.

With the preceding understood, here is the key property of this sequence:

\noindent {\bf Lemma A.9}\qua{\sl {The numbers $\{\l_n\}_{n=1,2,\dots}$
in} (A.29) {can be chosen so as to be strictly decreasing to zero and
so that the resulting sequence
$\{\{\underline{\eta}_{j,n}\}\}_{n=1,2,\dots}$  has an infinite
subsequence which converges weakly in the Sobolev $L^2_1$ topology and
strongly in the $C^0$ topology on the $u\geq 0$ radius $s_0/16$ half-disk
with center $w=0$ to a continuous non-trivial map into} \ Sym${}^p\Bbb C$.}

\medskip Let $D'$ denote the radius $s_0/16$ half-disk in the lemma.
Also, let $\{\{\underline{\eta}_{j,\infty}\}_{1\leq j\leq p}$ denote the
limit in Lemma A.9. Then let $C'\subset D'\times\Bbb C$ denote the closure
of the set of points of the form $(w,\eta)$ where both
$\eta\in\{\underline{\eta}_{j,\infty}\}|_w$, and $\eta\neq 0$.
Now use this $C'$ to define the function $\tau$ on $D'$ whose value at
$w$ is $\tau(w)\equiv\prod_{(w,\eta)\in C'}\eta$.
Here is the key property of $\tau$:

\noindent {\bf Lemma A.10}\qua{\sl{The function $\tau$ as defined above is not
identically zero, it is continuous and Sobolev class $L^2_1$.
Furthermore, it obeys an equation having the schematic form}
\begin{equation}
\bar\p\tau +\frac{3i}{4u} \ p'(1-{\cal P})\tau=0,
\end{equation}
{where $p'\in\{ 1,\dots ,p\}$ and where ${\cal P}$ is continuous almost
everywhere on $D'$ and has norm less than or equal to 1.}}

\medskip The preceding two lemmas are proved below starting in Step 3, so
accept them for now. 

The proof of Proposition A.8 proceeds from this point with a study of the
function $\tau$. In particular, consider:

\noindent {\bf Lemma A.11}\qua{\sl{Let $\tau$ be a non-trivial solution to}
(A.30) {on $D'$ which is continuous and Sobolev class $L^2_1$.
Then $\tau$ has isolated zeros where $u >0$ and $\tau\neq 0$ on a dense,
open subset of the $u=0$ line near the origin.}}

\medskip This last lemma is proved below in Step 2.

With the previous two lemmas understood, the proof of Proposition A.8 ends
with the following two remarks: First, given $\e >0$, there is, according
to Lemma A.11, a point $t'$ on the $u=0$ line with distance $\e$ or less
from the origin where $\tau\neq 0$.  Thus, $|\tau| > c > 0$ at this
point. Second, according to Lemma A.9, there is some (infinitely many)
value of $n$ for which there is a $j$ with
$|\underline{\eta}_{j,n}| > c/(2(1+c))$ at the point $t'$.
Thus at $2^{-n}t'$ on the $u=0$ line, one has
$|\underline{\eta}_j| > c\l_n/(2(1+c))$. Hence the point $2^{-n}t'$ on
the $u=0$ line is not in $A_p$.

\medskip
{\bf Step 2}\qua
This step contains the proof of Lemma A.11.

\medskip{\bf Proof of Lemma A.11}\qua Write $\tau=u^{-3p'/2}\l$. Then (A.30)
implies the following equation for $\l$:
\begin{equation}
\bar\p\l +\frac{3i}{4u} \ p'(1-{\cal P})\l=0.
\end{equation}
This last equation can be used to analyze $\l$ by writing the latter as
$e^{\kappa}\l_0$, where $\kappa$ is a solution to the equation
\begin{equation}
\bar\p\kappa = \frac{3i}{4u} \ p'{\cal P} .
\end{equation}
and where $\l_0$ is holomorphic where $u >0$ on $A$. Of particular
interest is the solution to (A.32) which is given by the integral
expression
\begin{equation}
\kappa|_w =\frac{3i}{4\pi} \ p'\int_{\Bbb C}\underline{\theta}
((w-w')^{-1}-(w-\bar w')^{-1}){u'}^{-1}{\cal P} d^2w' .
\end{equation}
Here $\underline{\theta}$ denotes the characteristic function of the 
radius $s_0/16$ half disk on which $\tau$ is defined.

This last expression for $\kappa$ can be used to estimate $|\kappa|$.
For this purpose, first write $(w'-w)^{-1}-(\bar w'-w)^{-1}=
2u'/((w'-w)(\bar w-w))$ and then take absolute values to find that
\begin{equation}
|\kappa|_w\leq 3p'(2\pi)^{-1}\int_{\Bbb C} \underline{\theta}
1/(|w'-w| \ |\bar w'-w|) d^2w' .
\end{equation}
In particular, suppose that the $t$ coordinate of $w$ is less than
$s_0/64$.  Then it is an exercise to use (A.34) to bound
\begin{equation}
|\kappa|_{(t,u)} \leq \z_1 +\textstyle{\frac 32} p'|\ln(u)|.
\end{equation}
where $\zeta_1$ is independent of $t$ and $u$.  With regard to the
derivation of (A.35), note that the $\ln(u)$ term in (A.35) arises from
that part of the integral in (A.34) where the coordiantes $t'$ and $u'$
of $w'$ are constrained by $|t'-t|^2 + u'{}^2 \ge 100 u^2$. Indeed,
this portion of the integral differs from
\begin{equation}
3p'(2\pi)^{-1}\int_{\ss \d\geq |w'|\geq 10u\atop0\ss \leq
{\rm arg}(w')\leq\pi} |w'|^{-2} d^2w'
\end{equation}
by a $u$--independent constant.  Likewise, the remaining
portion of the integral is also bounded by a $u$--independent
constant.  Here, this last assertion follows from the fact that
$1/|{\bar w}' - w| \le 1/u$ on $A$.

In any event, with (A.35) understood, it follows that $|\lambda| \ge
\zeta^{-1}u^{3p/2}|\lambda_0|$ where $\zeta \ge 1$ is a fixed,
$u$--independent constant.  On the other hand, because $|\tau| =
u^{-3p/2}|\lambda|$, it follows that
\begin{equation}
|\tau| \geq\z^{-1} |\l_0| .
\end{equation}
To end the story, note that (A.37) implies that $\tau$ can not
vanish on a non-trivial line segment of the $u = 0$ line near $(0,0)$.
Indeed, if it did, then the holomorphic function $\lambda_0$ would
vanish on a non-trivial line segment and thus vanish identically.
Then, $\tau$ would vanish identically, which is assumed not to be the
case.

\medskip
\noindent
{\bf Step 3}\qua  This step and the subsequent steps contain
the arguments which lead to the proof of Lemmas A.9 and A.10.  This step
serves to identify the number $\lambda_n$ in (A.29).

For this last purpose, introduce, for positive and integral $n$, the
half disk, $D(n)$, of radius $2^{n}s_0$ and center at $w = 0$.  Now,
define $\lambda_n$ by the equality
\begin{equation}
\l^2_n\equiv \int_{D(n)} \sum_j (|d\underline{\eta}_j|^2 +
u^{-2}|\mbox{im}(\underline{\eta}_j)|^2).
\end{equation}
Note that the sequence $\{\lambda_n\}$ is strictly decreasing
to zero as advertised in Lemma A.9.

\medskip
\noindent
{\bf Step 4}\qua  Here are some key properties of
$\{\underline{\eta}_{j,n}\}$:

\medskip
\noindent
{\bf Lemma A.12}\qua  {\sl{For each positive integer $n$, define
$\{\underline{\eta}_{j,n}\}$ as in} (A.29) {with $\lambda_n$ as
in} (A.38).  {There is a constant $\zeta \ge 1$ which is
independent of $n$ and is such that for all $n$ sufficiently large,
the following is true}:

\begin{itemize}
\item[$\bullet$] $1 = \int_{D(0)} \sum_j (|d\underline{\eta}_{j,n}|^2 +
u^{-2}|\mbox{im}(\underline{\eta}_{j,n})|^2)$

\item[$\bullet$] $\int_{D(4)} \sum_j (|d\underline{\eta}_{j,n}|^2 +
u^{-2}|\mbox{im}(\underline{\eta}_{j,n})|^2) \le \zeta \int_{D(2)}
|\underline{\eta}_{j,n}|^2$

\item[$\bullet$] {Let $s \in (0,s_0/4)$ and let $D' \subset D(0)$
denote the intersection of $D(0)$ with a disk of radius $s < s_0/4$ and
center in $D(2)$.  Then, $\int_{D'} \sum_j |d\underline{\eta}_{j,n}|^2
\le \zeta \cdot s^{1/\zeta}$}.

\item[$\bullet$] {The function $\theta_n \equiv (\sum_j
|\underline{\eta}_{j,n}|^2)^{1/2}$ is bounded by $\zeta$ on $D(2)$ and
on $D(2)$ satisfies a Holder estimate of the following form:
$|\theta_n(w) - \theta_n(w')| \le \zeta \cdot |w-w'|^{1/\zeta}$ when
both $w$ and $w'$ are in $D(2)$}.

\item[$\bullet$] {Each $\underline{\eta}_{j,n}$ can be viewed as a
well defined function on some neighborhood of all but countably many
$u > 0$ points of $D(1)$.  On such a neighborhood, each
$\underline{\eta}_{j,n}$ obeys an equation of the form}
\begin{equation}
\bar\p\underline{\eta}_{j,n}-\textstyle{\frac 32} u^{-1}
{\mbox{im}}(\underline{\eta}_j)+\s_{j,n}\dt\p \underline{\eta}_j\dt
{\mbox{im}}(\underline{\eta}_{k,n})=0, 
\end{equation}
{where $|\sigma_{j,n}| \le \zeta \cdot 2^{-n}$ and where
$\{\sigma_{j,n}^{k'}\}$ also tend to zero as $n$ tends to infinity, but
in the following weaker sense:  Let $s > 0$ and let $D' \subset D(0)$
be the intersection of $D(0)$ with a disk of radius $s$ in $A$.  Then,
$\int_{D'} \sum_{j,k} |\sigma_{j,n}^k|^2 \le \zeta
s^{1/\zeta}2^{-n/\zeta}$}.
\end{itemize}}

\medskip
\noindent
{\bf Proof of Lemma A.12}\qua  The first point of the lemma follows
directly from the scaling relation.  The proof of the second and third
points are lengthy and given below.  Meanwhile, the fourth point
follows from the fifth point and the first point using Theorem 3.5.2
of \cite{Mo}.  And, the fifth point follows driectly from (A.2) with some
straightforward algebraic manipulations and with the appropriate
rescaling.  (Note that (8.4) and (8.12) have been invoked to obtain
the estimates for $\{\sigma_{j,n}\}$ and $\{\sigma_{j,n}^{k'}\}$.)
The fifth point of the lemma is used in the proof of the second and
third points.

To prove the second point of the lemma, introduce a bump function
$\underline{\chi}$ which is $1$ on $D(4)$, and zero on $D(0)-D(2)$
and which is non-increasing.  Multiply (A.39) by $\underline{\chi}$.
Write the term $\underline{\chi}{\bar \partial}\underline{\eta}_{j,n}$
as ${\bar \partial}(\underline{\chi}\underline{\eta}_{j,n}) - {\bar
\partial}\underline{\chi}\underline{\eta}_{j,n}$ and likewise write
$\underline{\chi} \partial\underline{\eta}_{j,n}$ in terms of
$\partial(\underline{\chi}\underline{\eta}_{j,n})$.  Put the terms
with derivatives of $\underline{\chi}$ on the right hand side.  Then,
square the resulting equation and after an integration by parts, one
obtains the inequality
\begin{eqnarray}
\int_{D(1)} \sum_j
(|d(\underline{\chi}\underline{\eta}_{j,n})|^2
+u^{-2}|{\mbox{im}}(\underline{\chi}\underline{\eta}_{j,n})|^2)-&\hspace{-5pt}
\zeta
\int_{D(1)} \sum_{j,k} |\s^k_{j,n}|^2 \sum_m
|{\mbox{im}}\underline{\chi}\underline{\eta}_{m,n})|^2 \nonumber
\\
\leq\zeta\dt
s^{-2}_0
\int_{D(1)}|\underline{\eta}_{j,n}|^2 .&
\end{eqnarray}
(There are no extraneous boundary terms in the integration by
parts due to the presence of points in $\Lambda$.  The argument here
for the absence of such terms is the same as that given in the proof
of Lemma 8.3 for the absence of similar terms in (8.14).)

The next claim is that the right most term on the left side of (A.40)
obeys the bound
\begin{equation}
\int_{D(1)} \sum_{j,k} |\s^k_{j,n}|^2 \sum_{j'}
|{\mbox{im}}(\underline{\chi}\underline{\eta}_{j',n})|^2 \leq\zeta
2^{-n/\zeta}
\int_{D(1)} \sum_j (|d(\underline{\chi}\underline{\eta}_{j,n})|^2
\end{equation}
where $\zeta \ge 1$ is independent of $n$.  Note that when $n$
is large, this last claim with (A.40) gives the second point of Lemma
A.13.  Meanwhile, the claim in (A.41) follows from Lemma 5.4.1 of \cite{Mo}
using the estimates for the integrals of $\sum_{j,k}
|\sigma_{j,n}^k|^2$ over disks in $D(0)$.  With regard to this appeal
to Morrey's lemma, take for Morrey's domain $G$ the half disk $D(1)$
in the $w$--plane, and take Morrey's function $u$ to equal $(\sum_j
|\mbox{im}(\underline{\chi}\underline{\eta}_{j,n})|^2)^{1/2}$.  (Note
that the norm of the differential of the latter is no greater than
$(\sum_j |d(\underline{\chi}\underline{\eta}_{j,n})|^2)^{1/2}$.)

The proof of the third point of Lemma A.12 requires a preliminary
digression.  To start the digression, take $s$ positive but less than
$s < s_0/4$.  Then, let $S(s)$ denote a $u \ge 0$ half-disk 
with diameter $s$ and center at some point $(t',0)$ on the $u = 0$
line in $D(2)$.  The purpose of this digression is to prove that
\begin{equation}
\int_{S(s)}\sum_j (|d \underline{\eta}_{j,n}|^2+u^{-2}
|{\mbox{im}}(\underline{\eta}_{j,n})|^2
\leq\zeta\dt s^{1/\zeta}  ,
\end{equation}
where $\zeta \ge 1$ is independent of $s$, $n$ and the center
point of the half disk.

What follows next is the proof of (A.42):  Write
$\underline{\varphi}_{j,n}$ for the real part of
$\underline{\eta}_{j,n}$ and $\underline{b}_{j,n}$ for the imaginary
part.  Now, square both sides of (A.39), sum over $j$ and integrate the
result over $S(s)$ to find (after an integration by parts) that
\begin{equation}
\begin{array}{l}
\displaystyle{\int_{S(s)}\sum_j (|d \underline{\eta}_{j,n}|^2+ 9u^{-2}|
\underline{b}_{j,n}|^2+6u^{-1}\underline{b}_{j,n}\p_u b_{j,n}-
6u^{-1}\underline{b}_{j,n}\p_t \underline{\phi}_{j,n})} \vspace{1\jot}
\\
\displaystyle{\leq -\int_{\p S(s)}\sum_j \underline{b}_{j,n}
d\underline{\phi}_{j,n} +\zeta\dt (2^{-2n}\int_{S(s)}
\sum_j |d \underline{\eta}_{j,n}|^2+\int_{S(s)}\sum_{j,k}
|\s^k_{j,n}|\dt \sum_j |\underline{b}_{j,n}|^2) .} \end{array}
\end{equation}
As with the manipulations in Section~8 for the proof of Lemma
8.2, the points in $\Lambda$ do not contribute extraneous boundary
terms to (A.43).  Next, integrate by parts on the third term from the
left on the top line in (A.43) to replace the latter by
\begin{equation}
\int_{S(s)} 3u^{-2} \sum_j |\underline{b}_{j,n}|^2 + \int_{\p S(s)}
3((t-t')^2+u^2)^{-\frac 12}\sum_j |\underline{b}_{j,n}|^2.
\end{equation}
Thus, an application of the triangle inequality shows that the
top line of (A.43) is no greater than
\begin{equation}
\zeta^{-1}\dt \int_{S(s)} \sum_j (|d\underline{\eta}_{j,n}|^2
+u^{-2}|\underline{b}_{j,n}|^2 ).
\end{equation}
Here, $\zeta$ is an $n$--independent constant.

Meanwhile, for all but a measure zero set of values for $s$, the left
most term on the bottom line of (A.43) is no greater than
\begin{equation}
\zeta\dt s\dt\int_{\p S(s)} \sum_j
(|d\underline{\eta}_{j,n}|^2
+u^{-2}|\underline{b}_{j,n}|^2 ).
\end{equation}
Indeed, this follows directly by writing $\underline{b}_{j,n}
= u
\cdot (u^{-1} \underline{b}_{j,n})$ and then completing the square.
To analyze the right most term on the second line in (A.43), introduce
$\theta_n \equiv (\sum_j |\underline{b}_{j,n}|^2)^{1/2}$.  As
$|d\theta_n| \le \sum_j |d\underline{\eta}_{j,n}|$, it follows from
the estimate for $\{\sigma_{j,n}^k\}$ and from Lemma 5.4.1 in \cite{Mo}
that this term is no greater than $\zeta \cdot 2^{-n/\zeta}
\int_{S(s)} (|d\underline{\eta}_{j,n}|^2 +
u^{-2}|\underline{b}_{j,n}|^2)$.  Here, one should apply Lemma 5.4.1
of Morrey using $G$ to denote the half disk $S'$ which is concentric
with $S(s)$ but has twice the radius.  Also, this lemma from Morrey is
stated for a function $u$; and to apply this lemma, take $u$ to be a
compactly supported function on $S'$ which equals the given $\theta_n$
on $S(s)$ and whose $L_1^2$ norm is no greater than $\zeta \int_{S(s)}
(|d\theta_n|^2 + u^{-2}|\theta_n|^2)$ for some $s$ and $n$ independent
constant $\zeta$.  (Such a function can be obtained by reflecting
$\theta_n$ through the boundary of $S(s)$ and then using a radial
cut-off function.)

With all of the preceding understood, then (A.43) is seen to imply (for
large $n$) the following inequality for almost all values of $s$:
\begin{equation}
\int_{S(s)}\sum_j(|d\underline{\eta}_{j,n}|^2 +u^{-2}
|\underline{b}_{j,n}|^2 ) \leq
\zeta s \int_{\p S(s)}\sum_j
(|d\underline{\eta}_{j,n}|^2 +u^{-2}|\underline{b}_{j,n}|^2 ).
\end{equation}
Here, $\zeta$ is independent of $s$ and of $n$.  This last
equation has the form $F \le \zeta s F'$ for a function $F$ on
$[0,s_0]$ which is almost everywhere differentiable; and such an equation
can be
integrated to find that $F(s) \le \zeta \cdot s^{1/\zeta}F(s_0)$.
Finally, given the first point of the lemma, this last expression is
(A.42).

End the digression.

With the estimate in (A.42) understood, turn now to the proof of the
third point of Lemma A.12.  For this purpose, there are two cases to
consider.  In the first case, the center point $w_1 = t_1 + i \cdot
u_1$ for a disk satisfies $u_1 \le s/4$.  Here, the required estimate
for the fourth point of Lemma A.13 follows directly from (A.42).

The second case to consider has $u_1 \ge 4s$.  Here, save for one
major point and a few minor ones, the argument follows that part of
the proof of Lemma 8.2 which begins in the paragraph after (8.17).
(These arguments from Section~8 can be borrowed when $n$ is large and
so both $|\sigma_{j,n}|$ and $|\sigma_{j,n}^k|$ from (A.39) are small
in the sense indicated in the statement of Lemma A.12.)  Here are the
necessary modifications:  First, replace $(\varphi_j,\nu_j)$ with
$(\underline{\varphi}_{j,n},\underline{b}_{j,n})$.  Second, the $\zeta
\cdot s_1^2$ term in (8.18) should be replaced by
\begin{equation}
\zeta\dt(u^{-2}_1\int_D \sum_j |\underline{b}_{j,n}|^2 +
\int_D\sum_{j,k} |\s^k_{j,n}|^2 \sum_m |\underline{b}_m|^2).
\end{equation}
The claim now is that right most term in (A.48) is no greater
than
$\zeta s_1^{1/\zeta}$ for some constant $\zeta \ge 1$ which is
independent of $u_1$, $s_1$ and $n$.  Accept this claim for the moment
and then the $\zeta s_1^2$ term in (8.18) can ultimately be replaced
by $\zeta s_1^{1/\zeta}$.  With this last replacement understood,
then (8.21) will be replaced by
\begin{equation}
f\leq \zeta s_1 f'+ \zeta s_1^{1/\zeta} .
\end{equation}
This last equation can be integrated from $s$ to $u_1/4$ to
find that
\begin{equation}
f(s)\leq\zeta ((s/u_1)^{1/\zeta}f(u_1/4)+s^{1/\zeta}|\ln
(s/u_1)|+1). 
\end{equation}
Since (A.42) already bounds $f(u_1/4)$ by $\zeta'
u_1^{1/\zeta'}$, for constant $\zeta' \ge 1$ which is independent of
$n$ and $u_1$, this last equation implies Lemma A.12's third point.

Thus, it remains only to establish the bound of (A.48) by $\zeta
s_1^{1/\zeta}$ for some constant $\zeta$ which is independent of
$u_1$, $s_1$ and $n$.  Here is the argument for such a bound on the
left most term in (A.48):  To begin, let $\theta \equiv (\sum_j
|\underline{b}_{j,n}|^2)^{1/2}$.  As $\theta$ vanishes on the $u = 0$
line, this function is no greater than $\int_{0 \le u' \le u}
|d\theta|du'$ at a point $(t,u)$.  Of course, this last expression is
no greater than $u^{1/2} (\int_{0 \le u' \le u}
|d\theta|^2du')^{1/2}$.  And, since $|d\theta| \le (\sum_j
|d\underline{b}_{j,n}|^2)^{1/2}$, this last expression is no greater
than $u^{1/2} (\int_{0 \le u' \le u} \sum_j
|d\underline{b}_{j,n}|^2)^{1/2}$.  Now, substitute the latter into the
left most term (A.48) to bound said term by
\begin{equation}
\zeta u_1^{-1} s_1\int_{S(2u_1)} \sum_j
|d\underline{b}_{j,n}|^2 . 
\end{equation}
(Remember that the disk $D$ is contained in the half disk
$S(2u_1)$.)  Finally, as $s_1 \le u_1/4$, one can use (A.42) to bound
this last expression by the desired $\zeta s_1^{1/\zeta}$.

Now consider the right most term in (A.48).  Here, Lemma 5.4.1 of
Morrey can be invoked once again given the estimates for $\sum_{j,k}
|\sigma_{j,n}^k|^2$ in the fifth point of Lemma A.12.  In this case,
take the domain $G$ in the statement of Morrey's lemma to be the disk
which is concentric to $D$ but has radius $u_1/2$.  Moreover, take
the function $u$ in the statement of Morrey's lemma to equal the
product of $\theta$ with a standard cut-off function; here, the
cut-off function should equal one on the disk which is concentric to
$D$ but has radius $u_1/4$ and it should vanish near the boundary of
the radius $u_1/2$ disk.  Furthermore, its derivatives should be
bounded by $\zeta/u_1$.  Finally, note that the $L_1^2$ norm of this
function $u$ over $G$ here can then be bounded apriori using (A.42).

\medskip
\noindent
{\bf Step 5}\qua The purpose of this step is to offer a proof of Lemma
A.9:

\medskip
\noindent
{\bf Proof of Lemma A.9}\qua  As remarked previously, the number
$\lambda_n$ which appears in (A.29) should be given as in (A.38).

To control the convergence of the sequence
$\{\underline{\eta}_{j,n}\}$ as maps into $\mbox{Sym}^p {\Bbb C}$, it
proves useful to digress for a moment to consider the standard
parametrization of the latter as ${\Bbb C}^{\,p}$.  Here, $a =
(a_1,\dots,a_p) \in {\Bbb C}^{\,p}$ determines the $p$'th order, monic
polynomial
\begin{equation}
{\cal P}_a(\phi)=\phi^p+a_1 \phi^{p-1} +\dots +a_p \end{equation}
and the roots of the latter give the corresponding point in
$\mbox{Sym}^p {\Bbb C}$.  Conversely, a point $\gamma =
\{\gamma_j\}_{1 \le j \le p} \in \mbox{Sym}^p {\Bbb C}$, determines
the polynomial ${\cal P}_a$, where the coefficients are given by
\begin{equation}
a_m=\sum_{K\in\Theta(m)} \prod_{j\in K} \g_j. \end{equation}
Here the sum runs over the set $\Theta(m)$ of $m$--element
subsets of
$\{1,\dots,p\}$.

Now consider the case where the coefficients in (A.52) are functions on
$D(0)$ given by setting $\gamma_j \equiv \underline{\eta}_{j,m}$ in
(A.53).  Denote these functions by $\{a_{m,n}\}_{1 \le m \le p}$.  (In
this regard, note that $a_{1,n} \equiv 0$.)  The fourth point of Lemma
A.12 implies that for each $m \in \{1,\dots,p\}$, the resulting
sequence $\{a_{m,n}\}$ is bounded as $n$ tends to $\infty$ and
satisfies a uniform bound on its Holder norm over $D(2)$.  Moreover,
because of the first point of Lemma A.12, this same sequence enjoys a
uniform bound on its Sobolev $L_1^2$ norm.

The existence of the uniform $L_1^2$ bound implies that any infinite
subsequence of the sequence $\{a_{m,n}\}$ has, itself, an infinite
subsequence which converges weakly in $L_1^2$ over the half disk
$D(0)$ to an $L_1^2$ function, $a_{m,\infty}$.  By taking a so-called
diagonal subsequence, one can arrange that the same subsequence
converges for each of the possible values for $m$.

Note that the existence of the uniform Holder bound over $D(2)$ for
each sequence $\{a_{m,n}\}$ implies that these subsequences converge
strongly in the $C^0$ topology on $D(2)$ to give Holder continuous
functions on $D(2)$.  With the preceding understood, it follows that
the limiting data $\{a_{m,\infty}\}_{1 \le m \le p}$ defines a Holder
continuous and Sobolev class $L_1^2$ map from $D(2)$ into
$\mbox{Sym}^p {\Bbb C}$.

Given the preceding, the proof of Lemma A.9 will be complete with a
demonstration that there exists $m \in \{1,\dots,p\}$ and an infinite
subsequence of $\{a_{m,n}\}$ whose limit, $a_{m,\infty}$, is not
identically zero.  This demonstration requires the following crucial
lemma:

\medskip
\noindent
{\bf Lemma A.13}\qua  {\sl{There exists a constant $c$ and an infinite
subset $\Omega \subset \{1,\dots\}$ with the following property:  If
$n \in \Omega$, then $0 < \lambda_n \le c\lambda_{n+4}$}.}

\medskip
The proof of this lemma is provided in Steps 8--12.  Assume its
validity for now.

Given Lemma A.13, here is how to complete the proof of Lemma A.9:  To
begin, take $a_{m,\infty}$ to be a limit of the sequence
$\{a_{m,n}\}_{n \in \Omega}$.  Now, suppose that $a_{m,\infty} \equiv
0$ for each $m$.  This assumption generates a contradiction.  Indeed,
were all $a_{m,\infty}$ identically zero, then the sequence $\{\sum_j
|\underline{\eta}_{j,n}|^2\}_{n=1,2,\dots}$ would have a infinite
subsequence which converged to zero uniformly over $D(2)$.  However,
if for some large $n$, $\sum_j |\underline{\eta}_{j,n}|^2$ is small on
$D(2)$, say less than $\e$ everywhere, then, because of the
second point of Lemma A.12,
\begin{equation}
\int_{D(4)} \sum_j (|d\underline{\eta}_{j,n}|^2 +u^{-2}|
{\mbox{im}}(\underline{\eta}_{j,n})|^2 \leq\zeta\dt\e  ,
\end{equation}
where $\zeta$ is independent of both $n$ and $\e$.  But, if
(A.54) is true, then Lemma A.13 forces $\int_{D(0)} \sum_j
(|d\underline{\eta}_{j,n}|^2 +
u^{-2}|\mbox{im}(\underline{\eta}_{j,n})|^2) \le c \cdot \zeta \cdot
\e$.  And, this last condition contradicts the first point of
Lemma A.12 when $\e$ is small.

\medskip
\noindent
{\bf Step 6}\qua  The purpose of this step and Step 7 is to
present the proof of Lemma A.10.

\medskip
\noindent
{\bf Proof of Lemma A.10}\qua  Define a map, $q\co  D(2) \rightarrow
\{1,\dots,p\}$ by assigning to $w \in D(2)$ the number of distinct
roots of ${\cal P}^w(\lambda) \equiv \phi^p +
a_{1,\infty}(w)\phi^{p-1} + \dots + a_{p,\infty}(w)$.  Say that $w$ is
a {\sl regular point} of $D(2)$ when $q(w)$ maximizes the function
$q(\cdot)$ on some neighborhood of $w$.  Since $q$ maps into a finite
set, it follows that the set of regular values in $D(2)$ is both open
and dense.  Note that on some neighborhood of a regular value, the
roots, $\{\underline{\eta}_{j,\infty}\}$, of the polynomial ${\cal
P}^{(\cdot)}$ can be viewed as a set of $p$ honest, complex valued
functions.  In fact, it follows from Lemma A.9 that each such
$\underline{\eta}_{j,\infty}$ is Holder continuous and Sobolev class
$L_1^2$ near a regular value point.

Now, let $w \in D(2)$ be a regular point.  The claim here is that each
$\{\underline{\eta}_{j,\infty}(w)\}$ obeys the equation
\begin{equation}
\bar\p\underline{\eta}_{j,\infty} -\textstyle{\frac 32}
u^{-1}{\mbox{im}}(\underline{\eta}_{j,\infty}) =0.
\end{equation}
on some neighborhood of a regular point in $D(2)$.  This claim will be
proved in Step 7, below, so accept it for now.

The proof of Lemma A.10 is completed as follows:  Let $w \in D(2)$ be a
regular point.  Then, on some open neighborhood $U$ of $w$ which
consists solely of regular points, the function $\tau$ is the product
of those elements in $\{\underline{\eta}_{j,\infty}\}$ which are not
identically zero on $U$.  This identifies $\tau$ on a neighborhood of
$w$ with one of the function $a_{m,\infty}$.  Indeed, $\tau|_U =
a_{m,\infty}$ where $m$ is the smallest of the integers $m' \in
\{2,\dots,p\}$ for which $a_{m',\infty}$ is not identically zero on
$U$.  Note that $\tau$ is continuous since it extends continuously as
zero to the complement of the set of regular points.  Furthermore, one
can prove that $\tau$ is Sobolev class $L_1^2$ with
\begin{equation}
\int_{D(2)} |d\tau|^2\leq\zeta\int_{D(2)}\sum_{1\leq m\leq p}
|da_{m,\infty}|^2 .
\end{equation}
(To prove (A.56), it is sufficient to exhibit a sequence of $L_1^2$
functions $\{\tau_k\}_{k=1,2,\dots}$ which converges to $\tau$ in
$L^2$ and such that each $\tau_k$ has $L_1^2$ norm bounded by the
right side of (A.56) for some $k$ independent constant $\zeta$.  To
construct such a sequence, let $\chi$ denote the standard bump
function on $[0,\infty)$ and set $\tau_k \equiv
\chi(|\mbox{ln}|\tau|/k)\cdot \tau$.  Thus, $|\tau_k| \le |\tau|$ and
when $k$ is large, then $\tau_k =\tau$ where $|\tau| \ge e^{-k}$ and
$\tau_k = 0$ where $|\tau| \le e^{-2k}$.  Thus, $\{\tau_k\}$ converges
to $\tau$ in the $L^2$ topology.  Furthermore, $|d\tau_k| = 0$ where
$|\tau| \le e^{-2k}$, and $|d\tau_k| \le |d\tau|(1 + \zeta/k)$
otherwise.  This last inequality implies the $\tau_k$ version of
(A.55) since $|d\tau| \le \sum_m |da_{m,\infty}|$ where $|\tau| >
0$.)  Next, remark that (A.55) implies that $\tau$ obeys (A.30) on the
open set of regular points with the proviso that the integer $p'$ and
the function ${\cal P}$, may not extend continuously from this set to
the whole of $D(2)$.  In any event, define a function $f$ on $D(2)$
by setting $f = 0$ at non regular points, and by setting $f = \frac
{3i}{4u} p'(1-{\cal P})$ on the set of regular points.  Thus, $f$ is
bounded, and measurable on any subset of $D(2)$ where $u$ is bounded
away from zero.

Now, let $D \subset D(2)$ be a disk on which $u$ is bounded away from
zero, and let $\kappa$ be a continuous and Sobolev class $L_1^2$
solution to the equation ${\bar \partial} \kappa = -f$ on $D$.
Since $f$ is bounded and measurable on $D$, a standard construction
with the kernel $\pi^{-1}(w-w')^{-1}$ for the operator ${\bar
\partial}^{-1}$ will give a continuous and $L_1^2$ function on $D$
which solves this last equation.

With $\kappa$ understood, observe that $\tau = e^{\kappa}\tau_0$ on
$D$, where $\tau_0$ is a continuous and Sobolev class $L_1^2$
function on $D$ which obeys the equation ${\bar \partial} \tau_0 = 0$
at all points where $\tau_0 \ne 0$.  Since $\tau_0$ is Sobolev class
$L_1^2$, this last condition implies that $\tau_0$ is holomorphic
everywhere on $D$.  Thus, $\tau_0$ vanishes at no more than a finite
set of points on any compact domain in $D$.  And, since the zeros of
$\tau_0$ are the same as those of $\tau$ in $D$, it follows that all
but finitely many points of $D$ are regular points.  In particular,
this implies that the integer $p'$ is in (A.30) is constant over $D(2)$
and that the function ${\cal P}$ is continuous on the complement of a
set which is at worst countable, and which has no accumulation points
where $u > 0$.

\medskip
\noindent
{\bf Step 7}\qua  This step completes the proof of Lemma
A.10 with a proof of (A.56).  To begin the argument, let $w$ denote the
regular point in question.  By definition, the polynomial ${\cal P}^w$
has $q
\equiv q(w)$ distinct roots, and it has this same number at all points
$w'$ which are close to $w$.  Each root of ${\cal P}^w$ also has a
multiplicity and this multiplicity is also constant on some neighborhood of
$w$. With the preceding understood, let $\underline{\eta}_{\infty}$ be a root
of $\cal P^{(\cdot)}$ thought of as a function,k on a neighborhood, $D$, of
$w$. Let $e$ denote the multiplicity of this root.
  Then, for large $n$ which labels Lemma
A.9's convergent subsequence, there will be precisely $e$ elements of
$\{\underline{\eta}_{j,n}\}$ whose values at points in $D$ are very
close to the corresponding values of $\underline{\eta}_{\infty}$.
Furthermore, all other elements in $\{\underline{\eta}_{j,n}\}$ will
have values at the points of $D$ which differ substantially from the
corresponding values of $\underline{\eta}_{\infty'}$.  In particular,
the $e$ elements of $\{\underline{\eta}_{j,n}\}$ which are close to
$\underline{\eta}$ form an unambiguous subset
$\{\underline{\eta}_{j,n}\}$ over $D$ which one can declare,
unambiguously, consists of the first $e$ elements of
$\{\underline{\eta}_{j,n}\}$.  Then, $\gamma_n \equiv \sum_{1 \le j
\le e} \underline{\eta}_{j,n}$ defines an unambiguous, complex valued
function on $D$.

It then follows from (A.39) that the function $\gamma_n$ obeys the
equation
\begin{equation}
\bar\p\g_n -\textstyle{\frac 32}
u^{-1}{\mbox{im}}(\g_n)+\s_{j,n}\dt\p\g_n +
\sum_{1\leq j\leq e}\sum_k
\s^k_{j,n}\dt{\mbox{im}}(\underline{\eta}_{k,n})=0. 
\end{equation}
It then follows from the fourth and fifth assertions of Lemma A.12 that
the last two terms on the right side of (A.57) converge to zero (in the
$L^2$ topology) as $n$ tends to infinity, and thus (A.55) holds.

\medskip
\noindent
{\bf Step 8}\qua This step begins the proof of Lemma A.13.

\medskip
\noindent
{\bf Proof of Lemma A.13}\qua  The fact that $\lambda_n > 0$ follows from
the assumption that the $\{\eta_j\}$ are not identical over the region
$A$.  Meanwhile, the argument used below to find the set $\Omega$ and
the constant $c$, is a modified version of the argument used by
Aronszajn in the proof of his unique continuation theorem in \cite{Ar}.

The argument given below is a proof by contradiction, so for this
purpose. assume at the outset that there is no constant $c$ which makes
Lemma A.13 true.  This means that for any $c > 0$, the numbers
$\{\lambda_n\}$, though not identically zero, obey the inequality
$\lambda_n \ge c \cdot \lambda_{n+4}$ when $n$ is sufficiently large.
In particular, this last inequality implies that for any $c > 0$ and
sufficiently large $n$,
\begin{equation}
0 < \l_n\leq c^{-n} .
\end{equation}
A contradiction will be shown to follow from this assumption.

\medskip
\noindent
{\bf Step 9}\qua The contradiction to the assumption in
(A.58) is derived with the help of (A.2) which is obeyed by each $\eta \in
\{\eta_j\}$. The first remark with regard to (A.2) is that the term with
$\alpha_1$ in (A.2) complicates the arguments that follow.  The approach
taken here for dealing with this term involves the replacement of the
coordinate $w = t + i \cdot u$ by a new complex coordinate, ${\varpi}$.  This
particular step serves as a digression to introduce this new coordinate.

In particular, the coordinate $\varpi$ is chosen to have the
following properties:

\begin{itemize}
\item $|\varpi-w|\leq\zeta\dt |w|^2$
where $\zeta$ is independent of $w$ and $\d$.
\item $|\p\varpi-1|\leq\zeta\dt |w|$
where $\zeta$ is independent of $w$ and $\d$ when $\d$ is small.
\item Im$(\varpi)=0$ when $u=0$.
\item $\bar\p\varpi+\a_1(\underline{\eta})\dt\p\varpi=0$
where $\a_1$ appears in (A.2) and where\newline
\hbox{}\qquad $\underline{\eta}\equiv
p^{-1}\dt\sum_j\eta_j$. \hspace\fill (A.59)
\end{itemize}
\addtocounter{equation}{1}

The point of introducing $\varpi$ is as follows:  When $\eta$
in (A.2) is written in terms of $\varpi$, then the latter
equation is equivalent to
\begin{eqnarray}
&& \bar\p_{\varpi}\eta-\textstyle{\frac 32}{\mbox{Im}}(\varpi)^{-1}\dt 
(1+\a_3)\dt b
\nonumber \\
&&\quad +(1-\a^2_1)^{-1}\dt(\overline{\p\varpi})^{-1}\dt
(\a_1(\eta)-\a_1(\underline{\eta}))\dt (\p_w\eta)=0 , 
\end{eqnarray}
where $\alpha_3 = -1 + (1 + \alpha_2) \cdot (1 - \alpha_1^2)^{-1}
\cdot (\overline{\p\varpi})^{-1} \cdot \mbox{im}(\omega)/b$ with
$\alpha_1$ and $\alpha_2$ as in (A.2).  (Here, $w$, $b$,
$\partial_w\eta$ and $\partial\varpi$ are considered as functions of
$\varpi$ by inverting the map $w \rightarrow \varpi(w)$.)  Note that
$\alpha_3$ obeys
\begin{equation}
|\a_3|\leq c'\dt |\varpi|.
\end{equation}
for some constant $c'$ whichis independent of $\varpi$.

Here is an existence proof for the complex function $\varpi$:  To
begin, reintroduce the bump function $\chi$ and then, with $\delta >
0$ chosen to define the set $A$ as in Step~1 of the first section of this
Appendix, set $\chi_1$ to denote the function on ${\Bbb C}$ given by
$\chi(2|w|/\delta)$.  Next, introduce $\alpha \equiv
\alpha_1(\underline{\eta})$.  Then, consider writing $\varpi = w \cdot
e^{\sigma}$, where $\sigma$ is constrained to obey the equation
\begin{equation}
\bar\p\s +\chi_1\dt w^{-1}\dt (\a\dt\p\s+\a)=0. \end{equation}
(Note that $|\alpha| \le \zeta \cdot |w|^2$ by virtue of (8.4) and
(A.1).  Thus, the coefficients in (A.62) are bounded.)

To be more precise, use the Greens kernel for the operator ${\bar
\partial}$ to find a solution to (A.62) as a fixed point of the mapping
$T$ which sends a complex valued function $\sigma$ to
\begin{eqnarray}
&&\pi^{-1}\dt w\dt\int_{u\geq 0}
d^2\l [(\chi_1\dt (\a\dt\p\s +\a)\dt \l^{-2} (\l-w)^{-1} \nonumber \\
&& +(\chi_1\dt
(\a\bar\p\s +\a)\dt\bar\l^{-2} (\bar\l-w)^{-1})].
\end{eqnarray}
Note that ${\bar \partial}T(\sigma) = -\chi_1 \cdot w^{-1} \cdot
(\alpha \cdot \partial\sigma + \alpha)$ where $u > 0$ so that $\varpi = w
\cdot e^{\sigma}$ satisfies the fourth point in (A.59) when $\sigma$ is
a fixed point of $T$.  Also, $T(\sigma)$ is real where $u = 0$ so that
$\varpi$ will also satisfy the third point of (A.59).  Thus, it remains
only to find a fixed point $\sigma$ to (A.63) for which the
corresponding function $\varpi$ obeys the first and second points of
(A.59).

The contraction mapping theorem will be used to find a fixed point of
$T$ with the required properties.  In this regard, introduce the
Banach space ${\cal H}$ which is defined by completing the compactly
supported functions on ${\Bbb C}$ with the norm $\|\cdot\|_*$ whose
square sends $\sigma$ to
\begin{equation}
\sup_{\Bbb C} |\s|^2+\sup_{\Bbb C} \d^2|d\s|^2+ \sup_{0 < r < 2d} \ 
\sup_{w\in\Bbb C} (\d/r)^c\dt \d^2\dt \int_{|\l-w|\leq r} |\nabla
d\s|^2 \ d^2\l . \end{equation}
Here, $c \in (0,1)$ is half the value of the exponent $1/\zeta$ which
appears in (8.12).

Since the map $T$ is linear inhomogeneous, the contraction mapping
theorem will find a unique fixed point for $T$ on ${\cal H}$ given that
\begin{eqnarray}
T(0) & \in &\cal H.\nonumber\\
\|T(\s)-T(0)\|_* & < & 2^{-1}\dt\|\s\|_*.
\end{eqnarray}
Moreover, to prove the first point in (A.65) it is sufficient to verify
that $\|T(0)\|_*$ is finite; while the second point need be tested
only with smooth functions having compact support.

Given that (A.65) holds then the contraction mapping theorem finds a
fixed point $\sigma$ to $T$ with $\|\sigma\|_* \le 2 \cdot
\|T(0)\|_*$.  Furthermore, if
\begin{equation}
\|T(0)\|_* \leq \zeta\dt\d  ,
\end{equation}
then it follows from the contraction mapping theorem that the fixed
point $\sigma$ obeys $\sigma(0) = 0$ and $|d\sigma| \le \zeta$.  These
last conditions imply that $|\sigma| \le \zeta \cdot |w|$ and thus the
first two points in (A.59) will also be satisfied.

Both of the points in (A.65) and (A.66) can be addressed with the help
of the following lemma:

\medskip
\noindent
{\bf Lemma A.14}\qua  {\sl {Let $\tau$ be a continuous, Sobolev class
$L_1^2$, complex valued function on ${\Bbb C}$ which vanishes where
$u\le 0$ and which has compact support on the disk of radius $\delta$
about the origin.  Suppose that}
\begin{equation}
E^2\equiv \sup_{\Bbb C} |\tau|^2+
\sup_{0 < r \leq \d} \ \sup_{w\in\Bbb C} (\d/r)^{2c}
\int_{|\l-w|\leq r} |d\tau|^2
\end{equation}
{is finite.  Then,}
\begin{equation}
Q\equiv (4\pi)^{-1}\dt w\dt\int_{u\geq 0}
[\tau\dt\l^{-1}(\l -w)^{-1}+\bar\tau\dt\bar\l^{-1}
(\bar\l-w)^{-1}]d^2\l
\end{equation}
{is in ${\cal H}$ and satisfies $\|Q\|_* \le \delta \cdot E$}.}

\medskip
Note that Lemma A.14 with $\tau = \chi_1 \cdot w^{-1}\alpha$ gives the
first point in (A.68) for in this case, $E$ is bounded by $\zeta \cdot
\delta^c$ by virtue of (A.1), (8.4) and (8.12).  Indeed, as $\alpha =
\zeta \cdot b^2 + {\cal O}(b^4)$, the latter imply that $|\alpha| \le
\zeta \cdot |w|^2$ and $|d\alpha| \le \zeta \cdot |w| \cdot \sum_j
|d\nu_j|$.  Thus, $|\tau|$ and $|d\tau|$ have support where $|w| \le
\delta$, while the former is bounded by $\zeta \cdot |w|$ and the
latter is bounded by $\zeta \cdot (1 + \sum_j |d\nu_j|)$.  In
particular, these last remarks with Lemma A.14 imply that $\|T(0)\|_*
\le \zeta \cdot \delta^{1+c}$ for some $\delta$ independent constant
$\zeta$.  This gives (A.66) for small $\delta$.

Meanwhile, Lemma A.14 with $\tau = \chi_1 \cdot w^{-1} \cdot \alpha
\cdot d\sigma$ gives the second point in (A.65) (for small $\delta$).
In this case, $|\tau| \le \zeta \cdot \delta \cdot |d\sigma|$ while
$|d\tau| \le \zeta \cdot (|d\sigma| + \delta \cdot |\nabla d\sigma|)$
and both have support in the radius $\delta$ disk on ${\Bbb C}$ with
center the origin.  Thus, $E \le \zeta \cdot \delta \cdot
\|\sigma\|_*$ in this case.

Thus, the verification of (A.59) requires only the proof of Lemma A.14.

\medskip
\noindent
{\bf Proof of Lemma A.14}\qua  To begin, use the definition of $Q$ as
an integral to find the pointwise bound $|Q| \le \zeta \cdot \delta
\cdot \sup_{\Bbb C} |\tau|$.  Note also that where $|w| > 4\delta$,
the definition of $Q$ gives $|dQ| + \delta \cdot |\nabla dQ| \le \zeta
\cdot (\sup_{\Bbb C} |\tau|) \cdot \delta^2/|w|^2$.

To estimate the last two terms in (A.64) where $|w| \le 4 \cdot
\delta$, it proves useful to employ the equation ${\bar \partial}Q =
-\tau$ which is valid where $u > 0$.  Where $u \le 0$, ${\bar
\partial}Q$ is equal to the function $\tau'$ who value at $t + iu$
is that of $-{\bar \tau}$ at the conjugate point $t - i \cdot u$.
Differentiating these equations finds that ${\bar \partial}(dQ) =
-d\tau$ where $u > 0$ and ${\bar \partial}(dQ) = -d\tau'$ where $u <
0$.  Note that there is no ``delta function'' mass for ${\bar
\partial}(dQ)$ along $u = 0$ because $\tau$ is zero there by
assumption.  With this last equation for ${\bar \partial}(dQ)$
understood, a bound on $\int_{|\lambda-w| \le r} |\nabla
dQ|^2d^2\lambda$ by $\zeta \cdot (\delta/r)^{2c} \cdot
\int_{|\lambda-w| \le r} |d\tau|^2d^2\lambda$ follows using Morrey's
Theorem 5.4.1 from \cite{Mo}.  Note that this same bound plus the
aforementioned bound on $|dQ|$ at points where $|w| \ge 4\delta$
implies the bound $|\nabla Q| \le \zeta \cdot E$ at all points.
Indeed, in this regard, one need only appeal to Theorem 3.5.2 in
\cite{Mo}.  (Note that the bound for $|dQ|$ can be obtained directly from
the integral equation for $Q$.)

\medskip
\noindent
{\bf Step 10}\qua  With the function $\underline{\eta} =
p^{-1} \sum_j
\eta_j$ understood, now re-introduce the difference functions
$\underline{\eta}_j = \eta_j - \underline{\eta}$.

Remark first that (A.60) implies that the $\{\underline{\eta}_j\}$ obey
an equation which has the schem\-atic form
\begin{equation}
\bar\p_{\varpi}\underline{\eta}_j-\textstyle{\frac 32}
{\mbox{im}}(\varpi)^{-1}
{\mbox{im}}(\underline{\eta}_j) +\sum_k R_{jk} \underline{\eta}_k=0  ,
\end{equation}
where $|R_{jk}| \le \zeta \cdot \sum_j (|b_j| + |b_j| \cdot |d\eta_j|)
\le \zeta \mbox{ im}(\varpi)(1 + \sum_j|d\eta_j|^2)$.  Here, the
$1$--form norm can be taken using the Euclidean metric from the
$\varpi$--plane.  (This is equivalent to the Euclidean metric from the
$w$--plane.)

To make use of this last equation, introduce the polar coordinates
$(r,\theta)$ on the $\varpi$--plane by writing $\varpi =
re^{i\theta}$.  Thus, $2{\bar \partial}_{\varpi} =
e^{i\theta}(\partial_r + i \cdot r^{-1} \cdot \partial_{\theta})$.  In
particular, (A.69) is equivalent to
\begin{equation}
\p_r\underline{\eta}_j+ ir^{-1}\p_{\theta}
\underline{\eta}_j+i\textstyle{\frac 32}
r^{-1}(\cot\theta-i)(\underline{\eta}_j-\bar{\underline{\eta}}_j)
=-2e^{-i\theta}\sum_k R_{jk} \underline{\eta}_k . \end{equation}
Now, write $\underline{\eta}_j \equiv r^{-3/2} \sin^{-3/2}\theta \
\alpha_j$ in which case (A.70) is equivalent to
\begin{equation}
\p_r\s_j +ir^{-1}\p_{\theta}\s_j-i\textstyle{\frac 32}
r^{-1}(\cot\theta -i)\bar\s_j =-2e^{-i\theta} \sum_k R_{jk}\s_k . 
\end{equation}
To proceed with (A.71), fix $\e \in (0,\delta/100)$ and a
positive integer $m$ and then introduce the bump function
$\chi_{\e,m}$ whose value at $\omega$ is $\chi(r/\e)[1 -
\chi(2^{m+1}r/\delta)]$.  This function equals $1$ where $2^{-m}\delta
\le r \le \e$ and it vanishes where $r \ge 2 \cdot \e$ and
where $r \le 2^{-m-1} \delta$.

Next, introduce $f_j \equiv \chi_{\e,m} \varpi^{-N}\sigma_j$,
and then multiply both sides of (A.69) by $\chi_{\e,m}
\varpi^{-N}$ to obtain the following equation:
\begin{eqnarray}
&&\p_r f_j+ir^{-1}\p_{\theta}f_j-i\textstyle{\frac 32} r^{-1}(\cot\theta -i)
e^{2iN\theta}\bar f_j \nonumber \\
&& =-2e^{-i\theta} (\sum_k R_{jk}
f_k-(\bar\p\chi_{\e,m})\varpi^{-N}\s_j). \end{eqnarray}
This ``master equation'' for $f_j$ will be analyzed further below.
However, during the subsequent discussion, keep in mind that
\begin{equation}
|f_j| =\chi_{\e,m} r^{-N+3/2} \sin^{3/2}\theta  \ |\underline{\eta}_j| .
\end{equation}

\medskip
\noindent
{\bf Step 11}\qua  With (A.72) understood, square both sides,
sum over the index $j$ and integrate the result over $A$.  The left hand
side of the resulting equality can be manipulated with an integration by
parts to read:
\begin{equation}
\int_A\sum_j (|\p_rf_j|^2 +r^{-2}|\p_{\theta}f_j-i\textstyle{\frac 32}
r^{-1}(\cot\theta -i) e^{2iN\theta}\bar f_j|^2). \end{equation}

\noindent
(When driving (A.74), remember that $f_j$ has support only where
$2^{-m-1}\delta \le r \le 2\e$, and $f_j$ also vanishes at
$\theta = 0,\pi$.  As with the integration by parts in the proof of
Lemma 8.2, there are no anomalous boundary terms from the points in
$\pi(\Lambda)$ here.)

Because of the form of the right hand side of (A.72), the integral in
(A.74) is no greater than
\begin{equation}
8\int_A(\sum_{j,k} |R_{j,k}|^2) (\sum_i |f_i|^2)
+2\int_A |d\chi_{\e,m}|^2 r^{-2N} (\sum_j |\s_j|^2). \end{equation}
The task now is to manipulate the expressions in (A.74) and (A.75) to
obtain a useful inequality.  For this purpose, the task ahead is to
replace (A.74) by smaller terms and (A.75) by larger terms.

The first step in this strategy concerns the left most integral on the
right side of (A.74).  Here is the key lemma for treating this term:

\medskip
\noindent
{\bf Lemma A.15}\qua  {\sl {There exists $\zeta \ge 1$ which is independent
of $N$, and $m$ and such that when $\e < \zeta$, then the
left most integral in} (A.75) {is no greater than $1/1000$ of the
integral in} (A.74).}

\medskip
This lemma is proved in Step 14, below.

\medskip
\noindent
{\bf Step 12}\qua  To continue with the argument for Lemma
A.13, observe that Lemma A.15 has the following consequence:  Let $\gamma
\equiv (\sum_j |f_j|^2)^{1/2}$.  As long as $\e$ is less than
$\zeta^{-1}$ from Lemma A.5, then
\begin{equation}
\int_A |\p_r\g|^2 \leq\zeta\int_A |d\chi_{\e,m}|^2
r^{-2N}(\sum_j |\s_j|^2) .
\end{equation}
where $\zeta$ is independent of $m$ and $N$ and $\e$.  (This
follows because of the inequality $|\partial_r\gamma| \le (\sum_j
|\partial_rf_j|^2)^{1/2}$.)  Moreover, because $\gamma$ has support
where $r < 2\e$, the left hand side of (A.76) is no greater than
\begin{equation}
\zeta^{-1}\e^{-2}\int_A\g^2 =\zeta^{-1}\e^{-2}
\int_A r^{-2N}\sum_j |\s_j|^2  ,
\end{equation}
where $\zeta \ge 1$ is independent of $\e$, $N$ and $m$.

Meanwhile, the right hand side of (A.76) has an integrand whose support
(for large $m$) is located in two disjoint regions.  The first region
has $r \in [\e,2\e]$ and the contribution from this region
to the right side of (A.76) is no greater than
\begin{equation}
\zeta\e^{-2}\int_{A\cap\{r,\e\leq r\leq 2\e} r^{-2N}\sum_j
|\s_j|^2 .
\end{equation}
The other region that contributes to the right hand side of (A.76) has
$r$ restricted to $[2^{-m-2}\delta,2^{-m}\delta]$.  The claim now is
that the contribution to the right hand side of (A.76) from this second
region tends to zero as $m$ is taken to infinity with $N$ and
$\e$ fixed.  To see that such is the case, first note that the
contribution from the region in question is no greater than
\begin{equation}
\zeta\dt 2^{2(N+1)m}\d^{-2(N+1)}
\int_{A\cap\{r_* 2^{-m-2}\d\leq r\leq 2^{-m}\d\}} \sum_j
|\underline{\eta}_j|^2r\sin^3\theta \ dr \ d\theta . \end{equation}
Now let $\kappa \equiv (\sum_j |\underline{\eta}_j|^2)^{1/2}$.  Note
that $\kappa$ vanishes at $r = 0$.  Thus, the fundamental theroem of
calculus finds that for any constant $r'$, one has
\begin{equation}
\int_{0\leq\theta\leq\pi} |\kappa(r',\theta)|^2\sin^3\theta \ d\theta \leq
\int_{0\leq\theta\leq\pi}(\int_{0\leq r\leq r'}\p_r\kappa \ dr)^2 \sin^3\theta \
d\theta .
\end{equation}
And, Holder's inequality asserts that the right hand side in (A.80) is
no greater than
\begin{equation}
\int_{A\cap\{r\leq r'\}} |d\kappa|^2 r\sin^3\theta \ dr d\theta .
\end{equation}
In particular, since $|d\kappa| \le (\sum_j
|d\underline{\eta}_j|^2)^{1/2}$, these last two inequalities imply
that (A.79) is no greater than
\begin{equation}
\zeta\dt 2^{2(N+1)m}\d^{-2(N+1)}\l_{m-1}^2 \end{equation}
when $m$ is large.  This last expression tends to zero as $m$ tends to
infinity because of the assumption that (A.58) holds for any $c$ if $n$
is sufficiently large.

\medskip
\noindent
{\bf Step 13}\qua  With (A.79) seen to vanish as $m$ tends
to infinity, one obtains from (A.77) and (A.78) and Lemma A.15 the following
conclusion:  There exists $\zeta \ge 1$ such that when $\e \le
\zeta^{-1}$, then
\begin{eqnarray}
&&\int_{A\cap\{r\leq \e/2\}} r^{-2N}\sum_j
|\underline{\eta}_j|^2r\sin^3\theta \ dr \ d\theta \nonumber \\
&&\leq \zeta
\int_{A\cap\{\e\leq r\leq2\e\}} r^{-2N}\sum_j
|\underline{\eta}_j|^2r\sin^3\theta \ dr \ d\theta . \end{eqnarray}
This last inequality can be exploited to complete the proof of Lemma
A.13 as follows:  Multiply both sides of (A.83) by $(3\e/4)^{2N}$
to find that
\begin{eqnarray}
&&\int_{A\cap\{r\leq \e/2\}}\sum_j
|\eta_j|^2r\sin^3\theta \ dr \ d\theta \nonumber \\
&&\leq
\zeta\dt (15/16)^N\dt\int_{A\cap\{r\leq \e\}}\sum_j |\eta_j|^2r\sin^3\theta \
dr \ d\theta .
\end{eqnarray}
Since $N$ here can be arbitrarily large, it follows that all $\eta_j
\equiv 0$.  As this is assumed not to be the case, Lemma A.13 is proved.

\medskip
\noindent
{\bf Step 14}\qua  This step contains the
proof of Lemma A.15.

\medskip
\noindent
{\bf Proof of Lemma A.15}\qua To begin the argument, recall that
$|R_{j,k}| \le \zeta r \sin\theta(1 + \sum_j |d\eta_j|)$, and thus the
left most term in (A.75) is no greater than
\begin{equation}
\zeta\int_A r^2\sin^2\theta(\sum_j |f_j|^2)+ \zeta\int_A \sum_j
r^2\sin^2\theta(\sum_j |d\eta_j|^2)(\sum_j |f_j|^2). \end{equation}
To deal with the left most term in (A.85), reintroduce $\gamma \equiv
(\sum_j |f_j|^2)^{1/2}$.  As previously remarked,
$|\partial_r\gamma|^2 \le \sum_j |\partial_rf_j|^2$ and $\int_A
|\partial_r\gamma|^2 \ge \zeta^{-1} \e^{-2} \int_A \gamma^2$.
Thus,
\begin{equation}
\int_A\sum_j |f_j|^2\leq \zeta\e^2\int_A \sum_j |\p_r f_j|^2  ,
\end{equation}
where $\zeta$ is independent of $\e$, $m$ and $N$.  Thus, the
left most term in (A.85) is no greater than
\begin{equation}
\zeta\e^4\int_A\sum_j (|\p_r f_j|^2 +r^{-2}|\p_{\theta}f_j
-i\textstyle{\frac 32} r^{-1}(\cot\theta -i)e^{2iN\theta}\bar f_j|^2), 
\end{equation}
where $\zeta$ is independent of $\e$, $m$ and $N$.

The strategy for bounding the right most term in (A.85) calls on Lemma
8.2 and Lemma 5.4.1 in \cite{Mo}.  However, certain manipulations must be
done first.  In this regard, the first step is to produce some special
function partitions for $(0,\pi)$ and also for $(0,2\e)$.

With regard to the interval $(0,\pi)$, introduce, for each $k \in
\{0,1,\dots\}$, the angle $\theta_k = 2^{-k}(\pi/2)$.  Then, let
$\chi$ denote the standard bump function which is zero on $[2,\infty)$
and $1$ on $[0,1]$).  However, here, assume $\chi \equiv 0$ on
$[3/2,\infty)$.  Now, set $\chi_k^{\theta}(\theta) \equiv \chi(|\theta
- \theta_k|/\theta_{k+1})$.  Thus, $\chi_k^{\theta} = 0$ if $\theta <
\theta_k/4$ or if $\theta > 7\theta_k/4$.  On the other hand,
$\chi_k^{\theta} = 1$ if $\theta_k/2 \le \theta \le 3\theta_k/2$.
(Note that $\theta_k/2 = \theta_{k+1}$.)  With the preceding
understood, set $\theta_k$ for $k \in\{\dots,-2,-1\}$ to equal $\pi -
\theta_{-k}$, and define $\chi_k^{\theta}(\theta)$ by reflecting
$\chi_{-k}^{\theta}$ through the point $\pi/2$.  Thus defined, the set
$\{\chi^{\theta}_k\}_{k \in {\Bbb Z}}$ defines a partition of
$(0,\pi)$ with a uniform bound on the number of elements which are
non-zero at any one point, but with at least one element equal to one
at each point.

As for the interval $(0,2\e)$, introduce $r_k \equiv
2^{-k+1}\e$, and set $\chi_k^r(r) \equiv
\chi(|r-r_k|/r_{k+1})$.  The set $\{\chi_k^r\}_{k \ge 0}$ defines a
partition of $(0,2\e)$ with a uniform bound on the number of
elements which are non-zero at each point, but with at least one
element equal to one at each point.

With the preceding understood, it follows that
\begin{equation}
\sum_j |f_j|^2 \leq \sum_{k,k'} |\chi^{\theta}_k\chi^r_{k'}f_j|^2 .
\end{equation}
Thus, the second term in (A.85) is no greater than
\begin{equation}
\zeta\sum_{k,k'} r_{k'}^2\sin^2\theta_k\int_A (\sum_j
|d\eta_j|^2)\g_{k,k'}^2  ,
\end{equation}
where $\gamma_{k,k'} \equiv (\sum_j
|\chi_k^{\theta}\chi_{k'}^rf_j|^2)^{1/2}$.

Now, use Lemma 8.2 and Lemma 5.4.1 in \cite{Mo} on each $(k,k')$ term in
(A.89) to find a constant $\zeta \ge 1$ which is independent of
$k,k',N,m$ and $\e$ and is such that (A.89) is no greater than
\begin{equation}
\zeta\sum_{k,k'} r_{k'}^{2+1/\zeta} \sin^{2+1/\zeta}\theta_k
\int_A |d\g_{k,k'}|^2 .
\end{equation}
Then, differentiate $\gamma_{k,k'}$ to find that this last expression
is no greater than
\begin{eqnarray}
& &\zeta\sum_{k,k'} r_{k'}^{2+1/\zeta} \sin^{2+1/\zeta}\theta_k \Biggl\{
\int_A (r^{-2}_{k'}\theta_k^{-2}|{\chi'}^{\theta}_k|^2 \chi^r_k  \nonumber \\
& &+r_{k}^{-2}\chi^{\th}_k |{\chi'}^r_{k'}|^2) \sum_j |f_j|^2 +\int_A
\chi^{\th}_k \chi^r_{k'} \sum_j |df_j|^2\Biggr\}.
\end{eqnarray}
Here, $\chi'$ denotes the derivative of the function $\chi$.

Proceeding further, the triangle inequality applied to $|df_j|^2$
finds the latter no less than
\begin{eqnarray}
&&\zeta(|\p_r f_j|^2+r^{-2}|\p_{\theta} f_j-i\textstyle{\frac 32}
r^{-1}(\cot\theta-i)e^{2iN\theta}\bar f_j|^2 \nonumber \\
&&+r^2_{k'} \sin^{-2}\theta_k |f_j|^2).
\end{eqnarray}
Together, (A.91) and (A.92) imply that the second term in (A.85) is no
greater than
\begin{eqnarray}
&&\zeta(\sum_{k,k'} r_{k'}^{1/\zeta} \sin^{1/\zeta}\theta_k)
(\int_A \sum_j |\p_rf_j|^2 \nonumber \\
&&+r^{-2}|\p_{\theta}f_j-i\textstyle{\frac 32}
r^{-1}(\cot\theta-i)e^{2iN\theta}\bar f_j|^2 +\int_A\sum_j |f_j|^2).
\end{eqnarray}
(Here, as in the previous lines, $\zeta$ is independent of $\e$,
$m$ and $N$.)  Since the double sum here is no greater than
\begin{eqnarray}
\zeta\e^{1/\zeta}(\sum_{k=1,2,\dots} 2^{-k/\zeta})^2\leq \zeta'
\e^{1/\zeta'},
\end{eqnarray}
the second term in (A.85) can not be greater than
\begin{eqnarray}
&&\zeta\e^{1/\zeta}\int_A\sum_j
(|\p_rf_j|^2+r^{-2} |\p_{\theta}f_j- \nonumber \\
&&i\textstyle{\frac 32}
r^{-1}(\cot\theta-i)e^{2iN\theta}\bar f_j|^2
+\int_A\sum_j |f_j|^2).
\end{eqnarray}
Here, $\zeta$ is independent of $\e$, $m$ and $N$.

Together, (A.86), (A.87) and (A.95) imply the claim in Lemma A.15.

\end{document}